\documentclass[leqno]{siamonline1116}

\usepackage{amsmath}
\usepackage{amssymb}
\usepackage{lipsum}
\usepackage{braket,amsfonts}
\usepackage{graphicx}
\usepackage{epstopdf}
\usepackage{algorithmic}
\usepackage{multirow}
\usepackage{bold-extra}
\usepackage{fancybox,fancyhdr}
\usepackage{morefloats}

\ifpdf
  \DeclareGraphicsExtensions{.eps,.pdf,.png,.jpg}
\else
  \DeclareGraphicsExtensions{.eps}
\fi
\usepackage{setspace,color}
\usepackage[caption=false]{subfig}

\numberwithin{theorem}{section}
\numberwithin{equation}{section}

\usepackage[mathscr]{eucal}

\newtheorem{remark}[theorem]{Remark}

\DeclareMathAlphabet{\mathpzc}{OT1}{pzc}{m}{it}

\def\R{{\mathbb R}}
\def\C{{\mathbb C}}
\def\N{{\mathbb N}}

\def\Z{{\mathbb Z}}

\def\KK{{\mathbb K}}
\def\MM{{\mathbb M}}

\def\OO{{\mathbb O}}

\def\rank{\mathrm{rank}}

\def\rd{\mathrm{d}}
\def\Om{\Omega}

\def\f{\frac}
\def\p{\partial}

\def\na{\nabla}
\def\la{\langle}
\def\ra{\rangle}

\def\a{{\boldsymbol a}}

\def\bc{{\boldsymbol c}}

\def\e{{\boldsymbol e}}

\def\bsf{{\boldsymbol f}}

\def\bi{{\mathbf i}}

\def\bsj{{\boldsymbol j}}
\def\bk{{\boldsymbol k}}
\def\bsl{{\boldsymbol l}}

\def\bm{\boldsymbol{m}}
\def\bn{\boldsymbol{n}}
\def\bp{\boldsymbol{p}}

\def\bu{\boldsymbol{u}}
\def\bsv{\boldsymbol{v}}

\def\w{{\boldsymbol w}}

\def\x{\boldsymbol{x}}
\def\y{{\boldsymbol y}}

\def\bsA{{\boldsymbol A}}

\def\bC{{\boldsymbol C}}

\def\bsD{{\boldsymbol D}}

\def\bsH{{\boldsymbol H}}

\def\bsI{{\boldsymbol I}}

\def\bQ{{\boldsymbol Q}}

\def\bV{{\boldsymbol V}}

\def\bX{{\boldsymbol X}}
\def\bY{{\boldsymbol Y}}
\def\bsW{{\boldsymbol W}}

\def\bZ{{\boldsymbol Z}}

\def\mC{{\mathcal C}}

\def\mH{{\mathcal H}}
\def\mI{{\mathcal I}}

\def\mP{{\mathcal P}}

\def\mR{{\mathcal R}}
\def\mS{{\mathcal S}}
\def\mT{{\mathcal T}}

\def\mV{{\mathcal V}}
\def\mW{{\mathcal W}}

\def\bmH{{\boldsymbol \mH}}
\def\bmI{{\boldsymbol \mI}}

\def\bmP{{\boldsymbol \mP}}

\def\bmR{{\boldsymbol \mR}}
\def\bmS{{\boldsymbol \mS}}
\def\bmT{{\boldsymbol \mT}}

\def\bmW{{\boldsymbol \mW}}

\def\msF{{\mathscr F}}

\def\msH{{\mathscr H}}

\def\bmsF{{\boldsymbol \msF}}

\def\dde{\boldsymbol{\delta}}

\def\ssi{\boldsymbol{\sigma}}

\def\bvphi{\boldsymbol{\varphi}}

\def\bi{\begin{itemize}} \def\ei{\end{itemize}}
\def\be{\begin{eqnarray*}}
\def\ee{\end{eqnarray*}}

\def\0{{\mathbf 0}}

\newcommand{\beq}{\begin{equation}}
\newcommand{\eeq}{\end{equation}}

\def\xxi{{\boldsymbol{\xi}}}

\def\La{{\boldsymbol\Lambda}}

\def\Ph{\boldsymbol{\Phi}}
\def\Sig{\boldsymbol{\Sigma}}

\def\wh{\widehat}

\def\Na{\boldsymbol \nabla}

\newcommand{\eps}{\varepsilon}
\def\la{\langle}
\def\ra{\rangle}

\def\XXint#1#2#3{{\setbox0=\hbox{$#1{#2#3}{\int}$ }
\vcenter{\hbox{$#2#3$ }}\kern-.55\wd0}}

\newcommand{\TheTitle}{Data Driven Tight Frame for Compressed Sensing MRI Reconstruction via Off-the-Grid Regularization}
\newcommand{\TheAuthors}{Jian-Feng Cai, Jae Kyu Choi, and Ke Wei}
\newcommand{\argmin}{\operatornamewithlimits{argmin}}

\headers{DDTF for CS-MRI via Off-the-Grid Regularization}{\TheAuthors}

\title{{\TheTitle}\thanks{Submitted to the editors DATE. \funding{The research of the first author is supported by the Hong Kong Research Grant Council Grant 16306317. The research of the second author is supported in part by the NSFC Youth Program 11901436. The research of the third author is supported by the NSFC Youth Program 11801088.}}}

\author{
  Jian-Feng Cai\thanks{Department of Mathematics, Hong Kong University of Science and Technology, Clearwater Bay, Kowloon, Hong Kong, China (\email{jfcai@ust.hk}).}
  \and
  Jae Kyu Choi\thanks{Corresponding Author. School of Mathematical Sciences, Tongji University, 1239 Siping Road, Shanghai, 200092 China (\email{jaycjk@tongji.edu.cn}).}
  \and
  Ke Wei\thanks{Corresponding Author. School of Data Science, Fudan University, 220 Handan Road, Shanghai, 200433 China (\email{kewei@fudan.edu.cn}).}}

\usepackage{amsopn}
\DeclareMathOperator{\diag}{diag}

\begin{document}
\maketitle

\begin{abstract} Recently, the finite-rate-of-innovation (FRI) based continuous domain regularization is emerging as an alternative to the conventional on-the-grid sparse regularization for the compressed sensing (CS) due to its ability to alleviate the basis mismatch between the true support of the shape in the continuous domain and the discrete grid. In this paper, we propose a new off-the-grid regularization for the CS-MRI reconstruction. Following the recent works on two dimensional FRI, we assume that the discontinuities/edges of the image are localized in the zero level set of a band-limited periodic function. This assumption induces the linear dependencies among the Fourier samples of the gradient of the image, which leads to a low rank two-fold Hankel matrix. We further observe that the singular value decomposition of a low rank Hankel matrix corresponds to an adaptive tight frame system which can represent the image with sparse canonical coefficients. Based on this observation, we propose a data driven tight frame based off-the-grid regularization model for the CS-MRI reconstruction. To solve the nonconvex and nonsmooth model, a proximal alternating minimization algorithm with a guaranteed global convergence is adopted. Finally, the numerical experiments show that our proposed data driven tight frame based approach outperforms the existing approaches.
\end{abstract}



\section{Introduction}\label{Introduction}

Magnetic resonance imaging (MRI) is one of the most widely used medical imaging modality in clinical diagnosis \cite{E.M.Haacke1999}. It is non-radioactive, non-invasive, and has excellent soft tissue contrasts such as T1 and T2 with high spatial resolution \cite{Y.Liu2015}. Among these merits, the availability of high spatial resolution images, which will be the focus of this paper, facilitates early diagnosis by enabling the detection and characterization of clinically important lesions \cite{G.Ongie2015,E.VanReeth2012}. However, since the so-called k-space data acquisition is limited due to physical (gradient amplitude and slew-rate) and physiological (motion artifacts due to the patient's respiratory motion which can generate the outlier of the k-space data) constraints \cite{C.M.Hyun2018,K.H.Jin2017,Y.Liu2015,M.Lustig2007}, there has been increasing demand for methods which can reduce the amount of acquired data without degrading the image quality \cite{M.Lustig2007}.

When the k-space data is undersampled, the Nyquist sampling criterion is violated, and the \emph{direct restoration} using this undersampled k-space data will lead to the aliasing in the reconstructed image \cite{M.Lustig2007}. In the literature, the famous compressed sensing (CS) MRI can be viewed as a sub-Nyquist sampling method which exploits the sparse representation of an image to compensate the undersampled k-space data \cite{E.J.Candes2006,Donoho2006,C.M.Hyun2018,M.Lustig2007}. More precisely, a typical $\ell_1$ norm based ``on-the-grid'' CS-MRI model restores the MR image on the grid, by regularizing smooth image components while preserving image singularities such as edges, ridges, and corners.

Even though the conventional on-the-grid approaches have shown strong ability to reduce the data acquisition time and thus received a lot of attention over the past few years \cite{Y.Liu2015}, they have to be further improved as the k-space data is a discrete (and truncated) sampling of a Fourier transform of an underlying function in the \emph{continuous domain}, under the assumption of a single-coil MRI \cite{E.M.Haacke1999}. This means that, the sparse regularization based image restoration will work well when the singularities of the image are well aligned with the grid. However, even in the case of piecewise constant function whose (distributional) gradient is sparse in the continuous domain, its singularities may not necessarily agree with the discrete grid, leading to the problem of \emph{basis mismatch} \cite{Y.Chi2010,G.Ongie2016,J.Ying2017}. Such a basis mismatch between the true singularities and the discrete grid may result in the loss of sparse structure of the image, and thus degrades the restoration quality \cite{J.Ying2017}.

\subsection{Continuous domain regularization for image restoration}

The continuous domain regularization is emerging as a powerful alternative to the discrete domain sparse regularization \cite{B.N.Bhaskar2013,E.J.Candes2014,Y.Chen2014,G.Ongie2018}. By such an ``off-the-grid'' scheme, the regularization on the discrete data exploits the sparsity in the continuous domain to alleviate the basis mismatch, which is especially attractive in signal restorations from partial measurements \cite{G.Ongie2018}. One of the most successful examples is the finite-rate-of-innovation (FRI) framework which corresponds a superposition of a few sinusoids to a low rank Hankel matrix \cite{T.Blu2008}. Similarly, we can adopt the so-called \emph{structured low rank matrix completion (SLMC)} \cite{Y.Chen2014,G.Ongie2017,J.C.Ye2017} for the CS restoration based on the continuous domain regularization. The SLMC approach can also be used to restore Fourier samples of a one dimensional piecewise constant signal \cite{T.Blu2008,M.Vetterli2002} by considering the Fourier samples of a derivative. However, while this framework works well in the case of isolated singularities \cite{E.J.Candes2013,W.Xu2014}, the extension to the (piecewise smooth) image restoration is not straightforward \cite{G.Ongie2016}. This is because the image singularities such as the edges and ridges in general form a continuous curve in a two dimensional domain, which makes it challenging to construct a structured low rank matrix from the Fourier samples.

Recently, there are several extensions of the FRI framework to two dimensional continuous domain regularization for image restoration. Such extensions include the piecewise holomorphic complex image restoration \cite{H.Pan2014} and the piecewise constant real image restoration \cite{G.Ongie2018,G.Ongie2015a,G.Ongie2015,G.Ongie2016}. These approaches are commonly based on the assumption that the singularity curves, i.e. the supports of the (real/complex) derivatives of a target image, lie in the zero level set of a band-limited periodic function, called the \emph{annihilating polynomial}. Indeed, given that the edges lie in the zero level set of an annihilating polynomial, the structured matrix (Hankel/Toeplitz matrix) constructed from the Fourier samples of piecewise constant images becomes low rank \cite{G.Ongie2018,G.Ongie2015,G.Ongie2016}. Hence, in CS-MRI, we can first restore the fully sampled k-space data (the discrete sample of the Fourier transform of a piecewise constant function) via the SLMC, and then restore MR image by the inverse DFT \cite{Haldar2014,K.H.Jin2016}.


Though the rank minimization problem is an NP-hard problem in general \cite{E.J.Candes2009}, numerous tractable relaxation approaches have been proposed for the SLMC. One of them is the convex nuclear norm relaxation method (e.g. \cite{J.F.Cai2016,M.Fazel2013}), together with the theoretical restoration guarantees (e.g. \cite{G.Ongie2018}). In addition, the iterative reweighted least squares (IRLS) approaches for the Schatten $p$-norm minimization are proposed in \cite{M.Fornasier2011,K.Mohan2012,G.Ongie2017}, which can avoid the high computational cost of the SVD related to the rank minmization and the convex nuclear norm relaxation. Apart from these approaches, the nonconvex alternating projection methods based on the different parametrizations of the underlying low rank matrix structure are proposed and studied in \cite{J.F.Cai2018a,J.F.Cai2019}. These nonconvex methods are reported to be superior to the other relaxation methods in terms of the computational efficiency while theoretical restoration guarantees are still available.


\subsection{Motivations and contributions of our approach}

Our proposed off-the-grid approach for the CS-MRI restoration is based on the fact that the SVD of a Hankel matrix induces tight frame filter banks due to its underlying convolutional structure. This means that, if we can associate a signal with a low rank Hankel matrix, its right singular vectors form tight frame filter banks which allow us to represent the signal with sparse canonical coefficients. Hence, instead of directly minimizing the rank of a two-fold Hankel matrix, we can develop the sparse regularization model via data driven tight frames \cite{J.F.Cai2014} for the CS-MRI restoration. More specifically, using the $\ell_0$ norm of the frame coefficients, we propose the \emph{balanced approach} \cite{J.F.Cai2008,R.H.Chan2003} as a relaxation of the low rank two-fold Hankel matrix approximation which is based on the framework of annihilating filters. Finally, the numerical experiments show that our data driven tight frame approach outperforms the structured low rank matrix approaches (and their relaxations) as well as the existing on-the-grid approaches, leading to the state-of-the-art performance.

Note that the data driven tight frame as well as the dictionary learning has been discussed for the various image restoration tasks. More precisely, the data driven tight frame is first used for the image denoising \cite{J.F.Cai2014}, the multi-channel image denoising and inpainting \cite{J.Wang2015}, the joint spatial-Radon domain regularization for the sparse view CT \cite{R.Zhan2016}, and the PET-MRI joint reconstruction \cite{J.K.Choi2018}. In addition, the authors in \cite{S.Ravishankar2015} explored the general framework of the dictionary learning based image restoration model with the convergence analysis. However, while these aforementioned approaches mainly concern the sparse approximation of the discrete images, to the best of our knowledge, our approach for the first time applies the data driven tight frame regularization to the MRI restoration based on a continuous image model.

Finally, we would like to mention that the proposed off-the-grid approach for the CS-MRI restoration is not limited to the two dimensional single-coil MRI restoration. Even though further numerical studies are required, it is not hard to generalize the proposed model to the three dimensional cases. In addition, since the multiplication in the image domain becomes the convolution in the frequency domain via the Fourier transform, we can extend our approach to a SENSE setup \cite{K.P.Pruessmann1999}, a multiple-coil MRI restoration whose inverse problem corresponds to the multi-channel deconvolution in the frequency domain with a known coil sensitivity map.

\subsection{Organization and notation of paper}

The rest of this paper is organized as follows. We first briefly review the tight wavelet frames including the wavelet frame based image restorations and the data driven tight frames in \cref{TightFrames}. \cref{ModelandAlgorithm} begins with the review on the structured low rank matrix approaches for CS-MRI. Then we present the off-the-grid CS-MRI reconstruction approach based on data driven tight frames, followed by the alternating minimization algorithm. In \cref{Experiments}, experimental results are reported to demonstrate the performance of our new CS-MRI restoration method, and \cref{Conclusion} concludes this paper with a few future directions.

\begin{figure}[tp!]
\centering
\includegraphics[width=0.95\textwidth]{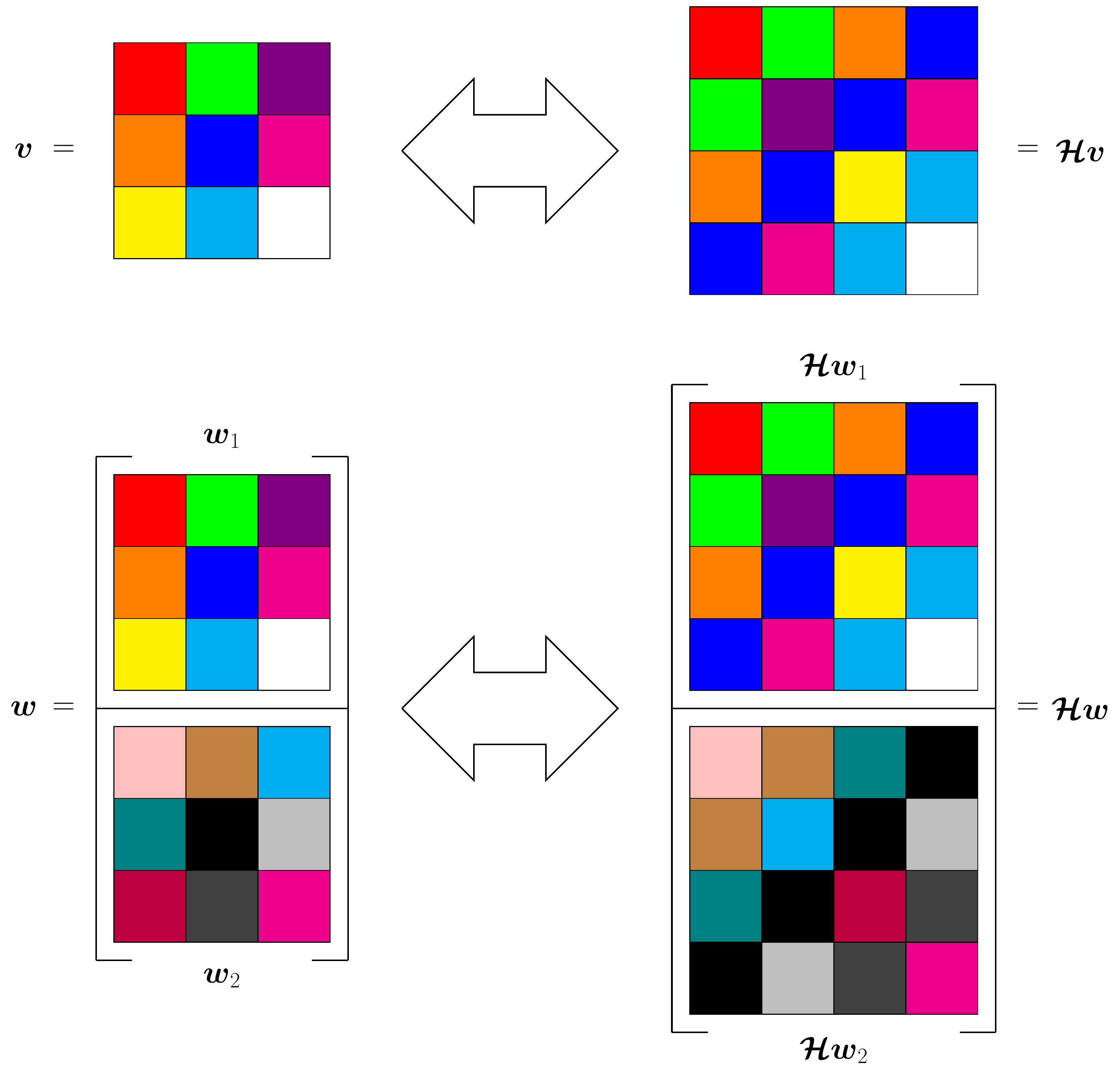}\vspace{-0.20cm}
\caption{Illustration of generating a Hankel matrix $\bmH\bsv$ from $\bsv\in\mV_2$ and a two-fold Hankel matrix $\bmH\w$ from $\w=\left(\w_1,\w_2\right)\in\mV_2\times\mV_2$. Here, $K_1=K_2=2$.}\label{BlockHankelIllustrate}
\end{figure}

Throughout this paper, all two dimensional discrete images will be denoted by the bold faced lower case letters while the functions in a two dimensional continuous domain will be denoted by regular characters. Note that a two dimensional image can also be identified with a vector whenever convenient. All matrices will be denoted by the bold faced upper case letters, and the $j$th row and the $k$th column of a matrix $\bZ$ will be denoted by $\bZ^{(j,:)}$ and $\bZ^{(:,k)}$, respectively. We denote by
\begin{align}\label{ImageGrid}
\OO=\left\{-\lfloor N/2\rfloor,\cdots,\lfloor N/2-1\rfloor\right\}^2,
\end{align}
with $N\in\N$, the set of $N\times N$ grid, and the space of complex valued functions on $\OO$ is denoted by $\mV_2\simeq\C^{|\OO|}$. Given two rectangular grids $\KK$ and $\MM$, the contraction $\KK:\MM$ is defined as
\begin{align*}
\KK:\MM=\left\{\bk\in\KK:\bk-\bsj\in\KK~\text{for all}~\bsj\in\MM\right\}.
\end{align*}
Finally, operators are denoted as the bold faced calligraphic letters. For instance, let $\bsv\in\mV_2$ and let $\KK$ be a rectangular $K_1\times K_2$ grid. The corresponding Hankel matrix $\bmH\bsv$ is an $M_1\times M_2$ matrix ($M_1=|\OO:\KK|$ and $M_2=|\KK|$) generated by concatenating $K_1\times K_2$ patches of $\bsv$ into row vectors. In the sense of two dimensional multi-indices, we can formally write $\bmH\bsv$ as
\begin{align*}
\left(\bmH\bsv\right)[\bk,\bsl]=\bsv[\bk+\bsl],~~~~~~\bk\in\OO:\KK~~~\text{and}~~~\bsl\in\KK.
\end{align*}
With a slight abuse of notation, we also use
\begin{align}\label{TwoFoldHankel}
\bmH\w=\left[\begin{array}{cc}
\left(\bmH\w_1\right)^T&\left(\bmH\w_2\right)^T
\end{array}\right]^T,
\end{align}
to denote the $2M_1\times M_2$ two-fold Hankel matrix constructed from $\w=\left(\w_1,\w_2\right)\in\mV_2\times\mV_2$; see \cref{BlockHankelIllustrate} for an illustration.


\section{Preliminaries on tight wavelet frames}\label{TightFrames}

Here we provide a brief introduction on tight wavelet frames and data driven tight frames. Interested readers may consult \cite{J.F.Cai2008,B.Dong2013,B.Dong2015,A.Ron1997,Shen2010} for detailed surveys on tight wavelet frames, and \cite{C.Bao2015,J.F.Cai2014} for details on data driven tight frames. For the sake of simplicity, we only discuss the real valued wavelet tight frame systems, but note that it is not difficult to extend the idea to the complex case.

Denote by $\msH$ a Hilbert space and let $\la\cdot,\cdot\ra$ be the inner product defined on $\msH$. A countable set $\{\bvphi_n:n\in\Z\}\subseteq\msH$ is called a tight frame on $\msH$ if
\begin{align}\label{TightFrame}
\|\bu\|^2=\sum_{n\in\Z}|\la\bu,\bvphi_n\ra|^2~~~~~\text{for all}~~~\bu\in\msH.
\end{align}
Given $\{\bvphi_n:n\in\Z\}\subseteq\msH$, we define the analysis operator $\bmW:\msH\to\ell_2(\Z)$ as
\begin{align*}
\bu\in\msH\mapsto \bmW\bu=\{\la\bu,\bvphi_n\ra:n\in\Z\}\in\ell_2(\Z).
\end{align*}
The synthesis operator $\bmW^T:\ell_2(\Z)\to\msH$ is defined as the adjoint of $\bmW$:
\begin{align*}
\bc\in\ell_2(\Z)\mapsto\bmW^T\bc=\sum_{n\in\Z}\bc[n]\bvphi_n\in\msH.
\end{align*}
Then $\{\bvphi_n:n\in\Z\}$ is a tight frame on $\msH$ if and only if $\bmW^T\bmW=\bmI$ where $\bmI$ is the identity on $\msH$. It follows that, for a given tight frame $\{\bvphi_n:n\in\Z\}$, we have the following canonical expression:
\begin{align*}
\bu=\sum_{n\in\Z}\la\bu,\bvphi_n\ra\bvphi_n,
\end{align*}
with $\bmW\bu=\{\la\bu,\bvphi_n\ra:n\in\Z\}$ being called the canonical tight frame coefficients. Hence, the tight frames are extensions of orthonormal bases to the redundant systems. In fact, a tight frame is an orthonormal basis if and only if $\|\bvphi_n\|=1$ for all $n\in\Z$.

One of the most widely used class of tight frames is the discrete wavelet frame generated by a set of finitely supported filters $\{\a_1,\cdots,\a_m\}$. Throughout this paper, we only discuss the two dimensional undecimated wavelet frames on $\ell_2(\Z^2)$, but note that it is not difficult to extend to $\ell_2(\Z^d)$ with $d\geq3$. For $\a\in\ell_1(\Z^2)$, define a convolution operator $\bmS_{\a}:\ell_2(\Z^2)\to\ell_2(\Z^2)$ by
\begin{align}\label{DiscreteConv}
(\bmS_{\a}\bu)[\bk]=(\a\ast\bu)[\bk]=\sum_{\bsl\in\Z^2}\a[\bk-\bsl]\bu[\bsl]~~~~~\text{for}~~~\bu\in\ell_2(\Z^2).
\end{align}
Given a set of finitely supported filters $\{\a_1,\cdots,\a_m\}$, define the analysis operator $\bmW$ and the synthesis operator $\bmW^T$ by
\begin{align}
\bmW&=\left[\bmS_{\a_1[-\cdot]}^T,\bmS_{\a_2[-\cdot]}^T,\cdots,\bmS_{\a_m[-\cdot]}^T\right]^T,\label{AnalConv}\\
\bmW^T&=\big[\bmS_{\a_1},\bmS_{\a_2},\cdots,\bmS_{\a_m}\big],\label{SyntConv}
\end{align}
respectively. Then, the direct computation can show that the rows of $\bmW$ form a tight frame on $\ell_2(\Z^2)$ (i.e. $\bmW^T\bmW=\bmI$) if and only if the filters $\{\a_1,\cdots,\a_m\}$ satisfy
\begin{align}\label{SomeUEP}
\sum_{j=1}^m\sum_{\bsl\in\Z^2}\a_j[\bk+\bsl]\a_j[\bsl]=\dde[\bk]=\left\{\begin{array}{cl}
1~&\text{if}~\bk=\0,\\
0~&\text{if}~\bk\neq\0,
\end{array}\right.
\end{align}
called the {\emph{unitary extension principle}} (UEP) condition \cite{B.Han2011}. Finally, for a two dimensional discrete image on the finite grid, $\bmS_{\a}$ in \cref{DiscreteConv} with a finitely supported $\a$ denotes the discrete convolution under the periodic boundary condition throughout this paper.

In the literature, wavelet frames are widely used for the sparse approximation of an image. To explore a better sparse approximation, the authors in \cite{J.F.Cai2014} proposed a \emph{data driven tight frame} approach, inspired by \cite{M.Aharon2006,S.Ravishankar2013}. Specifically, given a two dimensional image $\bu\in\mV_2$, a tight frame system $\bmW$ defined as in \cref{AnalConv}, which is generated by a set of finitely supported $p\times p$ filters $\{\a_1,\cdots,\a_{p^2}\}$ supported on a $p\times p$ grid $\KK$ and satisfying \cref{SomeUEP}, is constructed via the following minimization
\begin{align}\label{DDTFModel}
\min_{\bc,\bmW}~\|\bc-\bmW\bu\|_2^2+\gamma^2\|\bc\|_0~~~~~\text{subject to}~~~~\bmW^T\bmW=\bmI,
\end{align}
with the $\ell_0$ norm $\|\bc\|_0$ encoding the number of nonzero entries in the coefficient vector $\bc$. Note that in wavelet domain, there is usually a set of low pass coefficients that are not sparse. In fact, through the relaxation of \cref{DDTFModel} which will be discussed below, we can learn a tight frame system $\bmW$ by fixing a low pass filter and penalizing the $\ell_0$ norm of high pass coefficients only, following \cite{C.Bao2013}. However, since the main focus is to establish a relation between the tight frame filter banks and the Hankel matrix, we forgo further discussions on such a detailed construction.

To solve \cref{DDTFModel}, the authors in \cite{J.F.Cai2014} presented that if $\bsA\in\R^{p^2\times p^2}$ satisfies $\bsA\bsA^T=p^{-2}\bsI$, after reformulating into $p\times p$ filters supported on $\KK$, its column vectors generate filters satisfying the UEP condition \cref{SomeUEP}. Hence, under such an assumption on the filters $\a_1,\cdots,\a_{p^2}$, we can rewrite \cref{DDTFModel} as
\begin{align}\label{DDTFModelModi}
\min_{\bC,\bsA}~\|\bC-\left(\bmH\bu\right)\bsA\|_F^2+\gamma^2\|\bC\|_0~~~~~\text{subject to}~~~~\bsA\bsA^T=p^{-2}\bsI,
\end{align}
to construct filter banks. Here, $\bmH\bu\in\R^{(N-p+1)^2\times p^2}$ is a Hankel matrix generated by the $p\times p$ patches of $\bu$, $\bC$ is the corresponding matrix formulation of $\bc$, and $\left\|\cdot\right\|_F$ is the Frobenius norm of a matrix. To solve \cref{DDTFModelModi}, the alternating minimization method with closed form solutions for each stage is presented in \cite{J.F.Cai2014}. In addition, the proximal alternating minimization (PAM) scheme with global convergence to critical points is proposed in \cite{C.Bao2015}.

\section{Proposed off-the-grid regularization model}\label{ModelandAlgorithm}

\subsection{Structured low rank matrix approaches in CS-MRI}\label{SRLACSMRI}

We begin with the brief review on the structured low rank matrix approaches in CS-MRI. Interested readers can refer to \cite{G.Ongie2016} and references therein for more detailed surveys. Throughout this paper, we only consider the two dimensional single-coil MRI reconstruction for the sake of simplicity, so that the k-space data $\bsf$ can be modeled as
\begin{align}\label{MRForward}
\bsf[\bk]=\wh{u}(L^{-1}\bk)=\int_{\R^2}u(\x)e^{-2\pi i\x\cdot\bk/L}\rd\x,~~~~~\bk\in\MM,
\end{align}
where $\MM$ is the sampling region, $L>0$ is the length of field-of-view (FOV), and $u\in L_1(\R^2)$ is the proton spin density distribution in $\R^2$ \cite{E.M.Haacke1999}. We further assume that $u$ is supported on $[-L/2,L/2)^2$. This means that when $\bsf[\bk]$ in \cref{MRForward} is available for $\bk\in\Z^2$, we have
\begin{align*}
u(\x)=L^{-2}\sum_{\bk\in\Z^2}\bsf[\bk]e^{2\pi i\bk\cdot\x/L},~~~~~\x\in[-L/2,L/2)^2.
\end{align*}
In general, the k-space data can be sampled at ``arbitrary'' positions, i.e. $\MM\subseteq\R^2$. However, to simplify the discussion, we assume that $\MM\subseteq\OO\subseteq\Z^2$ where $\OO$ is defined as \cref{ImageGrid}, i.e. the discrete sampling on the Nyquist rate grid \cite{E.M.Haacke1999,T.H.Kim2019}.

The conventional CS-MRI approach aims to directly restore the MR image $\bu$ on $\OO$ from the k-space $\bsf$ undersampled on $\MM$. Note that when $\MM=\OO$, i.e. the fully sampled case, $\bu$ is obtained from the inverse DFT
\begin{align}\label{MRInverse}
\bu[\bp]=\left(\bmsF^{-1}\bsf\right)[\bp]=N^{-2}\sum_{\bk\in\OO}\bsf[\bk]e^{2\pi i\bp\cdot\bk/N},~~~~~~~\bp\in\OO.
\end{align}
Hence, the conventional sparse regularization based CS-MRI approach
\begin{align}\label{CSApproach}
\min_{\bu}~\f{1}{2}\left\|\bmR_{\MM}\bmsF\bu-\bsf\right\|_2^2+\gamma\left\|\Ph\bu\right\|_1
\end{align}
is based on the DFT
\begin{align}\label{DFT}
\bsf[\bk]=\left(\bmsF\bu\right)[\bk]=\sum_{\bp\in\OO}\bu[\bp]e^{-2\pi i\bp\cdot\bk/N},~~~~~~~\bk\in\MM.
 \end{align}
Here, $\bmR_{\MM}$ is a restriction onto $\MM$, and $\Ph$ is a sparsifying transformation. Note that \cref{CSApproach} is the on-the-grid scheme, as we promote the sparsity of the discrete image $\bu$ on the grid $\OO$. In other words, we can apply the sparse regularization effectively provided that the singularities of $u$ are well aligned with $\OO$. However, since the k-space data actually comes from the Fourier transform of a continuous domain function $u$, together with \cref{MRForward} and $\MM=\OO$, \cref{MRInverse} can be rewritten as
\begin{align*}
\bu[\bp]=\left(N^{-1}L\right)^2\left(D_N\ast u\right)(N^{-1}L\bp),~~~~~\bp\in\OO,
\end{align*}
where $D_N(\x)$ is the Dirichlet kernel defined as
\begin{align*}
D_N(\x)=L^{-2}\sum_{\bk\in\OO}e^{2\pi i\bk\cdot\x/L}.
\end{align*}
In other words, since the MR image is proportional to a discrete sampling of $D_N\ast u$ on $N^{-1}L\OO$ rather than that of $u$ \cite{E.M.Haacke1999}, there exists a basis mismatch between the singularities of $u$ in the continuous domain and the discrete grid $\OO$. Such a basis mismatch would destroy the sparsity structure, leading to a degradation of the restoration performance.

Instead of directly restoring the discrete MR image, the off-the-grid approaches attempt to first restore the fully sampled k-space data $\bsv$ from its undersampled version $\bsf$ modeled as in \cref{MRForward}. Formally, we can write the problem as follows:
\begin{align}\label{OurProblem}
\text{find}~~\bsv=\msF(u)\big|_{L^{-1}\OO}~~\text{subject to}~~\bmR_{\MM}\bsv=\bsf,
\end{align}
and then compute $\bu=\bmsF^{-1}\bsv$ by using the restored k-space data $\bsv$.

\begin{figure}[tp!]
\centering
\includegraphics[width=0.95\textwidth]{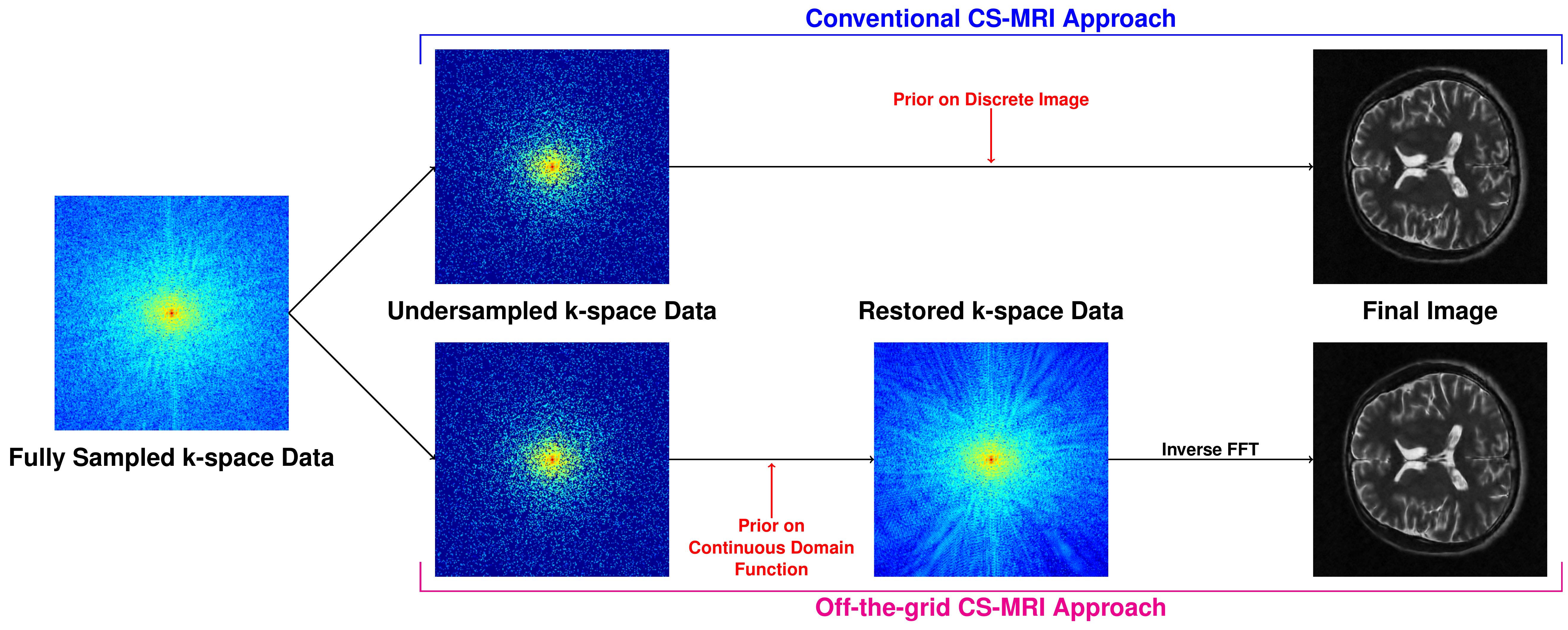}\vspace{-0.20cm}
\caption{Comparison between the conventional CS-MRI and the off-the-grid CS-MRI approach.}\label{CSMRIComparisons}
\end{figure}

Generally, it is impossible to directly solve \cref{OurProblem} without any further information. Nevertheless, by imposing the prior information of $u$ (in the continuous domain) into \cref{OurProblem}, we can exploit the continuous domain regularization, and reduce the basis mismatch; see \cref{CSMRIComparisons} for the schematic comparison. For this purpose, we follow \cite{G.Ongie2016} and consider the following piecewise constant function model of the proton density $u$:
\begin{align}\label{uModel}
u(\x)=\sum_{j=1}^J\alpha_j1_{\Om_j}(\x),~~~~~~~\x\in\R^2,
\end{align}
where $\alpha_j\in\C$ and $1_{\Om}$ denotes the characteristic function on a set $\Om$: $1_{\Om}(\x)=1$ if $\x\in\Om$, and $0$ otherwise. We further assume that each $\Om_j$ lies in $[-L/2,L/2)^2$, and \cref{uModel} is expressed with the smallest number of characteristic functions such that $\Om_j$'s are pairwise disjoint. Then the discontinuities of $u$ agrees with $\Gamma=\bigcup_{j=1}^J\p\Om_j$, called the \emph{edge set} of $u$ \cite{G.Ongie2016}. Under this setting, the authors in \cite{G.Ongie2015,G.Ongie2016} derived the following \emph{linear annihilation relation}.

\begin{proposition}\label{Prop1} Let $u(\x)$ be a piecewise constant function defined as \cref{uModel}. Assume that there exists a finite (rectangular and symmetric) index set $\KK$ such that
\begin{align}\label{MajorAssumption}
\Gamma\subseteq\left\{\x\in\R^2:\varphi(\x)=0\right\}~~~~~\text{with}~~\varphi(\x)=\sum_{\bk\in\KK}\a[\bk]e^{2\pi i\bk\cdot\x/L}.
\end{align}
Then for such $\varphi$, the Fourier transform of $\na u=\left(\p_1u,\p_2u\right)$ is annihilated by the convolution with the Fourier transform of $\varphi$. More precisely, we have
\begin{align}\label{FouGradConvSemiDiscrete}
\left(\msF(\na u)\ast\wh{\varphi}\right)(\xxi)=\sum_{\bk\in\KK}\msF(\na u)(\xxi-L^{-1}\bk)\a[\bk]=0,~~~~~~\xxi\in\R^2,
\end{align}
where the convolution acts on $\msF(\p_1 u)$ and $\msF(\p_2 u)$ separately.
\end{proposition}

The trigonometric polynomial $\varphi$ in \cref{MajorAssumption} is called the \emph{annihilating polynomial}, $\a$ is called the \emph{annihilating filter} with support $\KK$, and $\Gamma$ satisfying \cref{MajorAssumption} is called the \emph{trigonometric curves} \cite{G.Ongie2016}. It is proved in \cite{G.Ongie2016} that trigonometric curves can always be described as the zero set of a \emph{real valued} trigonometric polynomial, using the tool of algebraic geometries.

In MRI, the Fourier transform of $u$ is sampled on the grid $L^{-1}\OO$, so \cref{FouGradConvSemiDiscrete} becomes the finite system of linear equations
\begin{align}\label{FouGradConvDiscrete}
\sum_{\bk\in\KK}\msF(\na u)(L^{-1}(\bsl+\bk))\a[-\bk]=0,~~~~~~~\bsl\in\OO:\KK.
\end{align}
In the matrix-vector multiplication form, we have
\begin{align}\label{RankDeficientHankel}
\bmH\left(\msF(\na u)\big|_{L^{-1}\OO}\right)\a[-\cdot]=\0.
\end{align}
Moreover, if $\varphi$ in \cref{MajorAssumption} is the minimal polynomial with the support $\KK$ (i.e. $\KK$ has the smallest linear dimension) and $\KK'$ is the assumed filter support strictly containing $\KK$, then we have
\begin{align}\label{RankUB}
\rank\left(\bmH\left(\msF(\na u)\big|_{L^{-1}\OO}\right)\right)\leq|\KK'|-|\KK':\KK|.
\end{align}
Roughly speaking, \cref{RankUB} means that if $\varphi$ in \cref{MajorAssumption} is the minimal polynomial with $\a$ satisfying \cref{FouGradConvDiscrete,RankDeficientHankel}, so does the translation $\a[\cdot-\bm]$ for all $\bm\in\KK'\setminus\left(\KK':\KK\right)$, or equivalently, $e^{2\pi i\bm\cdot\x}\varphi(\x)$ is also an annihilating polynomial. See \cite{G.Ongie2015,G.Ongie2016} for details.

Note that \cref{RankUB} means the two-fold Hankel matrix constructed by $\msF(\na u)\big|_{L^{-1}\OO}$ is \emph{rank deficient} whenever the grid $\OO$ in \cref{ImageGrid} is large enough to satisfy the conditions in \cite[Proposition 5.1]{G.Ongie2016}. Based on this observation, we can formulate \cref{OurProblem} as the following low rank two-fold Hankel matrix completion \cite{Y.Chen2014,G.Ongie2017,J.C.Ye2017}
\begin{align}\label{LRSMCompletionRevist}
\min~~\rank\left(\bmH\left(\La\bsv\right)\right)~~~~~~~\text{subject to}~~\bmR_{\MM}\bsv=\bsf.
\end{align}
Here, $\La\bsv=\left(\La_1\bsv,\La_2\bsv\right)\in\mV_2\times\mV_2$ and $\La\in\C^{2|\OO|\times|\OO|}$ is a weight matrix defined as
\begin{align}\label{Lambda}
\La=\left[\begin{array}{cc}
\La_1&\La_2
\end{array}\right]^T=\left[\begin{array}{cc}
\diag\left(2\pi ik_1/L\right)&\diag\left(2\pi ik_2/L\right)
\end{array}\right]_{\bk=(k_1,k_2)\in\OO}^T,
\end{align}
which is derived from the sample of $\msF\left(\na u\right)(\xxi)=2\pi i\xxi\wh{u}(\xxi)=\left(2\pi i\xi_1\wh{u}(\xxi),2\pi i\xi_2\wh{u}(\xxi)\right)$ on $L^{-1}\OO$. Hence, the piecewise constant property of $u$ in the continuous domain can be transformed to the low rank property of $\bmH\left(\La\bsv\right)$.

\subsection{Proposed CS-MRI model}\label{ProposedModelandAlgorithm}

Denote by $\bsv\in\mV_2$ the k-space data on the grid $\OO$, which is to be restored from a given data $\bsf$ undersampled on $\MM$. Consider a two-fold Hankel matrix $\bmH\left(\La\bsv\right)\in\C^{2M_1\times M_2}$ defined as \cref{TwoFoldHankel} with $\La$ defined as \cref{Lambda}. To reconstruct $\bsv=\msF(u)\big|_{L^{-1}\OO}$ corresponding to a piecewise constant function $u$ in \cref{uModel}, we assume that $\rank\left(\bmH\left(\La\bsv\right)\right)=r\ll2M_1\wedge M_2$, following \cref{SRLACSMRI}. Considering its SVD, we present \cref{Th1}, which is the key idea of the proposed CS-MRI model. The proof can be found in \cref{ProofTh1}.

To do this, we let $\a_1,\cdots,\a_{M_2}$ be $K_1\times K_2$ filters supported on $\KK$, and we introduce
\begin{align}
\bmW&=\left[\bmS_{\a_1[-\cdot]}^T,\bmS_{\a_2[-\cdot]}^T,\cdots,\bmS_{\a_{M_2}[-\cdot]}^T\right]^T,\label{OurAnalysis}\\
\bmW^*&=\left[\bmS_{\overline{\a}_1},\bmS_{\overline{\a}_2},\cdots,\bmS_{\overline{\a}_{M_2}}\right],\label{OurSynthesis}
\end{align}
where $\bmS_{\a}$ is a discrete convolution under the periodic boundary condition defined as \cref{DiscreteConv}. In other words, both $\bmW$ and $\bmW^*$ are concatenations of discrete convolutions.

\begin{theorem}\label{Th1} Let $u(\x)$ be defined as in \cref{uModel} and $\bsv=\msF(u)\big|_{L^{-1}\OO}$. Let $\bmH\left(\La\bsv\right)\in\C^{2M_1\times M_2}$ be a two-fold Hankel matrix with $\La$ defined as \cref{Lambda}. Assume that $\rank\left(\bmH\left(\La\bsv\right)\right)=r\ll2M_1\wedge M_2$. Considering its full SVD $\bmH\left(\La\bsv\right)=\bX\Sig\bY^*$, we let $\a_j=M_2^{-1/2}\bY^{(:,j)}$ by reformulating $\bY^{(:,j)}\in\C^{M_2}$ into a $K_1\times K_2$ filter supported on $\KK$. Then $\bmW$ in \cref{OurAnalysis} defined by using these $\a_1,\cdots,\a_{M_2}$ satisfies
\begin{align}\label{LowRankHankelTightFrame}
\bmW^*\bmW\left(\La\bsv\right)=\sum_{j=1}^{M_2}\bmS_{\overline{\a}_j}\left(\bmS_{\a_j[-\cdot]}\left(\La\bsv\right)\right)=\La\bsv,
\end{align}
and for $j=r+1,\cdots,M_2,$ we have
\begin{align}\label{TightFrameSparse}
\left(\bmS_{\a_j[-\cdot]}\left(\La\bsv\right)\right)[\bk]=\0,~~~~~\bk\in\OO:\KK,
\end{align}
where the discrete convolution $\bmS_{\a}$ acts on $\La_1\bsv$ and $\La_2\bsv$ separately\footnote{In fact, \cref{TightFrameSparse} is independent from the boundary condition of the discrete convolution, as $\OO:\KK$ is the grid where the boundary condition is inactive.}. Consequently, if $\bmH\left(\La\bsv\right)$ is of low rank, then its right singular vectors construct a tight frame under which $\La\bsv$ is sparsely represented.
\end{theorem}

In words, \cref{Th1} implies that a tight frame under which the weighted k-space data $\La\bsv$ is sparse can be explicitly constructed from the SVD of the corresponding Hankel matrix, and the sparsity of the tight frame coefficients is grouped according to filters. Even though the Hankel matrix is not fully known in the CS-MRI setting since the k-space data is undersampled, \cref{Th1} does motivate us to learn a better tight frame from the undersampled and noisy k-space data such that $\La\bsv$ can be sparsely represented better. More precisely, we remove the group sparsity pattern in the canonical coefficients $\bmW\left(\La\bsv\right)$, and find the coefficients which are as sparse as possible. This leads to solving
\begin{align}\label{AnalysisVersion}
\begin{split}
&~~~\min_{\bsv,\bmW}~\f{1}{2}\left\|\bmR_{\MM}\bsv-\bsf\right\|_2^2+\gamma\left\|\bmW\left(\La\bsv\right)\right\|_0\\
&\text{subject to}~~\bsv\in\mC~~\text{and}~~\bmW^*\bmW=\bmI.
\end{split}
\end{align}
Here, $\bmW$ is a tight frame transform defined as in \cref{OurAnalysis}, i.e. the concatenation of the convolutions with filter banks $\left\{\a_1,\a_2,\cdots,\a_{M_2}\right\}$ which have to be learned, and $\mC$ is a constraint set imposing the boundedness of $\bsv$. By relaxing the sparsity of tight frame coefficients over the entire range of $\bmW$ (not necessarily in groups as in \cref{TightFrameSparse}), we expect to achieve more flexibility, thereby leading to the improvements in restoration performance.

However, the direct penalization on the canonical coefficients $\bmW\left(\La\bsv\right)$ requires an iterative solver, which will be too expensive in constructing an adaptive tight frame \cite{J.F.Cai2014}. Therefore, by introducing an auxiliary variable $\bc=\left(\bc_1,\bc_2\right)$, we approximate \cref{AnalysisVersion} by
\begin{align}\label{ProposedCSMRIModel}
\begin{split}
&~~~\min_{\bsv,\bc,\bmW}~\f{1}{2}\left\|\bmR_{\MM}\bsv-\bsf\right\|_2^2+\f{\mu}{2}\left\|\bmW\left(\La\bsv\right)-\bc\right\|_2^2+\gamma\left\|\bc\right\|_0\\
&\text{subject to}~~\bsv\in\mC~~\text{and}~~\bmW^*\bmW=\bmI,
\end{split}
\end{align}
to learn the tight frame system and reconstruct the fully sampled k-space data $\bsv$ on the grid $\OO$ simultaneously.

\begin{remark}\label{RK1} For the constraint set $\mC$, we choose
\begin{align}\label{ConstraintSet}
\mC=\left\{\bsv\in\mV_2:|\bsv[\bk]|\leq R~~~\text{for all}~~\bk\in\OO\right\}
\end{align}
with a sufficiently large $R>0$, and such a choice comes from the following: if $u\in L_1(\R^2)$ is modeled as in \cref{uModel}, then the (distributional) gradient satisfies
\begin{align*}
\na u=\left(\p_1 u,\p_2 u\right)=-\sum_{j=1}^J\alpha_j\bn_j\rd\ssi\big|_{\p\Om_j},
\end{align*}
where $\bn_j$ denotes the outward normal vector on $\p\Om_j$, and $\ssi$ denotes the surface measure. Since $\na u$ is a Radon vector measure on $\R^2$ supported on $\Gamma$, its Fourier transform is written as the Fourier transform of a measure (e.g. \cite{Folland1999}):
\begin{align*}
\msF(\na u)(\xxi)=-\sum_{j=1}^J\alpha_j\int_{\p\Om_j}e^{-2\pi i\xxi\cdot\x}\bn_j(\x)\rd\ssi(\x),~~~~~\xxi\in\R^2,
\end{align*}
which leads to
\begin{align*}
\left|\msF(\na u)(\xxi)\right|\leq\sum_{j=1}^J|\alpha_j|\ssi(\Om_j)<\infty~~~~~\text{for all}~~\xxi\in\R^2.
\end{align*}
Hence, $\wh{u}(\xxi)$ decays at infinity in the sense that
\begin{align*}
\sup_{\xxi\in\R^2}\left|\xxi\right|\left|\wh{u}(\xxi)\right|<\infty,
\end{align*}
and this implies that $\wh{u}(\xxi)$ (and whence $\bsv[\bk]$ since we want $\bsv=\msF(u)\big|_{L^{-1}\OO}$) has to be bounded. In addition, as mentioned in e.g. \cite{C.Bao2016}, the boundedness constraint is required for the convergence guarantee of the proximal alternating scheme to solve \cref{ProposedCSMRIModel}. That being said, numerically, the constraint set $\mC$ has a very minor effect on the restoration results provided that $R>0$ is sufficiently large. In fact, the proximal alternating minimization seems to converge even without using $\mC$.
\end{remark}

Since our data driven tight frame approach is inspired by the SVD of a low rank $\bmH\left(\La\bsv\right)=\bX\Sig\bY^*$, we expect the learned filters $\left\{\a_1,\a_2,\cdots,\a_{M_2}\right\}$ are related to the annihilating filters. More precisely, we expect $\bY=M_2^{1/2}\left[\begin{array}{cccc}
\a_1&\a_2&\cdots&\a_{M_2}
\end{array}\right]$, so that $M_2^{1/2}\a_{r+1},\cdots,M_2^{1/2}\a_{M_2}$ form an orthonormal basis for the null space of $\bmH\left(\La\bsv\right)$. Hence, following the arguments in \cite{G.Ongie2016}, we can expect that the edge sets can also be restored from the learned filter banks as well. To see this, we let $\varphi_j$ be the trigonometric polynomial defined as
\begin{align*}
\varphi_j(\x)=L^{-1}\sum_{\bk\in\KK}\a_j[-\bk]e^{2\pi i\bk\cdot\x/L},~~~~~j=r+1,\cdots,M_2,
\end{align*}
where the filters are flipped because we use the two-fold Hankel matrix formulation. Given that the conditions in \cite[Theorem 3.4]{G.Ongie2016} are satisfied, we can write $\varphi_j=\eta_j\varphi$ where $\varphi$ is the minimal polynomial in \cref{MajorAssumption} and $\eta_j$ is another $L$-periodic trigonometric polynomial. Following the similar arguments to \cite[Proposition 5.3]{G.Ongie2016}, for $\x\in\left[-L/2,L/2\right)^2$, we have
\begin{align}\label{EdgeEstimate}
\left\|\bV\bV^*\e_{\x}\right\|_2=\left(\sum_{j=r+1}^{M_2}\left|\varphi_j(\x)\right|^2\right)^{1/2}=\left(\sum_{j=r+1}^{M_2}\left|\eta_j(\x)\varphi(\x)\right|^2\right)^{1/2},
\end{align}
where
\begin{align*}
\bV=M_2^{1/2}\left[\begin{array}{ccc}
\a_{r+1}[-\cdot]&\cdots&\a_{M_2}[-\cdot]
\end{array}\right]~~~\text{and}~~~\e_{\x}[\bk]=L^{-1}e^{2\pi i\bk\cdot\x/L},~~~\bk\in\KK.
\end{align*}
Hence, \cref{EdgeEstimate} shows that the learned filter banks are indeed related to the annihilating polynomials in the image domain.

Finally, we would like to mention that the application of dictionary learning to the CS-MRI reconstruction is already explored in \cite{S.Ravishankar2015}. Specifically, the authors in \cite{S.Ravishankar2015} proposed the following (unitary) sparsifying transform learning based CS-MRI (TLMRI) reconstruction model
\begin{align}\label{TLMRI}
\begin{split}
&~~~\min_{\bu,\bsW,\bC}~\mu\left\|\bmR_{\MM}\bmsF\bu-\bsf\right\|_2^2+\sum_{j=1}^{N_p}\left\|\bsW\bmP_j\bu-\bC^{(:,j)}\right\|_2^2\\
&\text{subject to}~~\left\|\bu\right\|_2\leq R,~~\left\|\bC\right\|_0\leq r~~\text{and}~~\bsW^*\bsW=\bsI,
\end{split}
\end{align}
where $\bmP_j$ represents an operator that extracts an $K\times K$ patch of $\bu$ as a vector, $N_p$ denotes the number of patches, and $\bsW$ denotes a square sparsifying transform for the patches of $\bu$. Indeed, it seems that our proposed model \cref{ProposedCSMRIModel} can also be classified into the dictionary learning based approach due to the data driven tight frame regularization. However, unlike the above TLMRI model \cref{TLMRI} which considers the sparse structure of (the patches of) the underlying discrete image, the sparse regularization via data driven tight frame in \cref{ProposedCSMRIModel} is a relaxation of the low rank two-fold Hankel matrix completion in the frequency domain, reflecting the piecewise constant function in a \emph{continuous domain} via annihilating filters. Hence, \cref{TLMRI} is an \emph{on-the-grid} approach, whereas \cref{ProposedCSMRIModel} is an \emph{off-the-grid} approach.

\subsection{Alternating minimization algorithm}\label{PAMAlgorithm}

To solve \cref{ProposedCSMRIModel}, we use the proximal alternating minimization (PAM) algorithm introduced in \cite{H.Attouch2010}. We initialize $\bsv_0=\bmP_{\mC}\left(\bmR_{\MM}\bsf\right)$ where
\begin{align}\label{ProjectOntoC}
\bmP_{\mC}\left(\bsv\right)[\bk]=\min\left\{\left|\bsv[\bk]\right|,R\right\}\exp\left\{i\arg\bsv[\bk]\right\}.
\end{align}
For $\bmW_0$, we can use the SVD of $\bmH\left(\La\bsv_0\right)$. However, to reduce the computational burden, we use the following way. Noting that we in general densely undersample the k-space data on the low frequencies based on the calibrationless sampling strategy, we restrict $\La\bsv_0$ onto the smaller grid, for instance, $\OO'=\left\{-\lfloor N/4\rfloor,\cdots,\lfloor N/4-1\rfloor\right\}^2$. Then we compute the SVD $\bmH\left(\bmR_{\OO'}\La\bsv_0\right)=\bX_0\Sig_0\bY_0^*$, and set $\bmW_0$ via
\begin{align}\label{InitialFilterBanks}
\a_{0,j}=M_2^{-1/2}\bY_0^{(:,j)}~~~~~j=1,\cdots,M_2.
\end{align}
Using the initial filter banks as in \cref{InitialFilterBanks}, we initialize $\bc_0$ via
\begin{align}\label{Initialize}
\bc_0=\left[\begin{array}{cccccc}
\left[\bmS_{\a_{0,1}[-\cdot]}\left(\La\bsv_0\right)\right]^T&\cdots&\left[\bmS_{\a_{0,r}[-\cdot]}\left(\La\bsv_0\right)\right]^T&\0&\cdots&\0
\end{array}\right]^T,
\end{align}
where $r$ is an estimated rank.  After the initializations, we optimize $(\bsv,\bc,\bmW)$ alternatively, as summarized in \cref{Alg1}.

\begin{algorithm}[tp!]
\caption{Proximal Alternating Minimization Algorithm for \cref{ProposedCSMRIModel}}\label{Alg1}
\begin{algorithmic}
\STATE{\textbf{Initialization:} $\bsv_0$, $\bc_0$, $\bmW_0$}
\FOR{$n=0$, $1$, $2$, $\cdots$}
\STATE{\textbf{(1)} Optimize $\bsv$:
\begin{align}\label{vsubprob}
\bsv_{n+1}=\argmin_{\bsv\in\mC}~\f{1}{2}\left\|\bmR_{\MM}\bsv-\bsf\right\|_2^2+\f{\mu}{2}\left\|\bmW_n\left(\La\bsv\right)-\bc_n\right\|_2^2+\f{\beta_1}{2}\left\|\bsv-\bsv_n\right\|_2^2.
\end{align}
\textbf{(2)} Optimize $\bc$:
\begin{align}\label{csubprob}
\bc_{n+1}=\argmin_{\bc}~\gamma\left\|\bc\right\|_0+\f{\mu}{2}\left\|\bc-\bmW_{n}\left(\La\bsv_{n+1}\right)\right\|_2^2+\f{\beta_2}{2}\left\|\bc-\bc_n\right\|_2^2.
\end{align}
\textbf{(3)} Optimize $\bmW$:
\begin{align}\label{Wsubprob}
\bmW_{n+1}=\argmin_{\bmW^*\bmW=\bmI}~\f{\mu}{2}\left\|\bmW\left(\La\bsv_{n+1}\right)-\bc_{n+1}\right\|_2^2+\f{\beta_3}{2}\left\|\bmW-\bmW_n\right\|_F^2.
\end{align}}
\ENDFOR
\end{algorithmic}
\end{algorithm}

It is easy to see that each subproblem in \cref{Alg1} has a closed form solution. Since $\bmW_n^*\bmW_n=\bmI$ for all $n$, the solution to \cref{vsubprob} is given by
\begin{align}\label{vsubexplicit}
\begin{split}
\bsv_{n+1/2}&=\left(\bmR_{\MM}^T\bmR_{\MM}+\mu\La^*\La+\beta_1\bmI\right)^{-1}\left(\bmR_{\MM}^T\bsf+\mu\La^*\bmW_n^*\bc_n+\beta_1\bsv_n\right),\\
\bsv_{n+1}&=\bmP_{\mC}\left(\bsv_{n+1/2}\right),
\end{split}
\end{align}
where $\bmR_{\MM}^T$ is the adjoint of $\bmR_{\MM}$ (i.e. the zero-filling operator), and $\bmP_{\mC}$ is defined as \cref{ProjectOntoC}. It is worth noting that since $\bmR_{\MM}+\mu\La^*\La+\beta_1\bmI$ is a diagonal matrix, no matrix inversion is needed.

For the remaining subproblems \cref{csubprob,Wsubprob}, we note that while the coefficients in \cref{DDTFModel} are intermediate variables to construct tight frame filter banks, the coefficients in \cref{ProposedCSMRIModel} are used to update $\bsv$. Hence, we solve \cref{csubprob,Wsubprob} in the following way. The solution to \cref{csubprob} is expressed as
\begin{align}\label{csubexplicit}
\bc_{n+1}=\bmT_{\sqrt{2\gamma/(\mu+\beta_2)}}\left(\f{\mu\bmW_n\left(\La\bsv_{n+1}\right)+\beta_2\bc_n}{\mu+\beta_2}\right),
\end{align}
where the hard thresholding $\bmT$ for $\bc=\left(\bc_1,\bc_2\right)$ is defined as
\begin{align}\label{HardThresholding}
\bmT_{\gamma}\left(\bc\right)=\left(\bmT_{\gamma}\left(\bc_1\right),\bmT_{\gamma}\left(\bc_2\right)\right),~~~\text{where}~~~\bmT_{\gamma}\left(\bc_j\right)[\bk]=\left\{\begin{array}{cl}
\bc_j[\bk]~&\text{if}~|\bc_j[\bk]|>\lambda\\
0~&\text{otherwise.}
\end{array}\right.
\end{align}
To solve \cref{Wsubprob}, let $\bsH_{n+1}=\bmH\left(\La\bsv_{n+1}\right)\in\C^{2M_1\times M_2}$ and let $\bsA\in\C^{M_2\times M_2}$ be a matrix whose column vectors are concatenations of the filters $\a_1,\cdots,\a_{M_2}$. Then we have
\begin{align}
\bsA_{n+1}&=\argmin_{\bsA\bsA^*=M_2^{-1}\bsI}~\f{\mu}{2}\left\|\bsH_{n+1}\bsA-\bC_{n+1}\right\|_F^2+\f{\beta_3}{2}\left\|\bsA-\bsA_n\right\|_F^2,\label{Dupdate}
\end{align}
where $\bC=\left[\begin{array}{cc}
\bC_1^T&\bC_2^T
\end{array}\right]^T\in\C^{2M_1\times M_2}$, and each $\bC_j\in\C^{M_1\times M_2}$ is the corresponding matrix formulation of frame coefficients $\bc_j$. Then the closed form solution to \cref{Dupdate} is expressed as
\begin{align}\label{Wsubexplicit}
\bsA_{n+1}=M_2^{-1/2}\bX_n\bY_n^*~~\text{with}~~\bX_n\Sig_n\bY_n^*=\bsH_{n+1}^*\bC_{n+1}+\mu^{-1}\beta_3\bsA_n~~\text{is the SVD}.
\end{align}
Therefore, we only compute the SVD of $\bsH_{n+1}^*\bC_{n+1}+\mu^{-1}\beta_3\bsA_n\in\C^{M_2\times M_2}$, leading to the computational efficiency over directly minimizing the rank of $\bmH\left(\La\bsv\right)\in\C^{2M_1\times M_2}$.

Note that, based on the frameworks in \cite{H.Attouch2010,C.Bao2016,J.Bolte2014,Kurdyk1998,Lojasiewicz1995,Y.Xu2013}, it can be verified that \cref{Alg1} globally converges to the critical point, as the minimization of a real-valued complex variable objective function in \cref{ProposedCSMRIModel} is based on the Wirtinger calculus (i.e. by identifying $\C\simeq\R^2$).

For the computational cost, we first mention that we do not explicitly formulate $\bsH_{n+1}=\bmH\left(\La\bsv_{n+1}\right)$ at each iteration step. Without loss of generality, we assume the square patch, i.e. $K_1=K_2=K$, $M_1=(N-K+1)^2$ and $M_2=K^2$. Then $\bmW_n\left(\La\bsv_{n+1}\right)$ and $\bsH_{n+1}^*\bC_{n+1}$ are computed by using $4K^2$ ($2K^2$ for $\bmW_n\left(\La\bsv_{n+1}\right)$ and $2K^2$ for $\bsH_{n+1}^*\bC_{n+1}$) two dimensional convolutions/fast Hankel matrix-vector multiplications \cite{L.Lu2015} directly from $\La\bsv_{n+1}$, requiring $O(K^2N^2\log N)$ operations. In addition, since the hard thresholding $\bmT$ is an elementwise operation, we can update $\bsH_{n+1}^*\bC_{n+1}^{(:,j)}+\mu^{-1}\beta_3\bsA_n^{(:,j)}$ directly after updating the $j$th element of $\bc_{n+1}$. Hence, together with the SVD of $\bsH_{n+1}^*\bC_{n+1}+\mu^{-1}\beta_3\bsA_n\in\C^{K^2\times K^2}$ requiring $O(K^6)$ operations, we need $O(K^2N^2\log N+K^2N^2+K^6)$ operations for \cref{csubexplicit,Wsubexplicit}. For \cref{vsubexplicit}, note that $\bmW_n^*\bc_n$ can be computed by using $2K^2$ two dimensional convolutions, requiring $O(K^2N^2\log N)$ operations. Then since the remaining operations are all pointwise multiplications and additions, \cref{vsubexplicit} requires $O(K^2N^2\log N+N^2)$ operations. Combining these all together, the computational cost of \cref{Alg1} is $O(K^2N^2\log N+K^2N^2+K^6)$.

We further mention that our PAM algorithm seems to be similar to the block coordinate descent (BCD) algorithm presented in \cite{S.Ravishankar2015} to solve the TLMRI model \cref{TLMRI}. In fact, when we set $\beta_1=\beta_2=\beta_3=0$, \cref{Alg1} indeed reduces to the BCD algorithm. In addition, $\bsW$ subproblem in \cref{TLMRI} can be solved similarly to \cref{Wsubexplicit}. Nevertheless, there are differences between \cref{Alg1} and the BCD algorithm in \cite{S.Ravishankar2015}. First of all, since \cref{ProposedCSMRIModel} aims to restore fully sampled k-space data, we do not explicitly use fast Fourier transform (FFT) during the iteration, unlike \cref{TLMRI} which requires to use FFT several times in updating the image $\bu$. Second, as previously mentioned, it is not necessary to explicitly formulate the two-fold Hankel matrix for \cref{csubexplicit,Wsubexplicit}, while $\bC$ subproblem and $\bsW$ subproblem in \cref{TLMRI}, i.e. the so-called \emph{sparse coding} subproblem requires to extract each patch of the image and store it as a column vector. Finally, even though the subsequence convergence of the BCD algorithm for \cref{TLMRI} is established in the works including \cite{S.Ravishankar2015,B.Wen2017}, this does not necessarily imply the global convergence of the iterates. Meanwhile, as previously mentioned, \cref{Alg1} for \cref{ProposedCSMRIModel} guarantees the global convergence of the iterates. Even though the global convergence is not crucial to the applications in image restoration in general, we propose the PAM algorithm with a global convergence guarantee, to balance both theoretical interest and the needs from applications.

\section{Experimental results}\label{Experiments}

In this section, we present the experimental results on the phantom image and the real MR image used in \cite{G.Ongie2016}, to compare the proposed DDTF based CS-MRI model \cref{ProposedCSMRIModel} with several existing methods. We choose to compare with the total variation (TV) model \cite{M.Lustig2007}
\begin{align}\label{TVModel}
\min_{\bu}~\f{1}{2}\left\|\bmR_{\MM}\bmsF\bu-\bsf\right\|_2^2+\gamma\left\|\Na\bu\right\|_1,
\end{align}
the Haar framelet (Haar) model (e.g. \cite{J.F.Cai2012})
\begin{align}\label{HaarModel}
\min_{\bu}~\f{1}{2}\left\|\bmR_{\MM}\bmsF\bu-\bsf\right\|_2^2+\gamma\left\|\bmW\bu\right\|_1,
\end{align}
and the TLMRI model \cref{TLMRI} where a learned transform $\bsW$ is a unitary flipping and rotation invariant sparsifying transform \cite{B.Wen2017}. Both \cref{TVModel,HaarModel} are solved by the split Bregman method \cite{J.F.Cai2009/10,T.Goldstein2009}, and \cref{TLMRI} is solved by the BCD algorithm in \cite{S.Ravishankar2015,B.Wen2017}. We also compare with the following Schatten $p$-norm minimization model
\begin{align}\label{Schatten}
\min_{\bsv\in\mC}~\left\|\bmR_{\MM}\bsv-\bsf\right\|_2^2+\gamma\left\|\bmH\left(\La\bsv\right)\right\|_p^p,
\end{align}
with $p=0$, $0.5$, and $1$, solved by the generic iterative reweighted annihilating filters (GIRAF) method \cite{G.Ongie2017} based on the split Bregman algorithm, which are referred to as ``GIRAF$0$'', ``GIRAF$0.5$'', and ``GIRAF$1$'', respectively. All experiments are implemented on MATLAB $\mathrm{R}2014\mathrm{a}$ running on a laptop with $64\mathrm{GB}$ RAM and Intel(R) Core(TM) CPU $\mathrm{i}7$-$8750\mathrm{H}$ at $2.20\mathrm{GHz}$ with $6$ cores.

\begin{table}[tp!]
\centering
\caption{Parameter selection for each dataset. For TV \cref{TVModel}, Haar \cref{HaarModel}, and GIRAF \cref{Schatten}, $\mu$ refers to the internal parameter for the split Bregman algorithm.}\label{TableParameter}
\vspace{-0.20cm}
\begin{tabular}{|c|c|c|c|c|c|c|c|c|c|}
\hline
Dataset&Model&$K$&$r$&$\mu$&$\gamma$&$\beta_1$&$\beta_2$&$\beta_3$\\ \hline\hline
\multirow{5}{*}{Phantom}&TV \cref{TVModel}&$\cdot$&$\cdot$&$10$&$0.05$&$\cdot$&$\cdot$&$\cdot$\\ \cline{2-9}
&Haar \cref{HaarModel}&$\cdot$&$\cdot$&$10$&$0.1$&$\cdot$&$\cdot$&$\cdot$\\ \cline{2-9}
&TLMRI \cref{TLMRI}&$8$&$198403$&$15.2588$&$\cdot$&$\cdot$&$\cdot$&$\cdot$\\ \cline{2-9}
&GIRAF \cref{Schatten}&$25$&$\cdot$&$4$&$10^{-6}$&$\cdot$&$\cdot$&$\cdot$\\ \cline{2-9}
&DDTF \cref{ProposedCSMRIModel}&$25$&$500$&$0.1$&$10$&$10^{-4}$&$10^{-4}$&$10^{-4}$\\ \hline\hline
\multirow{5}*{Real MR}&TV \cref{TVModel}&$\cdot$&$\cdot$&$10$&$0.01$&$\cdot$&$\cdot$&$\cdot$\\ \cline{2-9}
&Haar \cref{HaarModel}&$\cdot$&$\cdot$&$10$&$0.025$&$\cdot$&$\cdot$&$\cdot$\\ \cline{2-9}
&TLMRI \cref{TLMRI}&$8$&$396806$&$15.2588$&$\cdot$&$\cdot$&$\cdot$&$\cdot$\\ \cline{2-9}
&GIRAF \cref{Schatten}&$45$&$\cdot$&$2$&$10^{-6}$&$\cdot$&$\cdot$&$\cdot$\\ \cline{2-9}
&DDTF \cref{ProposedCSMRIModel}&$45$&$1620$&$0.05$&$5$&$10^{-4}$&$10^{-4}$&$10^{-4}$\\ \hline
\end{tabular}
\end{table}

In \cref{ProposedCSMRIModel,Schatten}, we choose $R=|\bsf[\0]|$ if $\0\in\MM$ and $R=10^8$ otherwise for the constraint set $\mC$ in \cref{ConstraintSet}. In \cref{TLMRI}, we choose $R=10^5$. We use the forward difference for the discrete gradient $\Na$ in \cref{TVModel}, and $\bmW$ in \cref{HaarModel} is chosen to be the undecimated tensor product Haar framelet transform with $1$ level of decomposition \cite{B.Dong2013}. Both \cref{ProposedCSMRIModel,Schatten} use the $K\times K$ square patch to generate the two-fold Hankel matrix for simplicity. We choose $K$ to be no more than $0.2N$, based on \cite{G.Ongie2018}, and we choose $r\approx0.8K^2$ for the initialization \cref{Initialize} of \cref{ProposedCSMRIModel}, both of which depend on the geometry of the target image. For \cref{TLMRI}, we use $8\times 8$ patches. The detailed choice of the remaining regularization parameters are summarized in \cref{TableParameter}. Empirically, we observe that $\mu\approx0.01\gamma$ is an appropriate choice for \cref{ProposedCSMRIModel}, and we observe that when $\mu$ is large, the restored k-space data $\bsv$ has a faster decay than smaller $\mu$. Hence, the parameters are manually tuned so that we can achieve the optimal restoration of both the low frequencies and high frequencies. For the on-the-grid approaches  \cref{TVModel,HaarModel,TLMRI}, the stopping criterion is
\begin{align}\label{StoppingOntheGrid}
\f{\left\|\bu_{n+1}-\bu_n\right\|_2}{\left\|\bu_n\right\|_2}\leq\eps:=2\times10^{-4},
\end{align}
and for the off-the-grid approaches \cref{ProposedCSMRIModel,Schatten}, we use
\begin{align}\label{StoppingOfftheGrid}
\f{\left\|\bsv_{n+1}-\bsv_n\right\|_2}{\left\|\bsv_n\right\|_2}\leq\eps.
\end{align}
We also set the maximum allowable number of iterations to be $600$. To measure the quality of restored images, we compute the signal-to-noise ratio (SNR), the high frequency error norm (HFEN) \cite{S.Ravishankar2011}, the number of iterations (NIters) to reach the stopping criterion, and the CPU time. Note that for the off-the-grid approaches \cref{ProposedCSMRIModel,Schatten}, the restored image is computed via the inverse DFT of the restored k-space data; see \cref{CSMRIComparisons}.

\subsection{Phantom experiments}\label{PhantomExperiments}

\begin{figure}[tp!]
\centering
\hspace{-0.1cm}\subfloat[Fully sampled]{\label{PhantomOriginalk}\includegraphics[width=3.0cm]{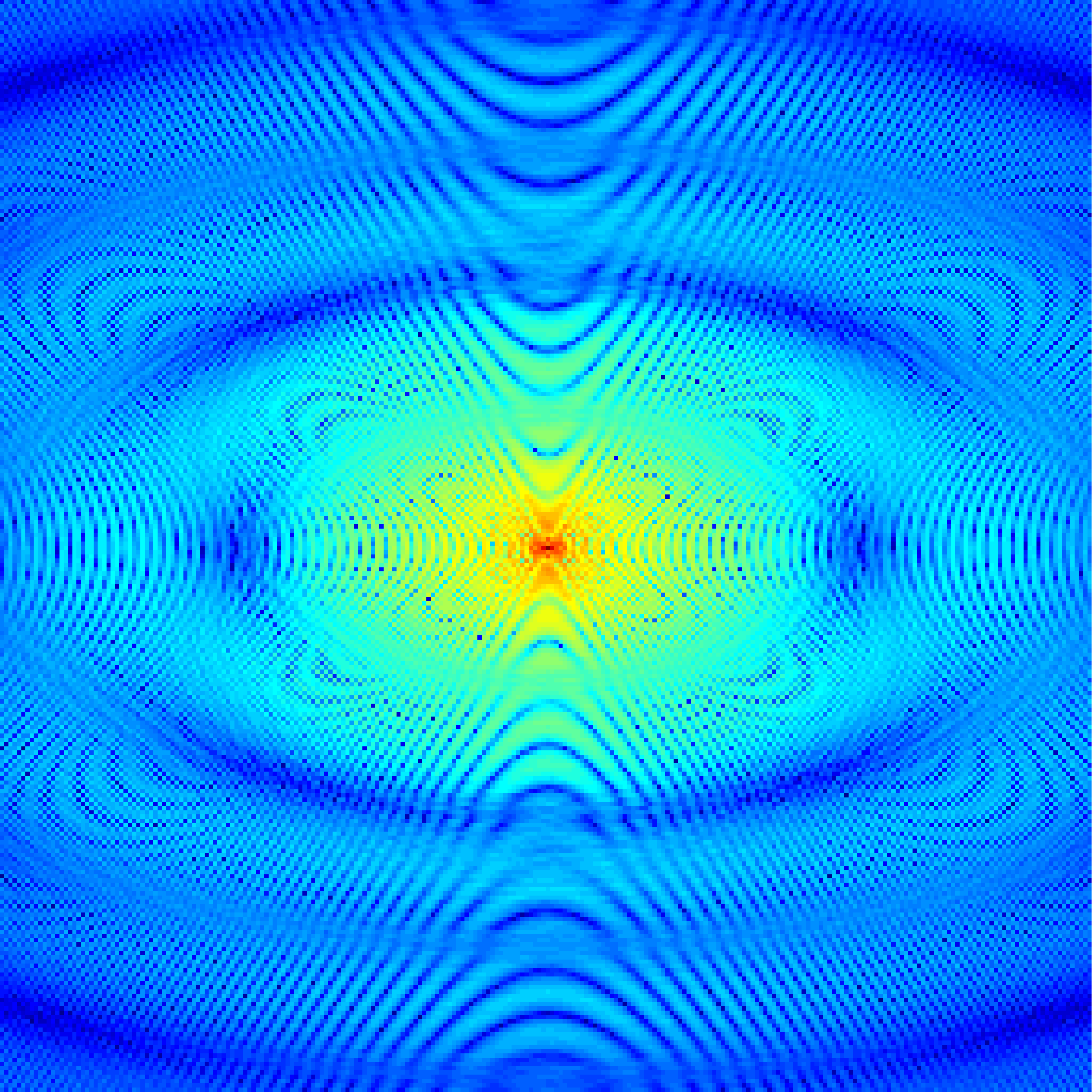}}\hspace{0.005cm}
\subfloat[Ground truth]{\label{PhantomOriginal}\includegraphics[width=3.0cm]{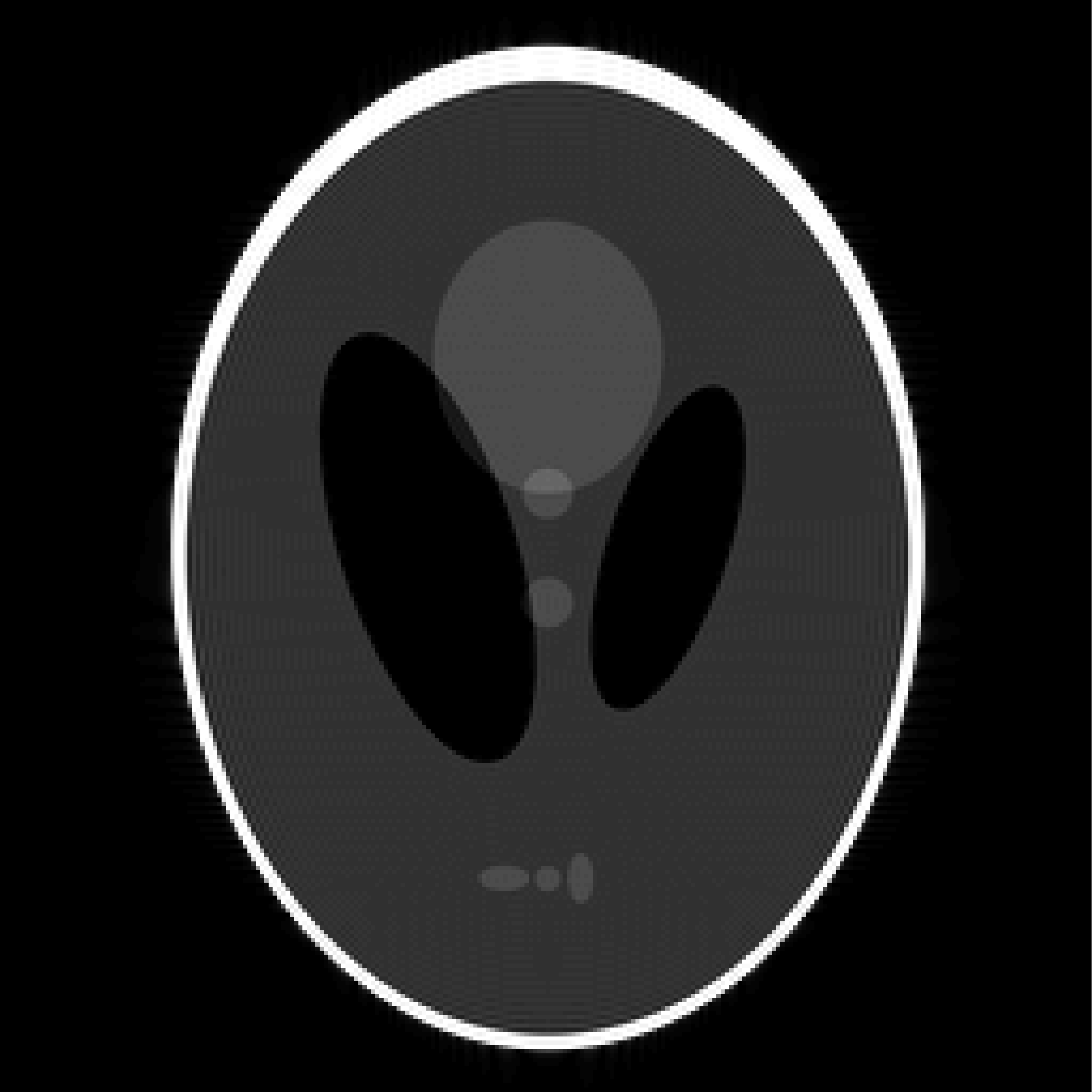}}\hspace{0.005cm}
\subfloat[Sample mask]{\label{PhantomMask}\includegraphics[width=3.0cm]{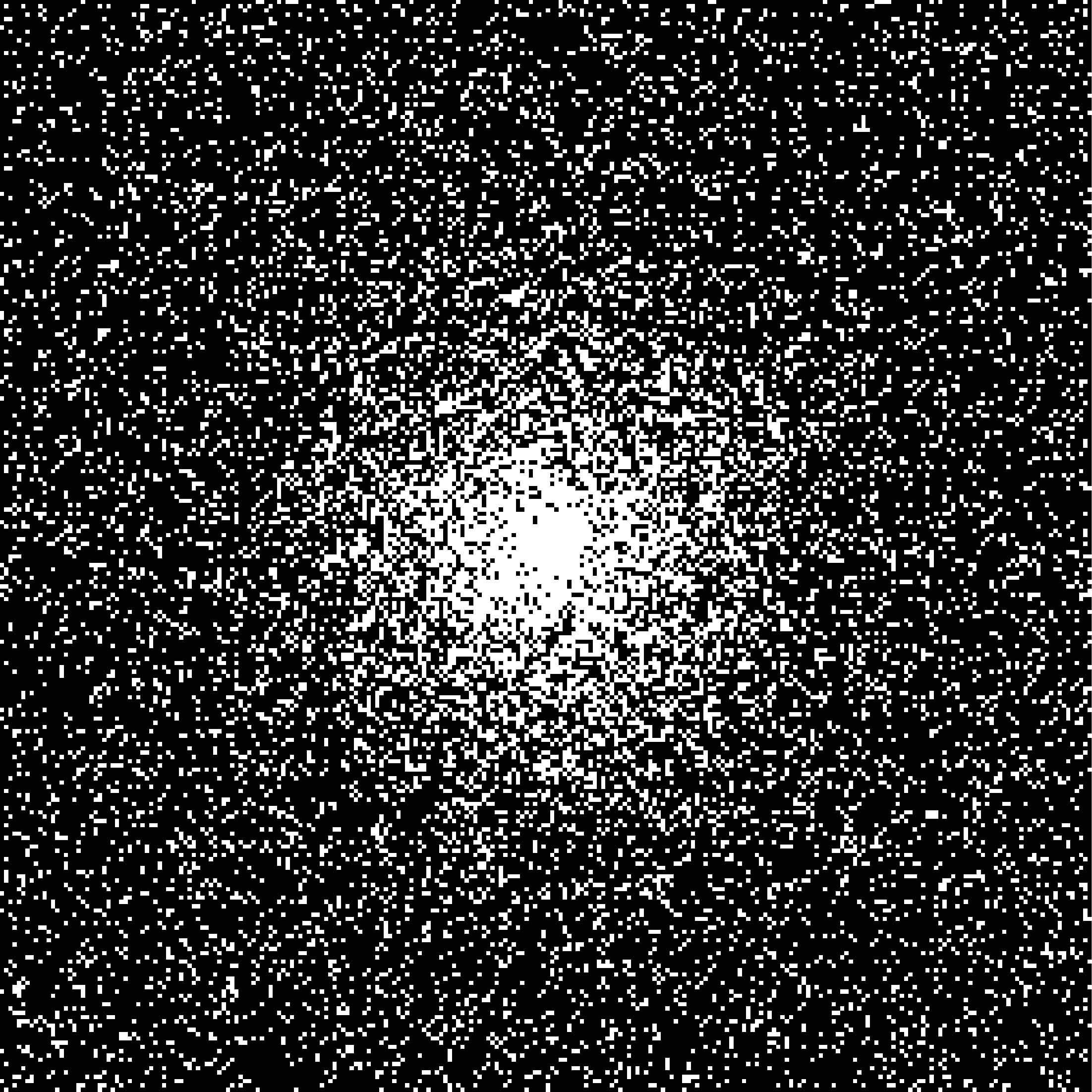}}\hspace{0.005cm}
\subfloat[Undersampled]{\label{PhantomUndersamplek2}\includegraphics[width=3.0cm]{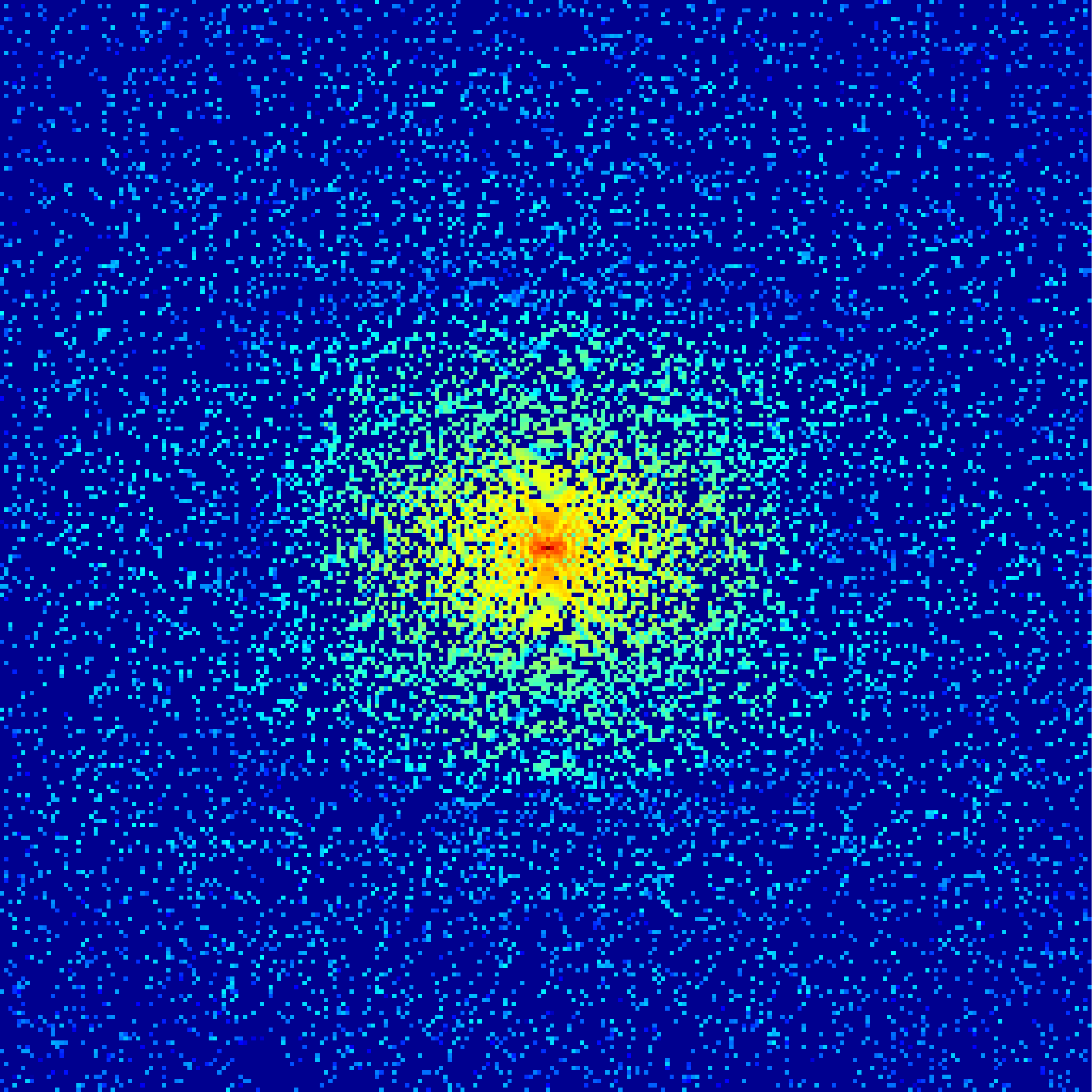}}\vspace{-0.20cm}
\caption{Dataset for the phantom experiments. Fully sampled k-space data, its inverse DFT as a ground truth, the undersampling mask, and the undersampled k-space data.}\label{PhantomDataSet}
\end{figure}

For the piecewise constant phantom experiments, we first compute the fully sampled k-space data using \cref{MRForward} with a piecewise constant function model \cref{uModel}. More precisely, given that $\Om_j$'s are ellipses, \cref{uModel} can be explicitly written as
\begin{align}\label{uModelExp}
u(\x)=\sum_{j=1}^J\alpha_j1_{B(\0,1)}\left(\bsD_j\bQ_j\left(\x-\y_j\right)\right),
\end{align}
where $\bsD_j$'s are diagonal matrices, $\bQ_j$'s are rotation matrices, and $B(\0,1)$ is the unit disk. Following \cite{M.Guerquin-Kern2012}, we can compute the Fourier transform of \cref{uModelExp} by noting that
\begin{align*}
\msF(1_{B(\0,1)})(\xxi)=\f{J_1(2\pi|\xxi|)}{|\xxi|}~~~\text{and}~~~\msF(1_{B(\0,1)})(\0)=\pi,
\end{align*}
where $J_1$ is the first kind Bessel's function of order $1$. Then using the variable density random sampling method in \cite{M.Lustig2007}, we generate $20\%$ undersampled k-space data. The complex white Gaussian noise is also added so that the resulting SNR of the samples is approximately $25\mathrm{dB}$ (See \cref{PhantomDataSet}).

\begin{figure}[tp!]
\centering
\hspace{-0.1cm}\subfloat[$r=550$]{\label{PhantomDDTFEdge551}\includegraphics[width=3.5cm]{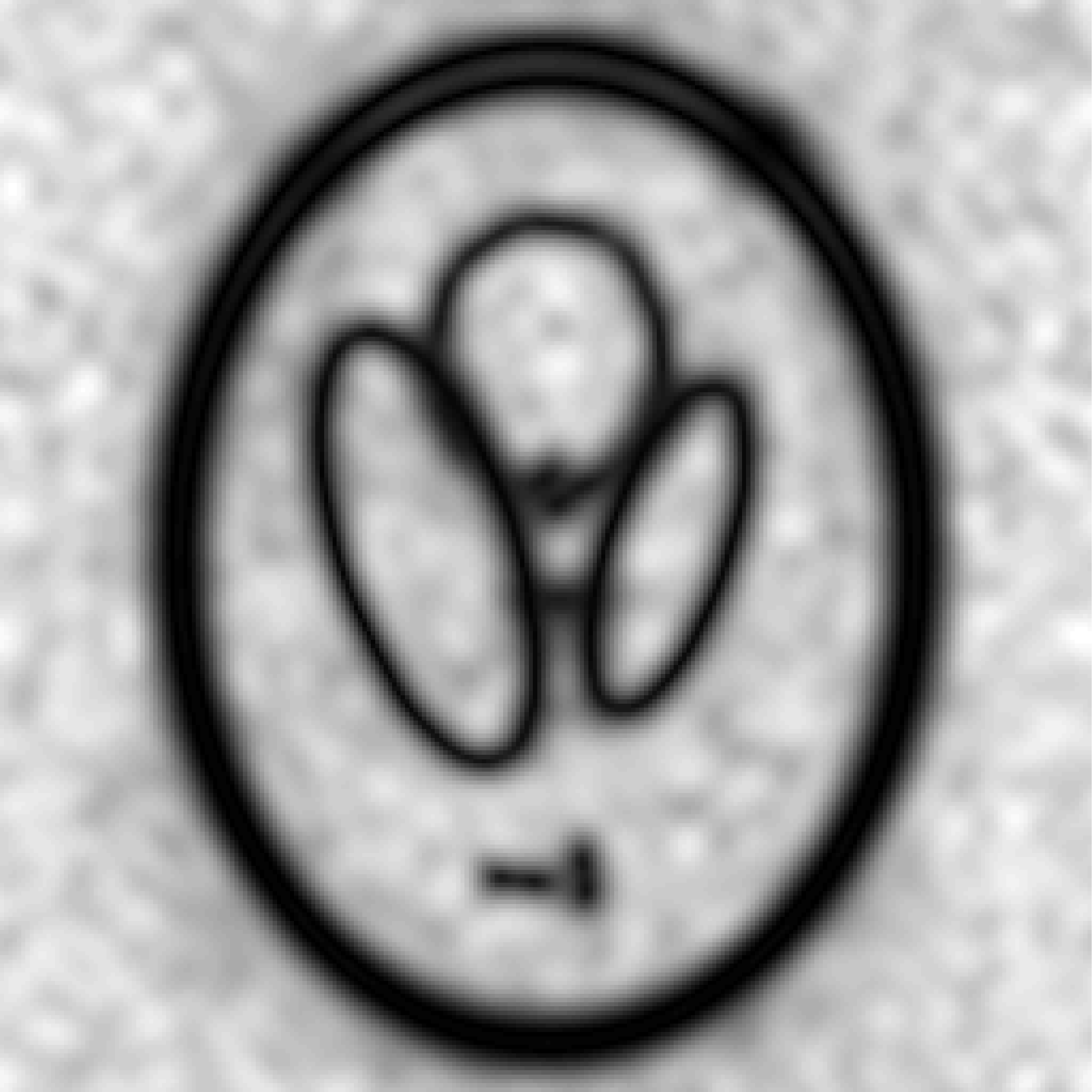}}\hspace{0.005cm}
\subfloat[$r=450$]{\label{PhantomDDTFEdge451}\includegraphics[width=3.5cm]{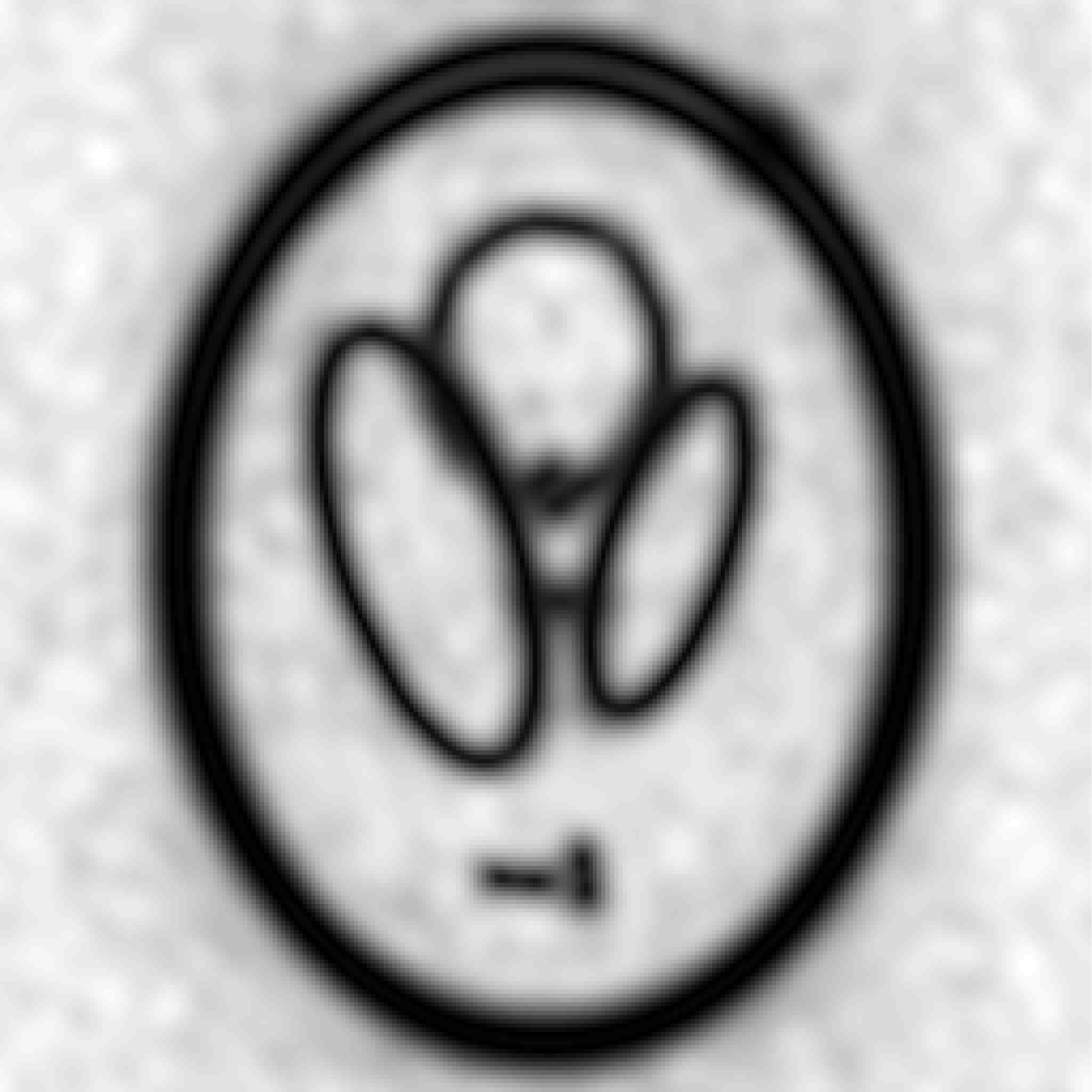}}\hspace{0.005cm}
\subfloat[$r=350$]{\label{PhantomDDTFEdge351}\includegraphics[width=3.5cm]{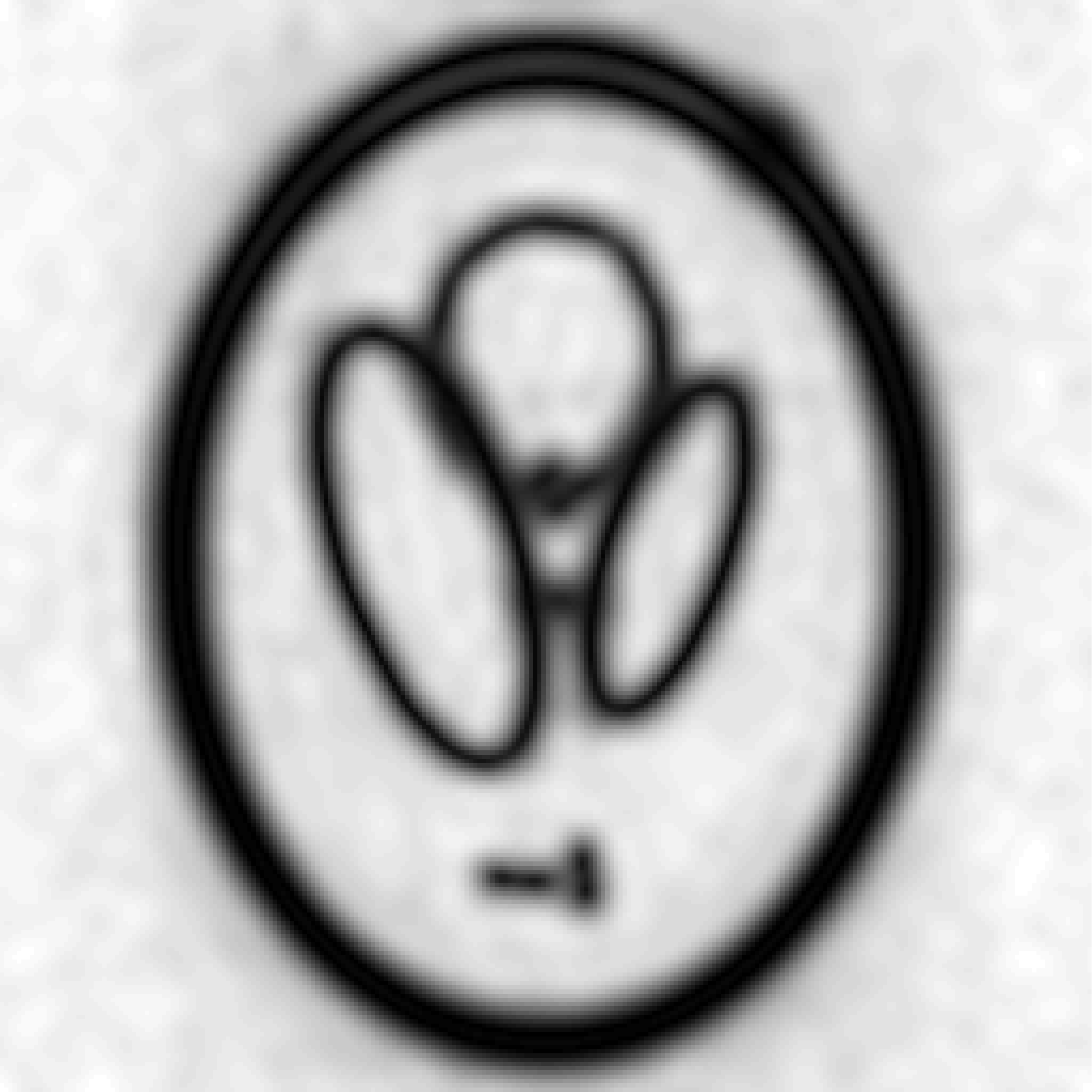}}\hspace{0.005cm}
\subfloat[$r=336$]{\label{PhantomDDTFEdge337}\includegraphics[width=3.5cm]{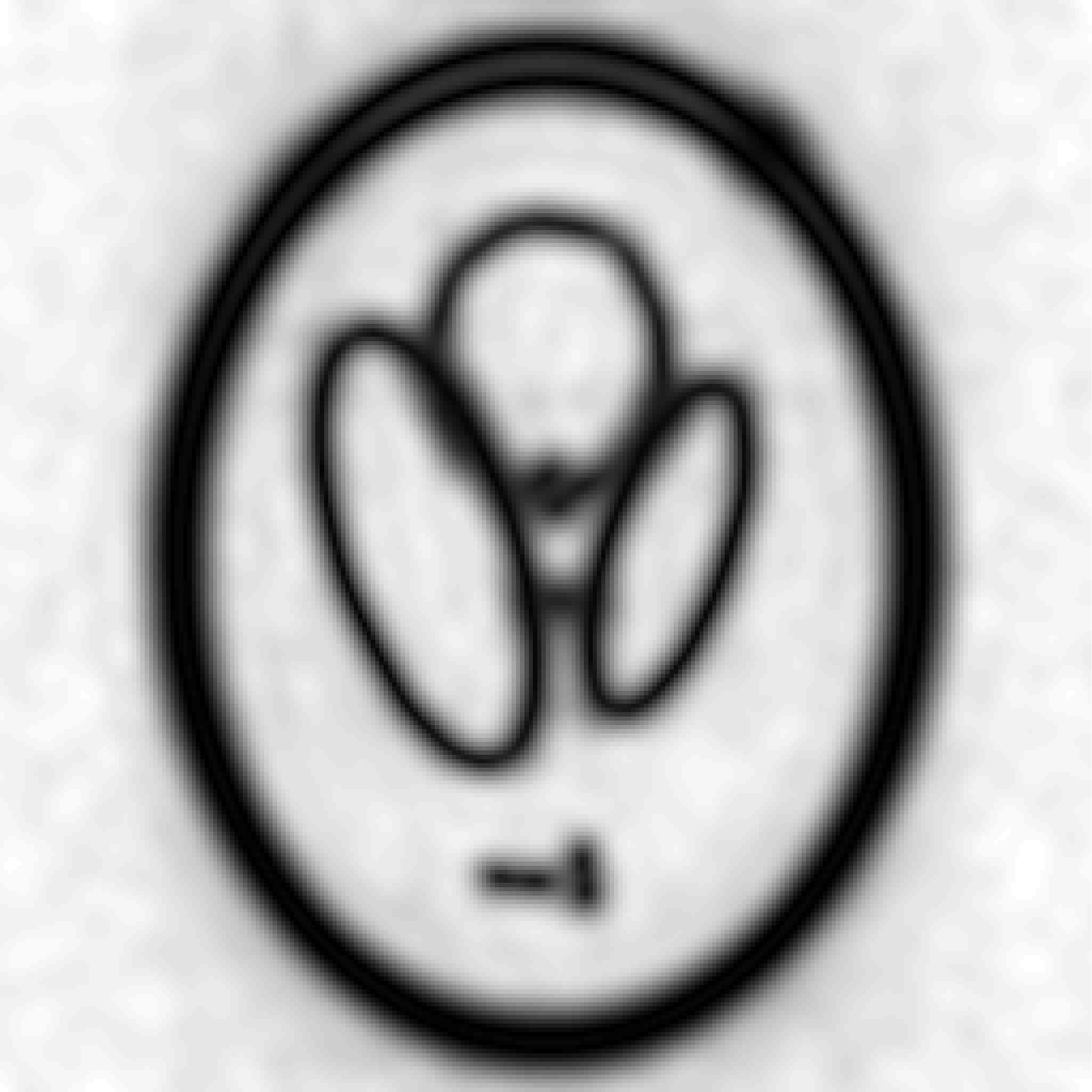}}\vspace{-0.20cm}
\caption{Images of edges estimated by using \cref{EdgeEstimate} with the filter banks in \cref{ProposedCSMRIModel} and various values of $r$. All estimated images are displayed in the window level $[0,1]$.}\label{PhantomEdges}
\end{figure}

\begin{table}[tp!]
\centering
\caption{Comparison of signal-to-noise ratio, high frequency error norm, iteration numbers, and CPU time for the phantom experiments. MaxIter denotes that the iteration reached the maximum number without satisfying the stopping criterion.}\label{PhantomTable}
\vspace{-0.20cm}
\begin{tabular}{|c|c|c|c|c|c|c|c|c|}
\hline
\multirow{2}{*}{Indices}&\multirow{2}{*}{Zero fill}&\multirow{2}{*}{TV \cref{TVModel}}&\multirow{2}{*}{Haar \cref{HaarModel}}&\multirow{2}{*}{TLMRI \cref{TLMRI}}&\multicolumn{3}{|c|}{GIRAF \cref{Schatten}}&\multirow{2}{*}{DDTF \cref{ProposedCSMRIModel}}\\ \cline{6-8}
&&&&&$p=0$&$p=0.5$&$p=1$&\\ \hline\hline
SNR&$9.23$&$20.29$&$20.99$&$21.62$&$21.12$&$20.88$&$17.13$&$\textbf{26.66}$\\ \hline
HFEN&$0.5490$&$0.0992$&$0.0904$&$0.0675$&$0.0948$&$0.0999$&$0.1828$&$\textbf{0.0572}$\\ \hline
NIters&$\cdot$&$210$&$244$&MaxIter&$18$&$20$&$15$&$154$\\ \hline
Time($\mathrm{s}$)&$\cdot$&$2.4$&$13.7$&$32612.9$&$15.0$&$16.2$&$12.2$&$1438.2$\\ \hline
\end{tabular}
\end{table}

\cref{PhantomEdges} displays the edges estimated by the proposed model \cref{ProposedCSMRIModel}. More precisely, using the obtained filter banks $\a_1,\cdots,\a_{K^2}$ at convergence, we compute $\left\|\bV\bV^*\e_{\x}\right\|_2$ defined as \cref{EdgeEstimate}, by varying the rank cutoff $r=550, 450, 350, 336$, respectively. Here, $r=336$ is the ground truth rank cutoff described in \cite{G.Ongie2018}. At first glance, we can see that as the rank cutoff $r$ in \cref{EdgeEstimate} decreases, the image of $\left\|\bV\bV^*\e_{\x}\right\|_2$ contains less artifacts. Most importantly, we can further observe that in any case, $\left\|\bV\bV^*\e_{\x}\right\|_2$ is close to $0$ near the edges. This demonstrates that the data driven tight frame filter banks used in the proposed model \cref{ProposedCSMRIModel} are indeed related to the annihilating filter, as they are also able to estimate the edges of the piecewise constant function via \cref{EdgeEstimate}. Therefore, we conclude that the proposed sparse regularization via data driven tight frame provides another relaxation of the low rank two-fold Hankel matrix completion based on the annihilating filter.

\begin{figure}[tp!]
\centering
\hspace{-0.1cm}\subfloat[Fully sampled]{\label{PhantomOriginal2}\includegraphics[width=3.00cm]{PhantomOriginal.pdf}}\hspace{0.005cm}
\subfloat[Zero fill]{\label{PhantomZeroPad}\includegraphics[width=3.00cm]{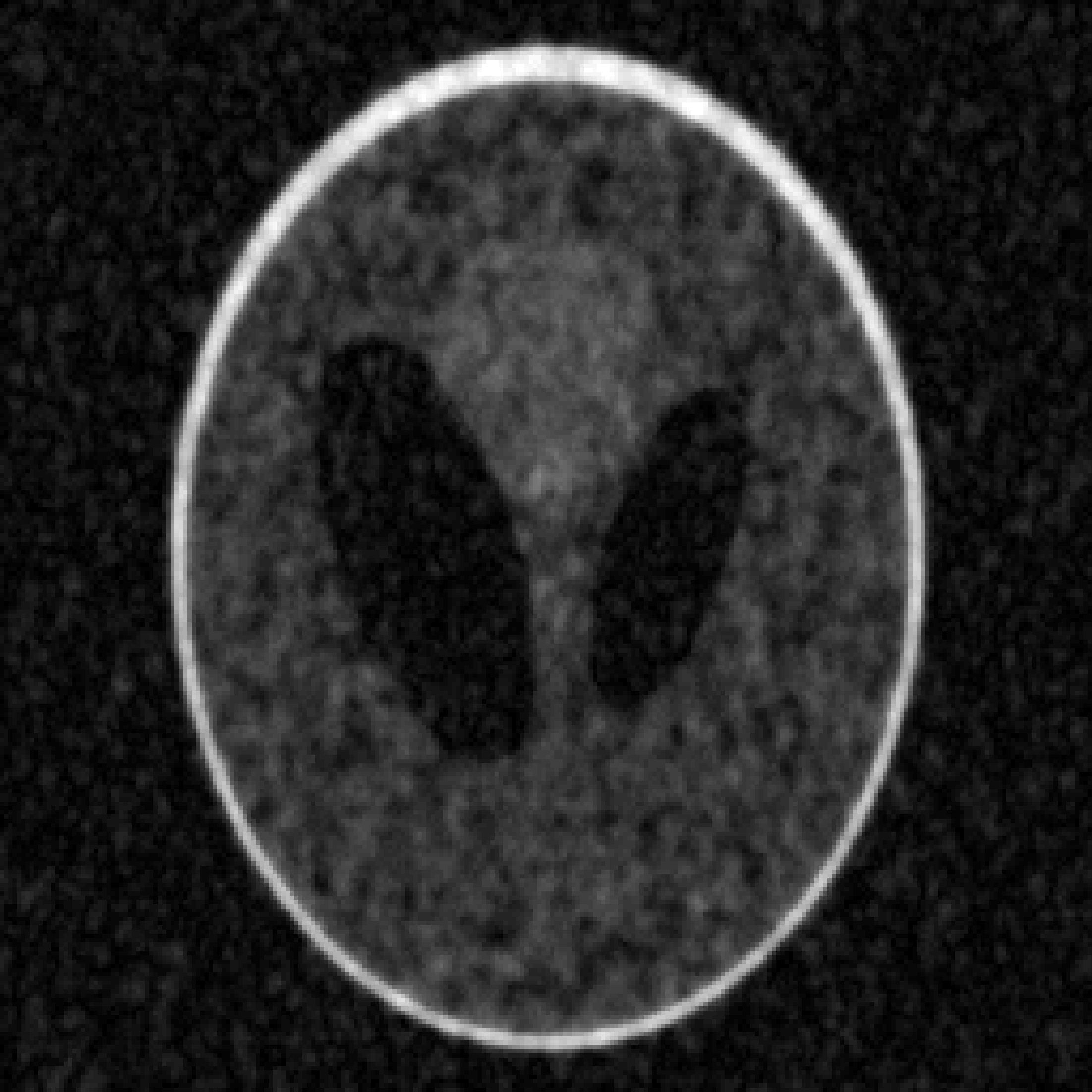}}\hspace{0.005cm}
\subfloat[TV \cref{TVModel}]{\label{PhantomTV}\includegraphics[width=3.00cm]{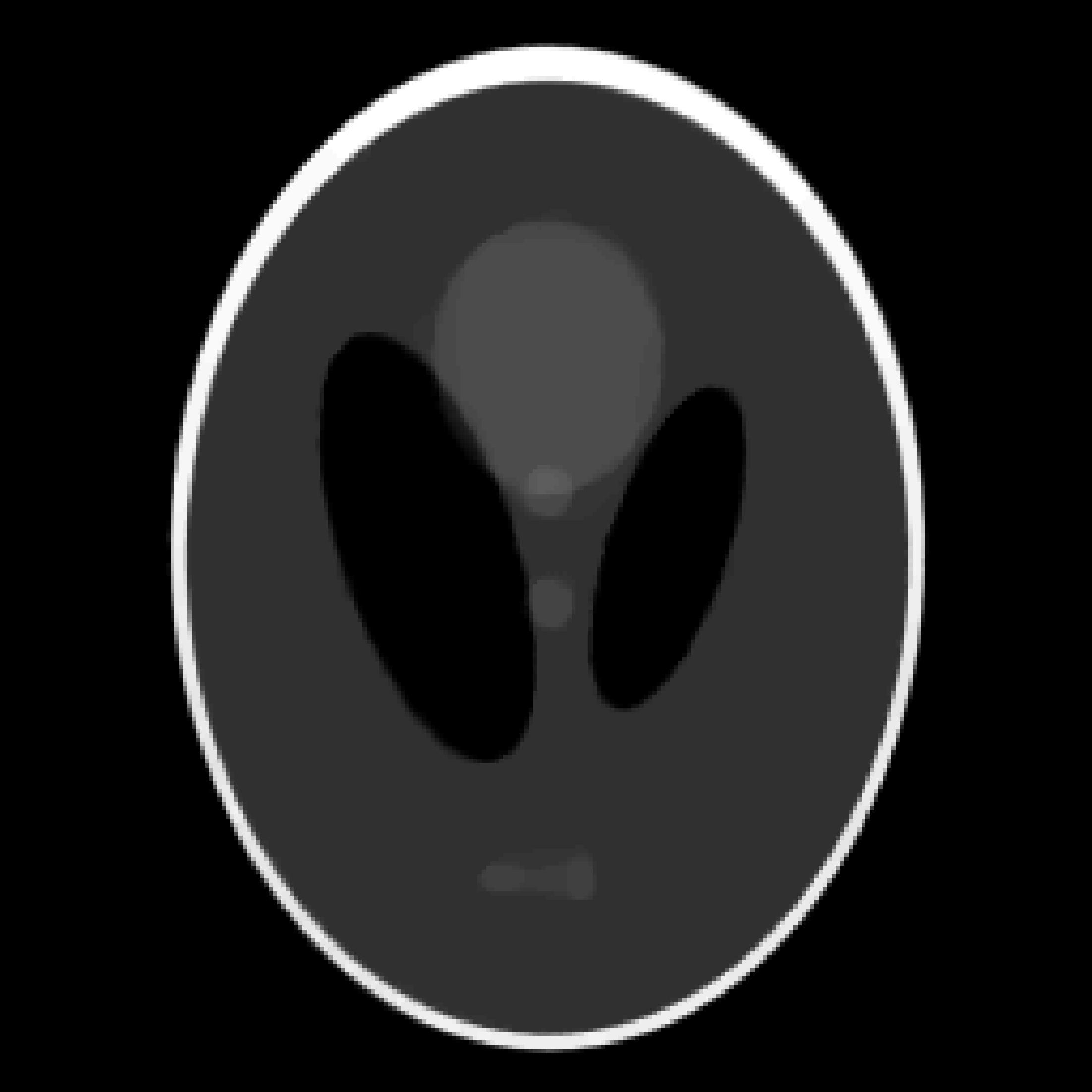}}\hspace{0.005cm}
\subfloat[Haar \cref{HaarModel}]{\label{PhantomHaar}\includegraphics[width=3.00cm]{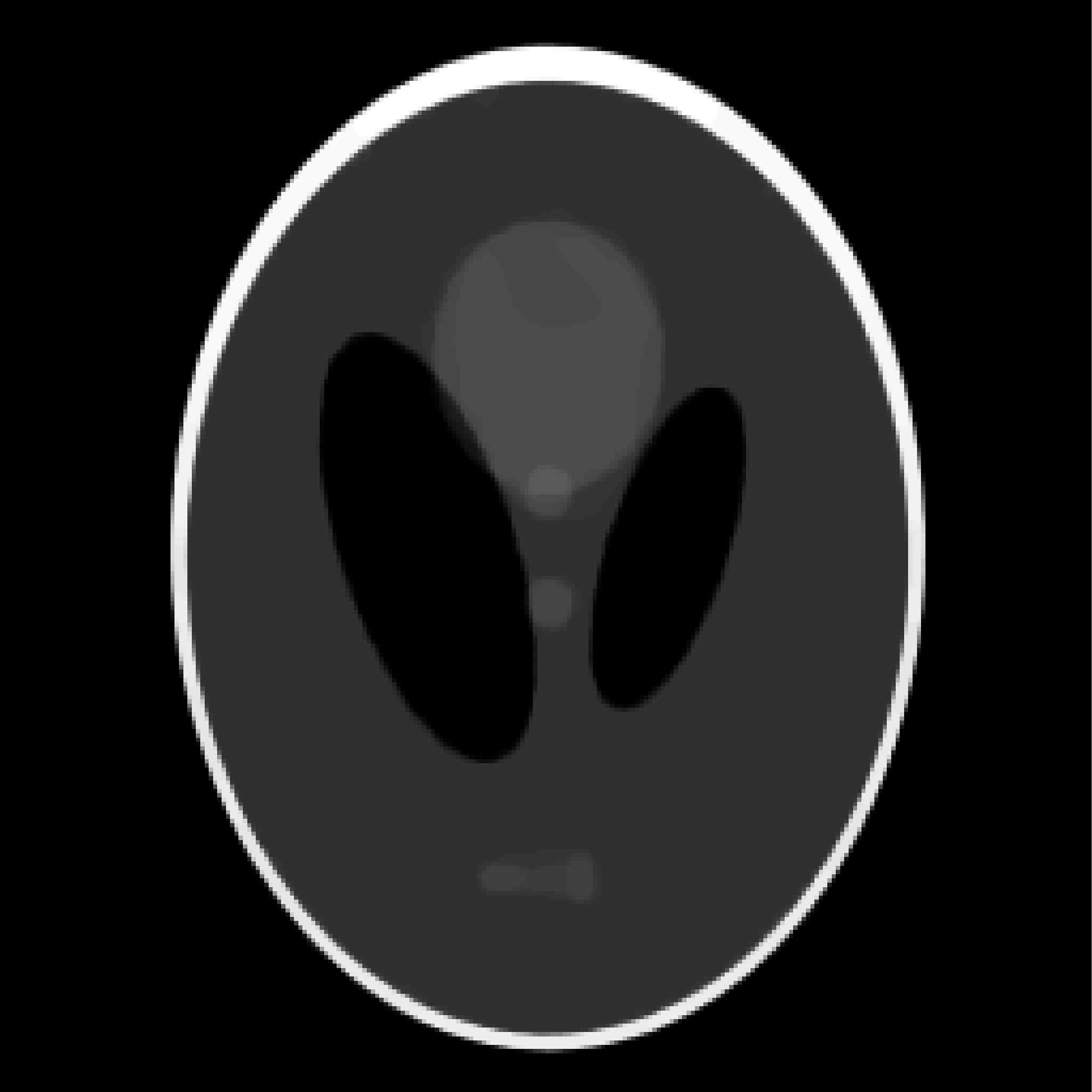}}\hspace{0.005cm}
\subfloat[TLMRI \cref{TLMRI}]{\label{PhantomFRIST}\includegraphics[width=3.00cm]{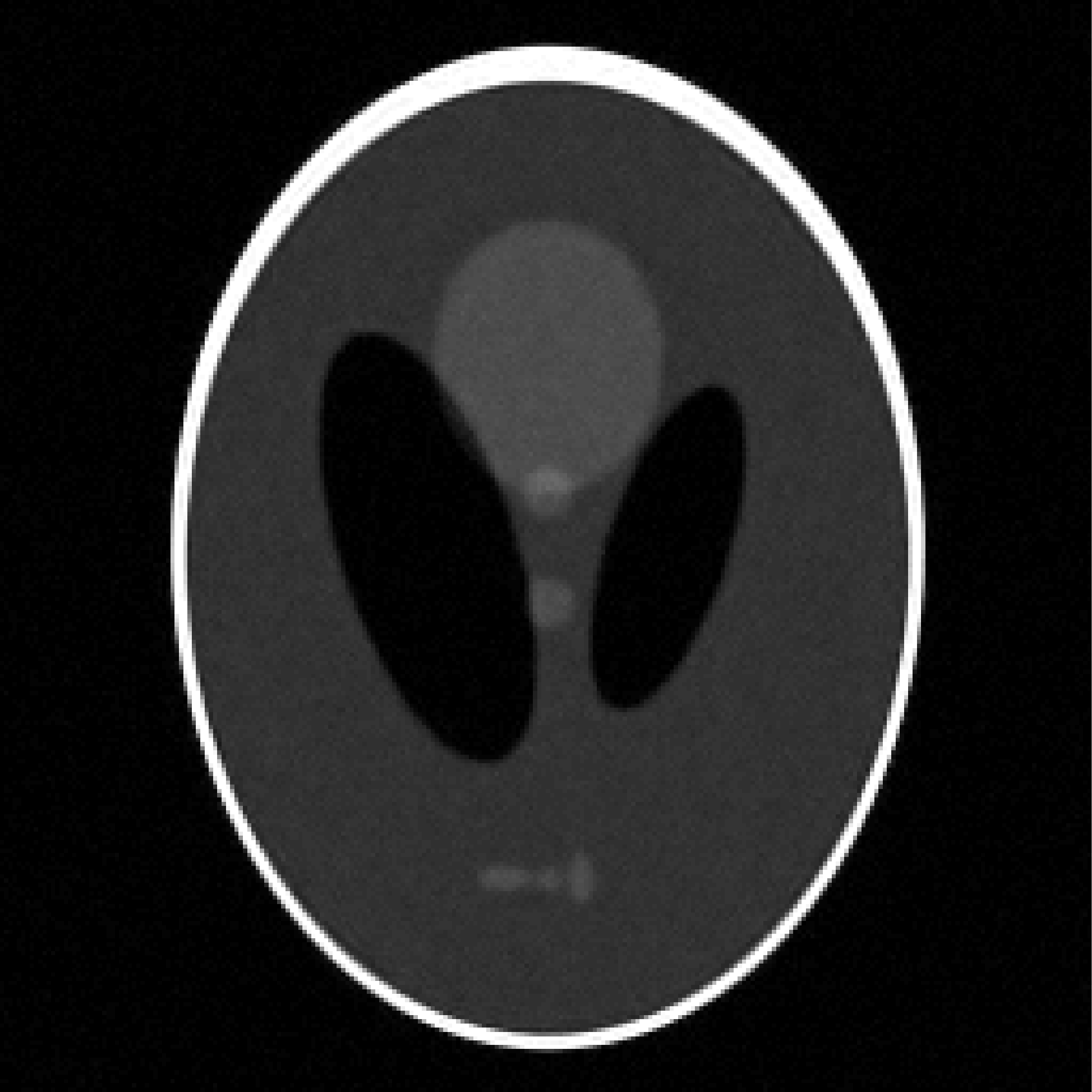}}\vspace{-0.20cm}\\
\subfloat[GIRAF$0$]{\label{PhantomGIRAF0}\includegraphics[width=3.00cm]{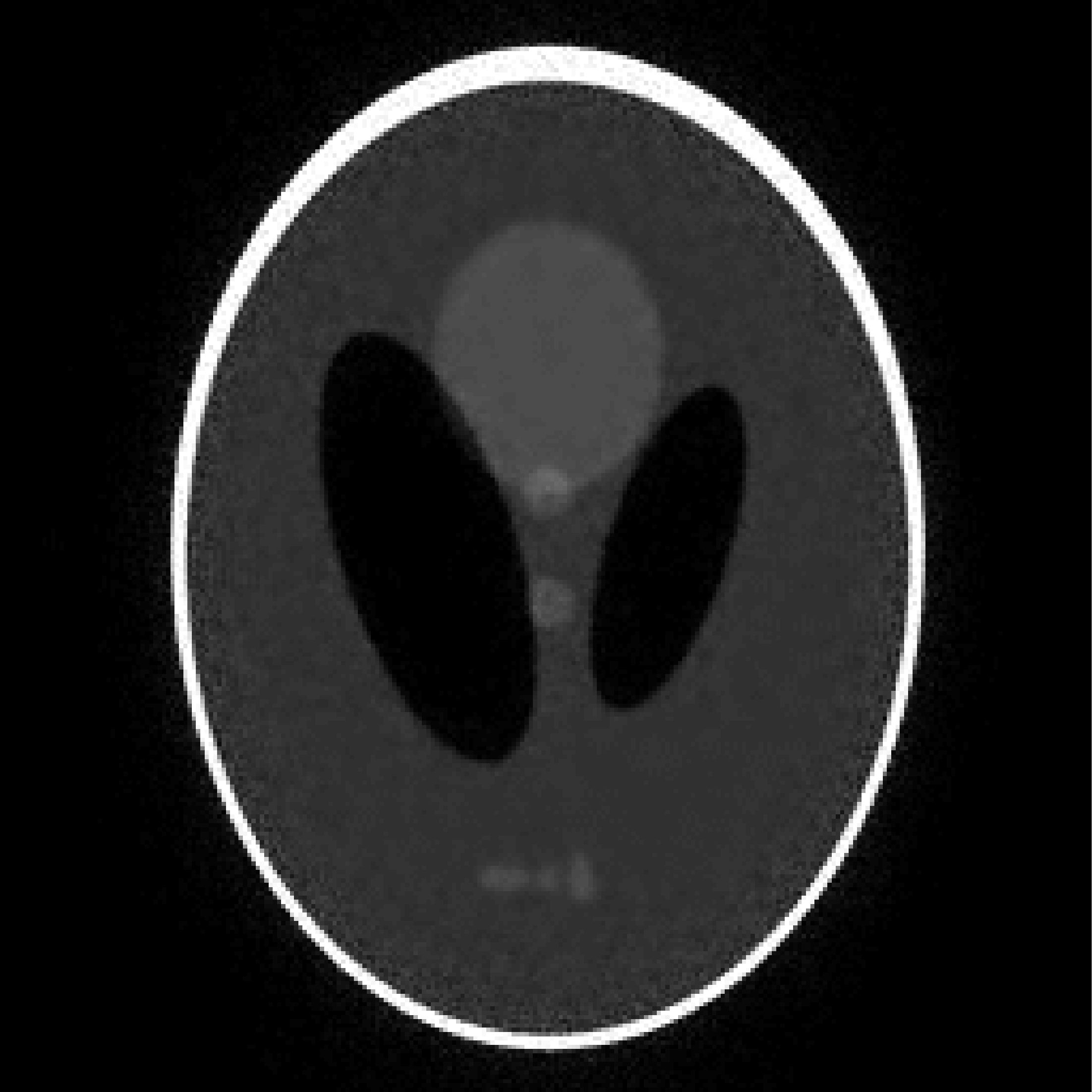}}\hspace{0.005cm}
\subfloat[GIRAF$0.5$]{\label{PhantomGIRAFHalf}\includegraphics[width=3.00cm]{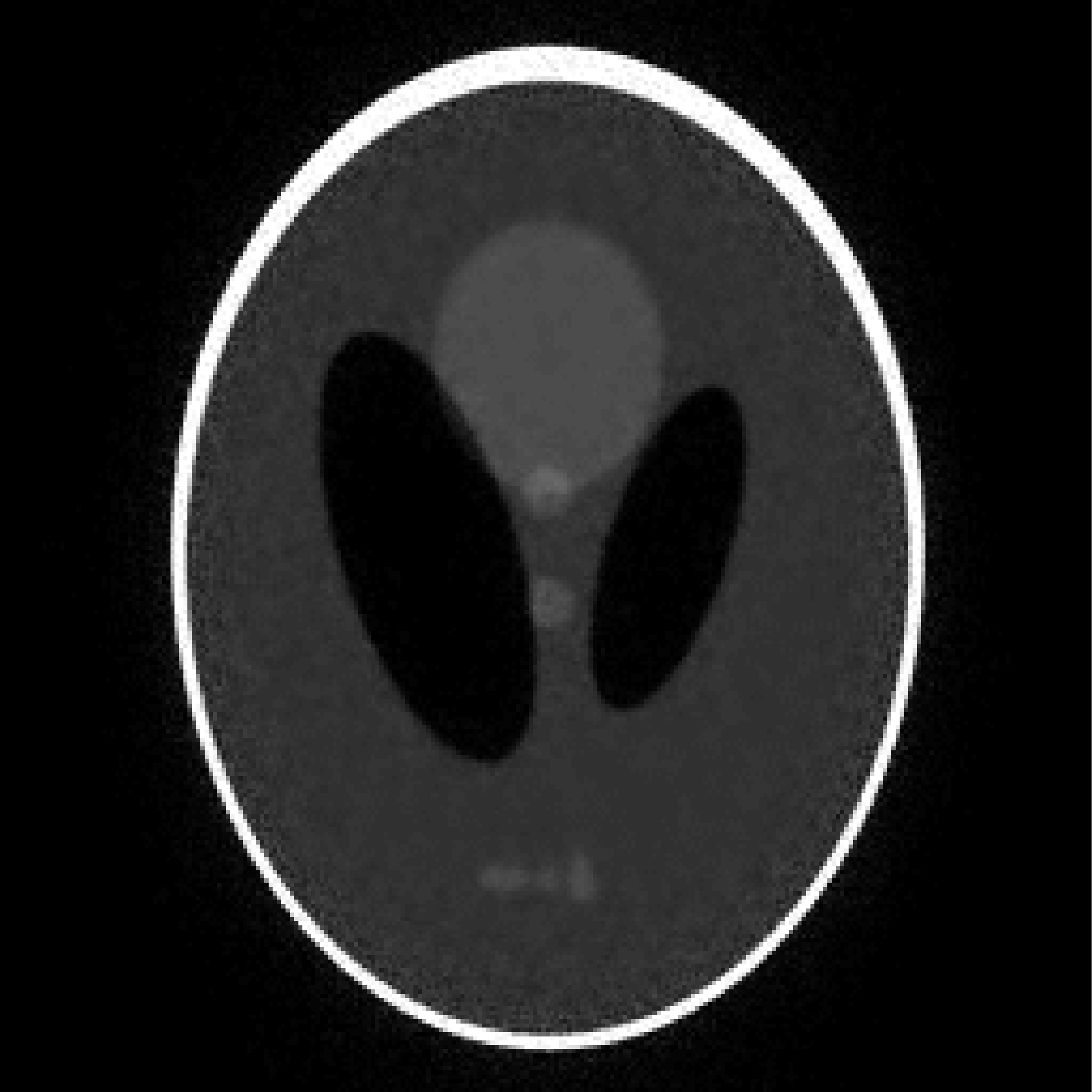}}\hspace{0.005cm}
\subfloat[GIRAF$1$]{\label{PhantomGIRAF1}\includegraphics[width=3.00cm]{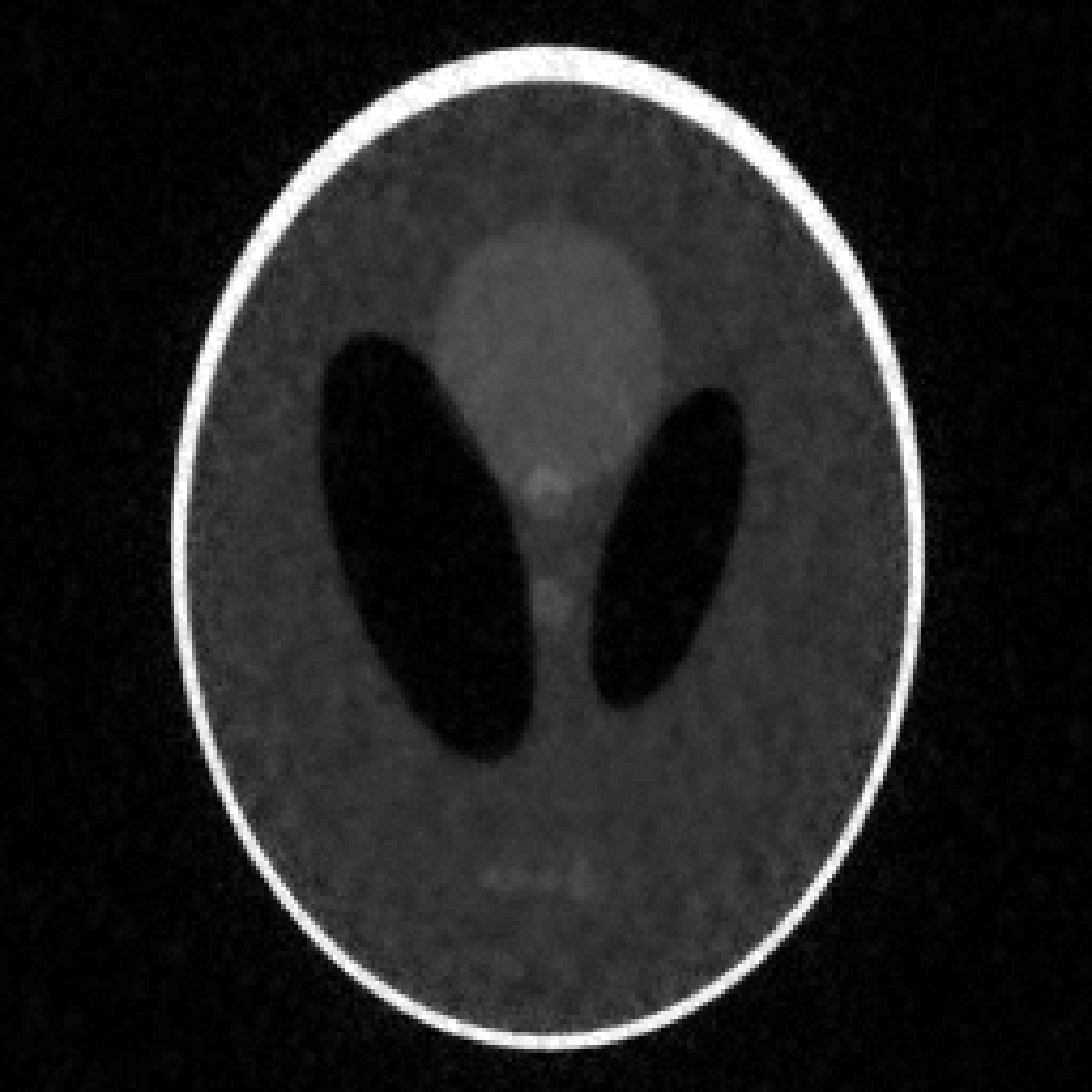}}\hspace{0.005cm}
\subfloat[DDTF \cref{ProposedCSMRIModel}]{\label{PhantomDDTF}\includegraphics[width=3.00cm]{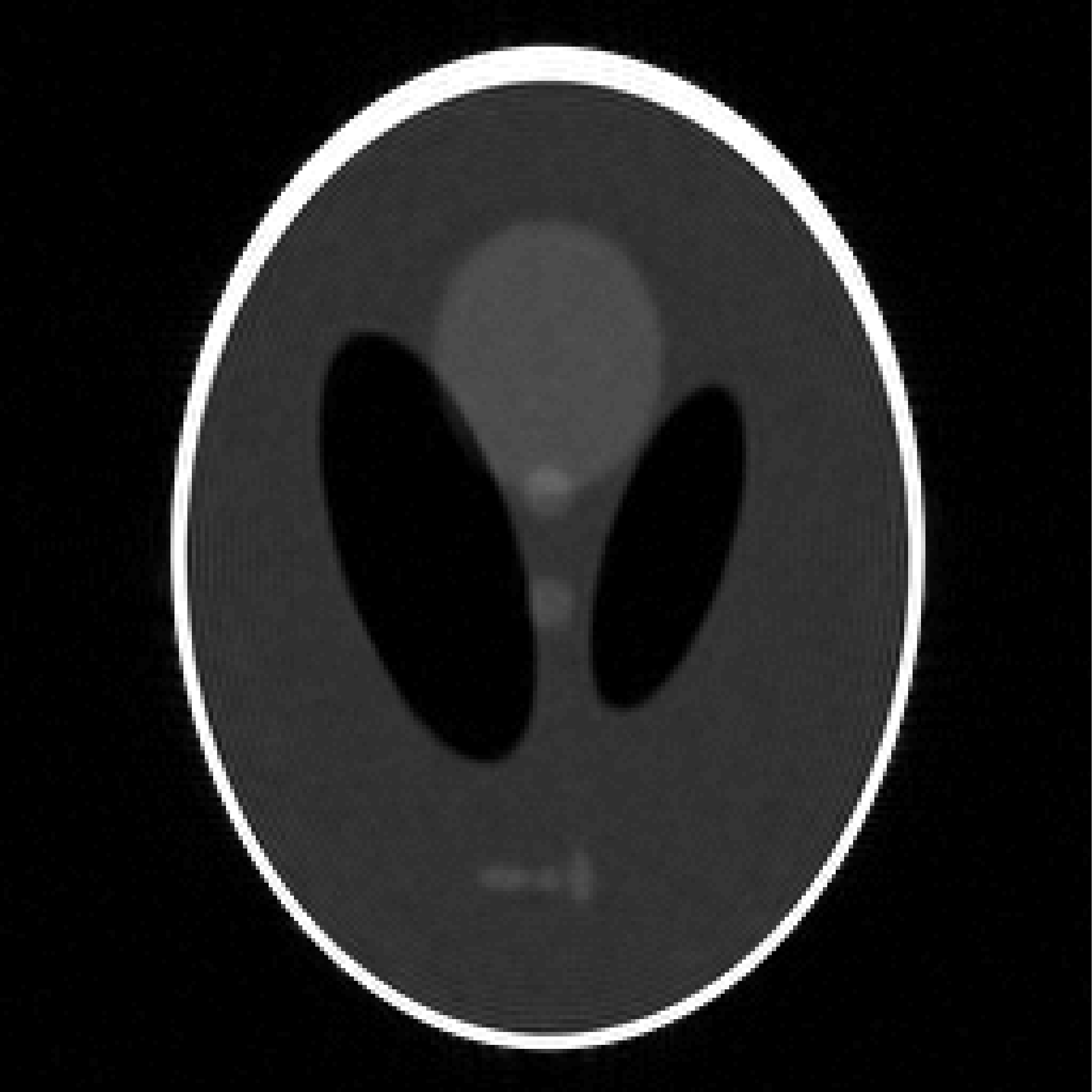}}\vspace{-0.20cm}
\caption{Visual comparisons of each restoration models for the phantom experiments. All restored images are displayed in the window level $[0,1]$ for fair comparisons.}\label{PhantomResults}
\end{figure}

\begin{figure}[tp!]
\centering
\hspace{-0.1cm}\subfloat[Sample mask]{\label{PhantomMask2}\includegraphics[width=3.00cm]{PhantomMask.pdf}}\hspace{0.005cm}
\subfloat[Zero fill]{\label{PhantomZeroPadError}\includegraphics[width=3.00cm]{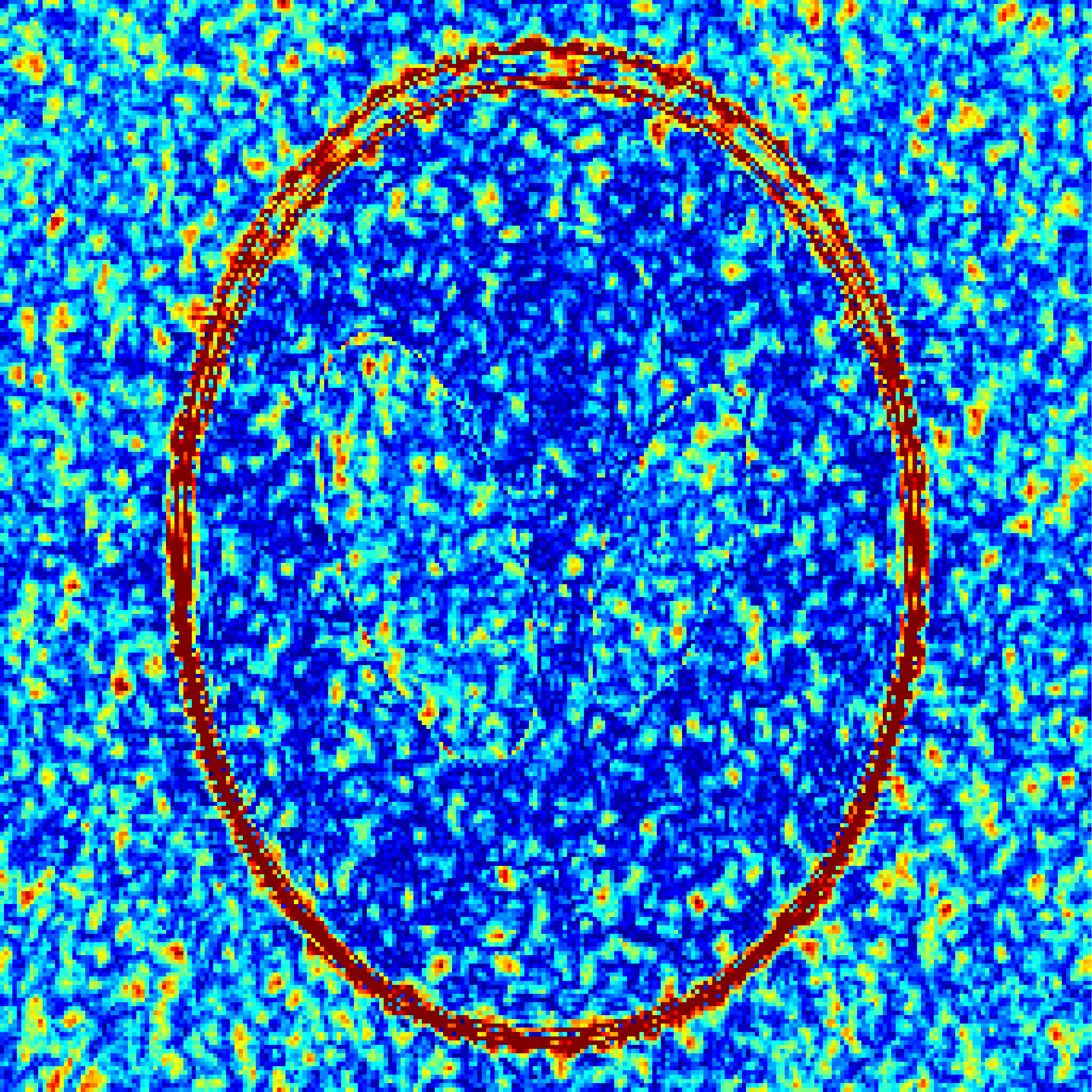}}\hspace{0.005cm}
\subfloat[TV \cref{TVModel}]{\label{PhantomTVError}\includegraphics[width=3.00cm]{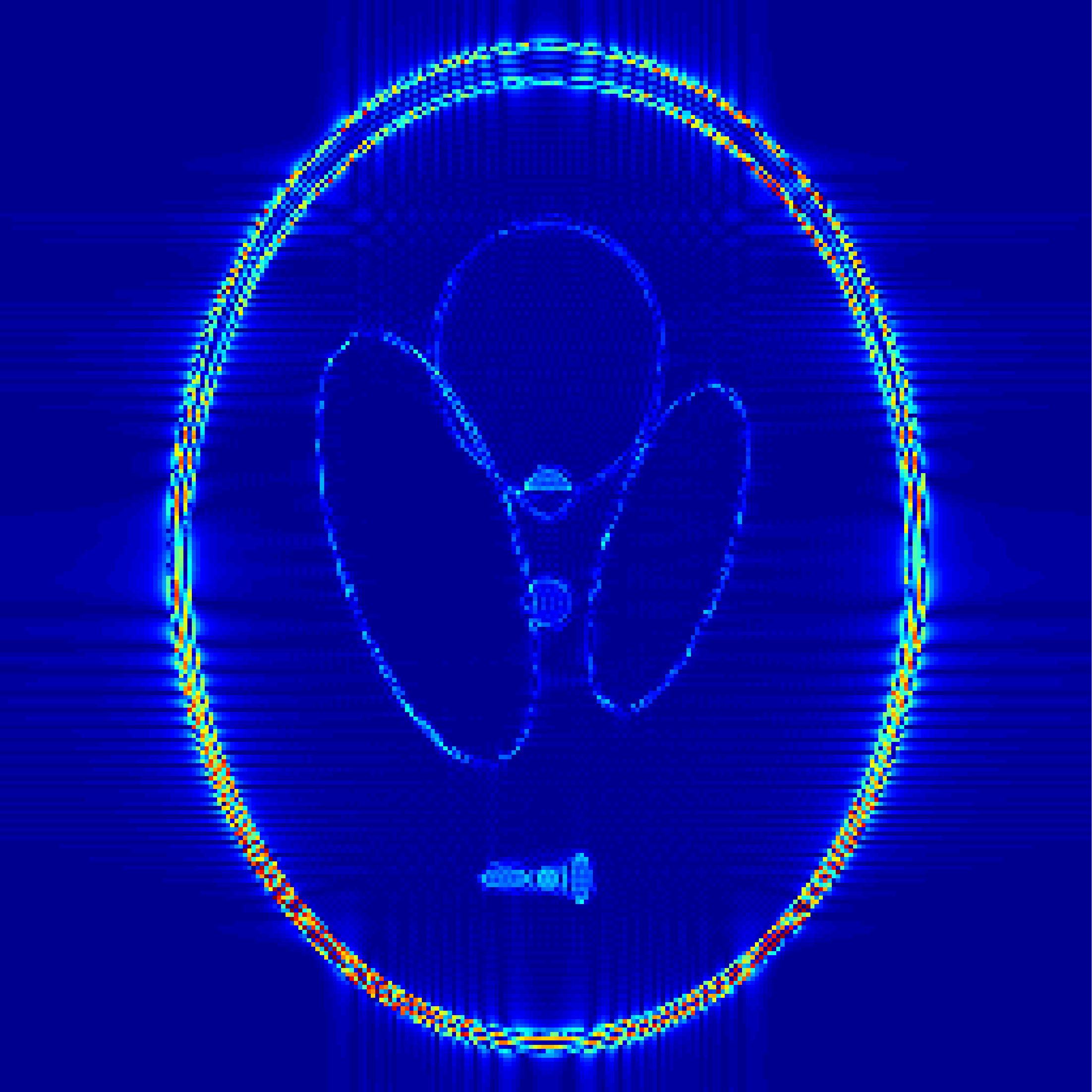}}\hspace{0.005cm}
\subfloat[Haar \cref{HaarModel}]{\label{PhantomHaarError}\includegraphics[width=3.00cm]{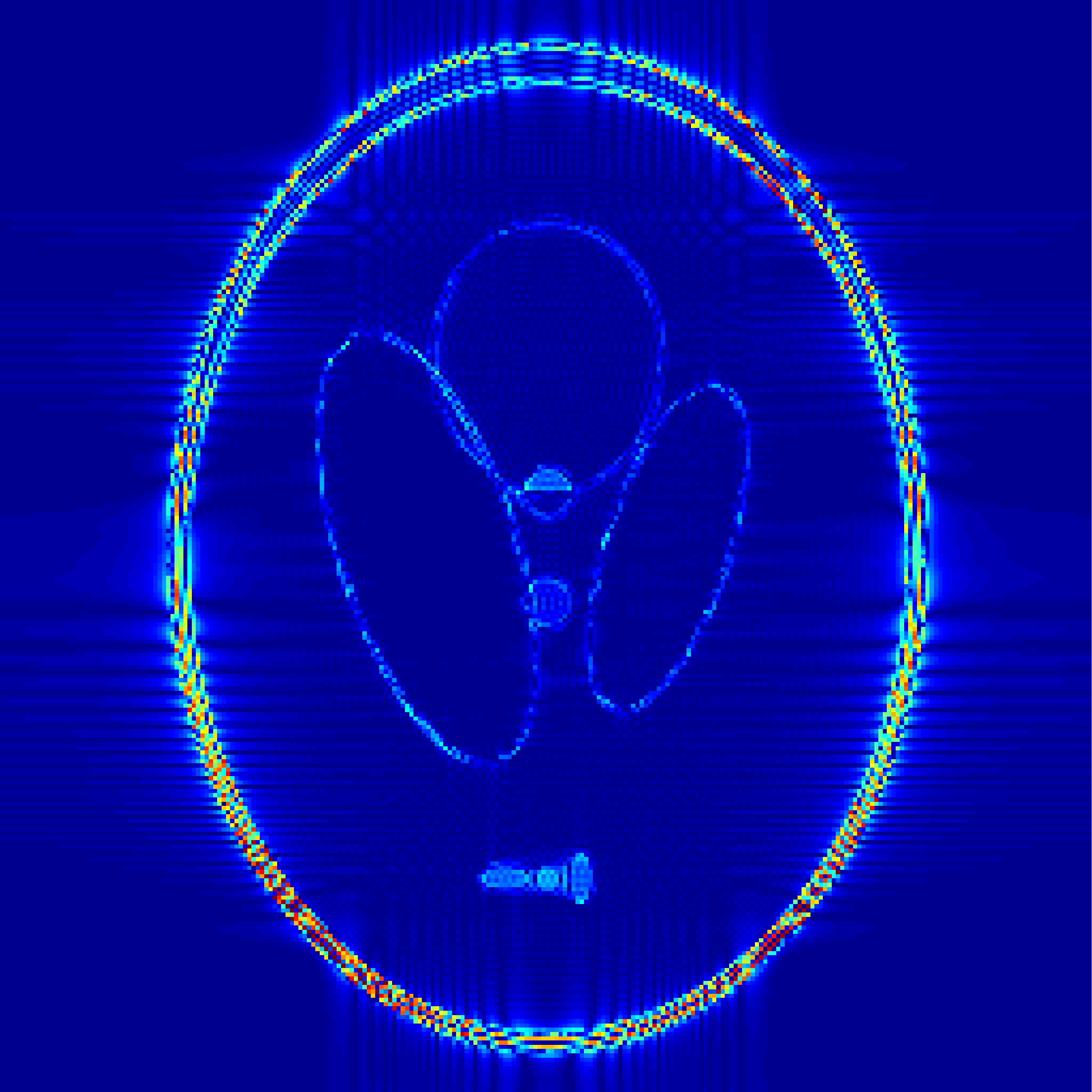}}\hspace{0.005cm}
\subfloat[TLMRI \cref{TLMRI}]{\label{PhantomFRISTError}\includegraphics[width=3.00cm]{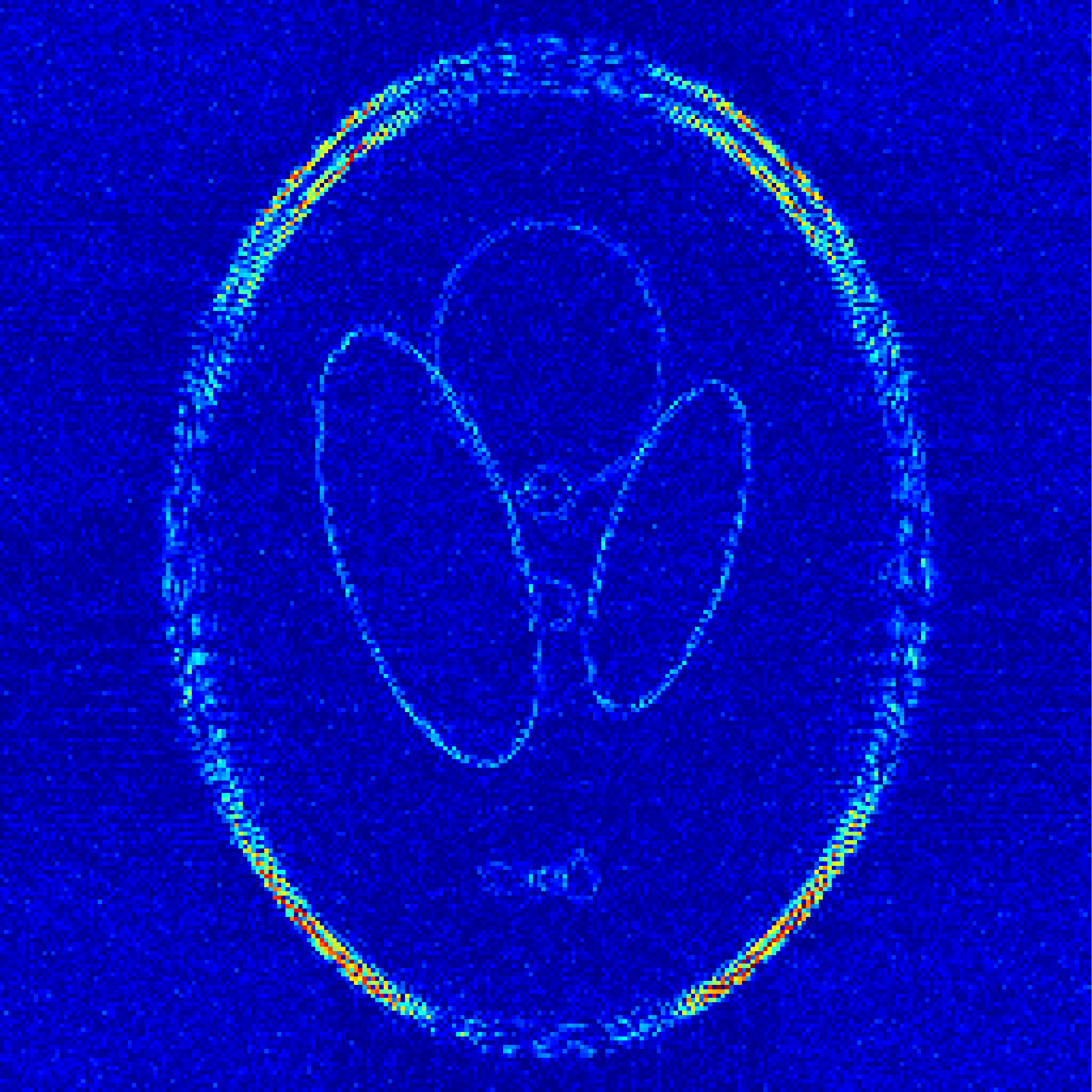}}\vspace{-0.20cm}\\
\subfloat[GIRAF$0$]{\label{PhantomGIRAF0Error}\includegraphics[width=3.00cm]{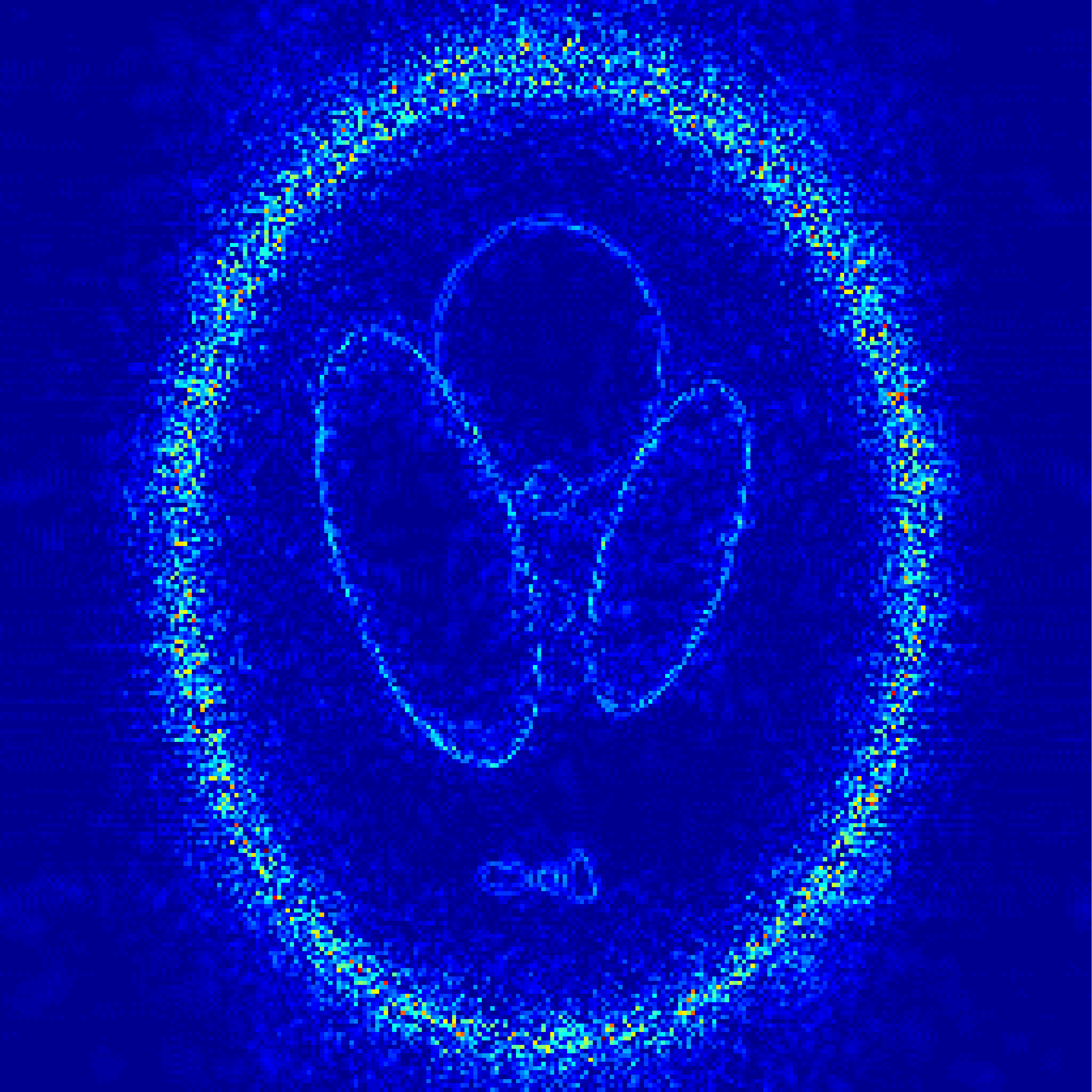}}\hspace{0.005cm}
\subfloat[GIRAF$0.5$]{\label{PhantomGIRAFHalfError}\includegraphics[width=3.00cm]{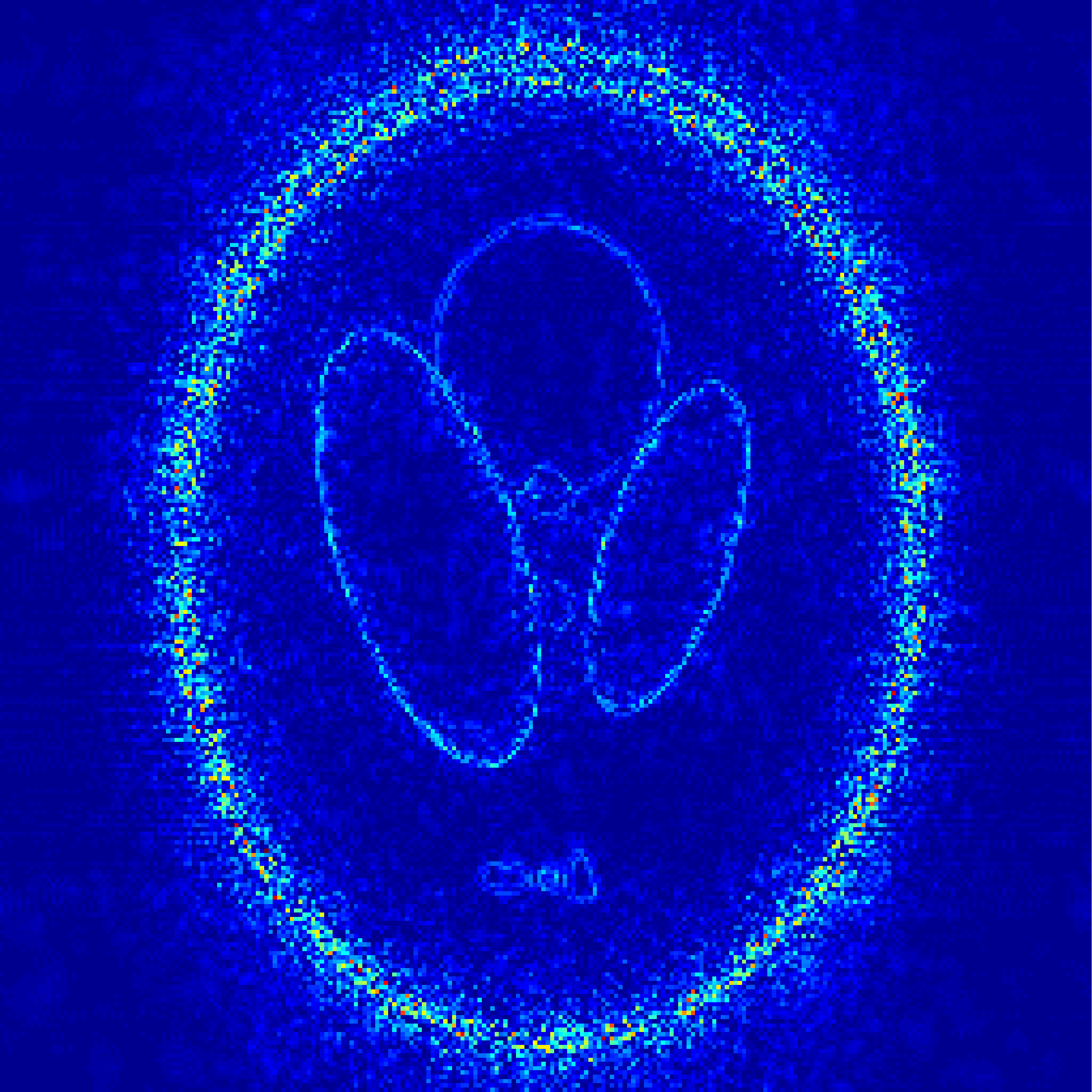}}\hspace{0.005cm}
\subfloat[GIRAF$1$]{\label{PhantomGIRAF1Error}\includegraphics[width=3.00cm]{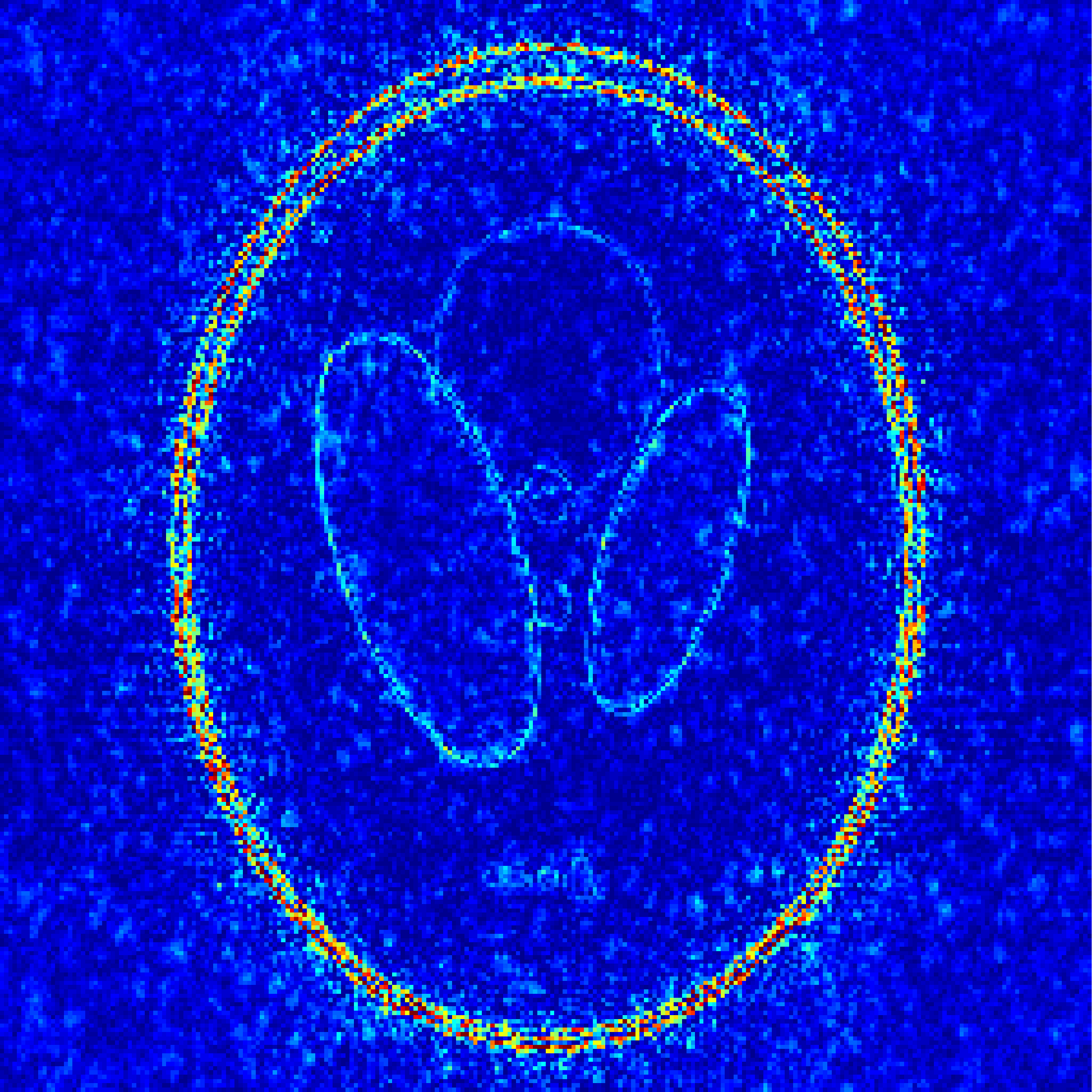}}\hspace{0.005cm}
\subfloat[DDTF \cref{ProposedCSMRIModel}]{\label{PhantomDDTFError}\includegraphics[width=3.00cm]{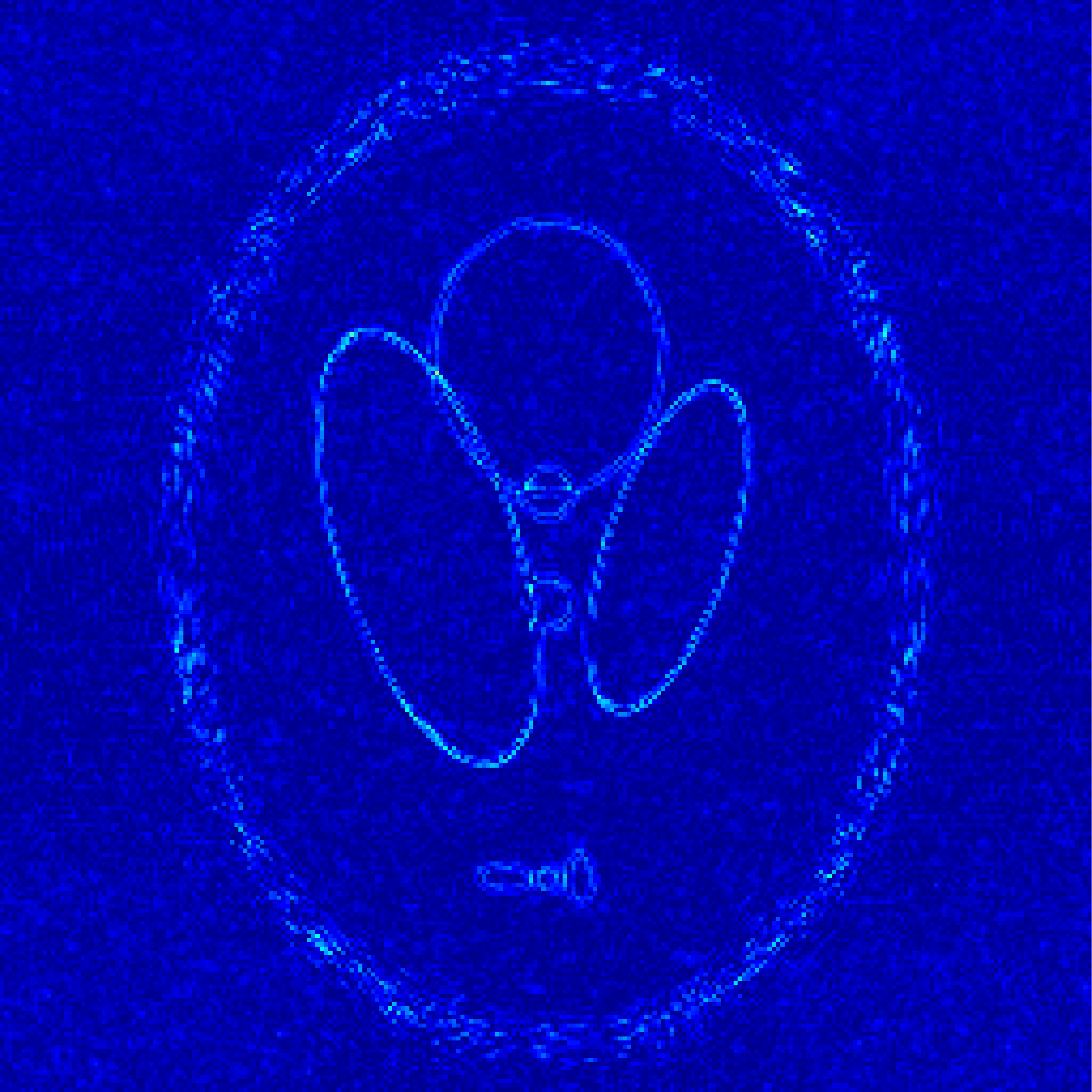}}\vspace{-0.20cm}
\caption{Comparisons of error maps for the phantom experiments. All error maps are displayed in the window level $[0,0.2]$ for fair comparisons.}\label{PhantomResultsErrorMap}
\end{figure}

\begin{figure}[tp!]
\centering
\hspace{-0.1cm}\subfloat[Fully sampled]{\label{PhantomOriginalZoom}\includegraphics[width=3.00cm]{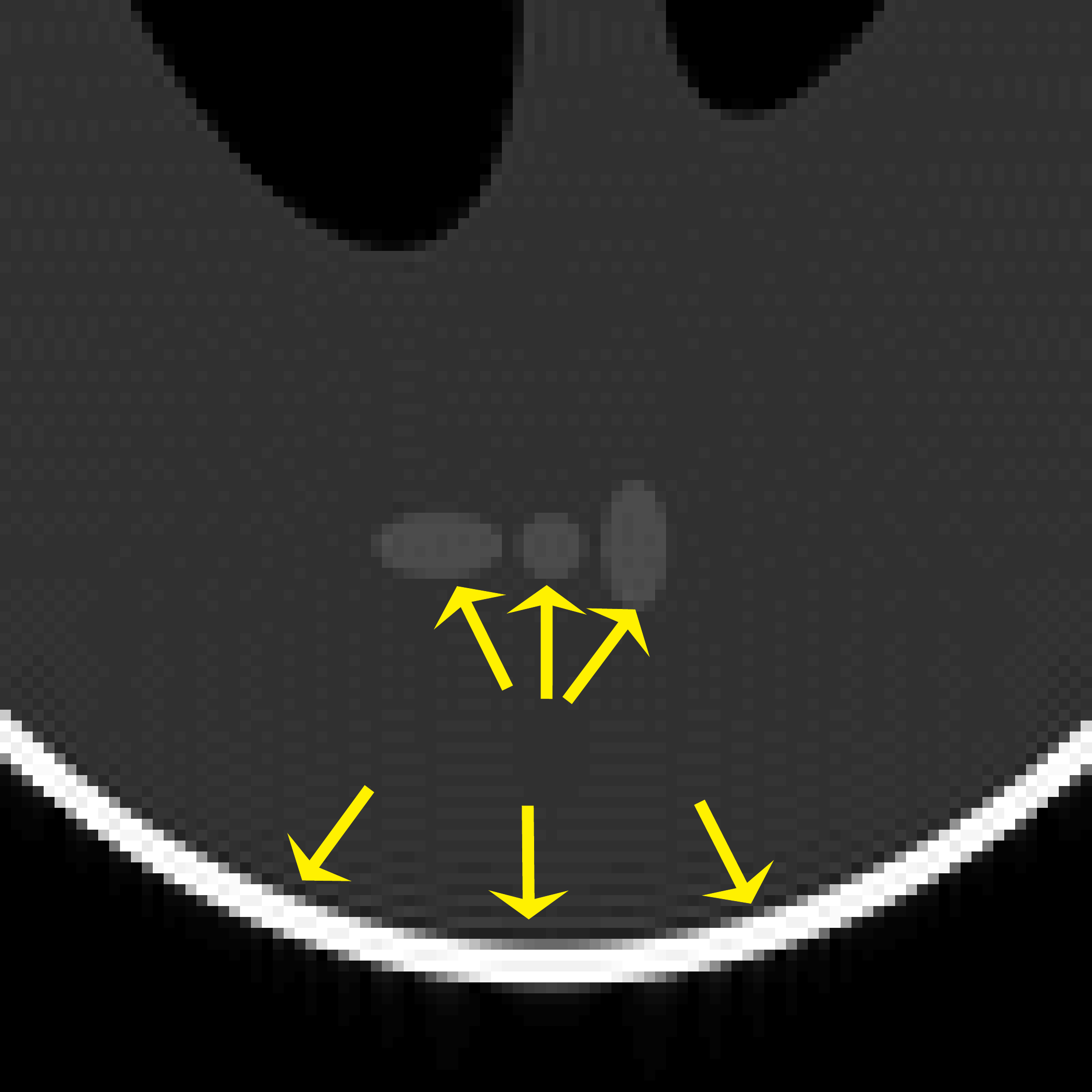}}\hspace{0.005cm}
\subfloat[Zero fill]{\label{PhantomZeroPadZoom}\includegraphics[width=3.00cm]{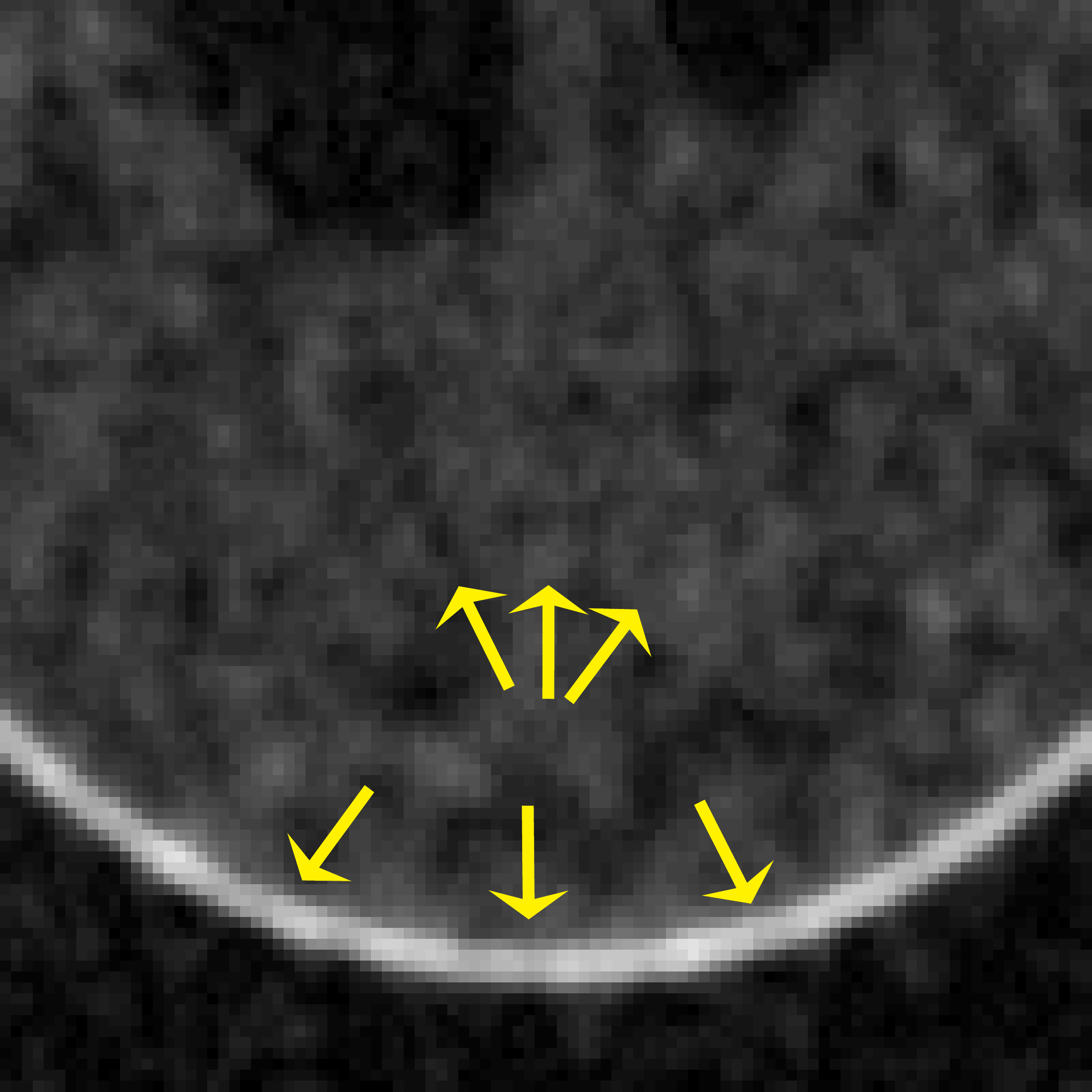}}\hspace{0.005cm}
\subfloat[TV \cref{TVModel}]{\label{PhantomTVZoom}\includegraphics[width=3.00cm]{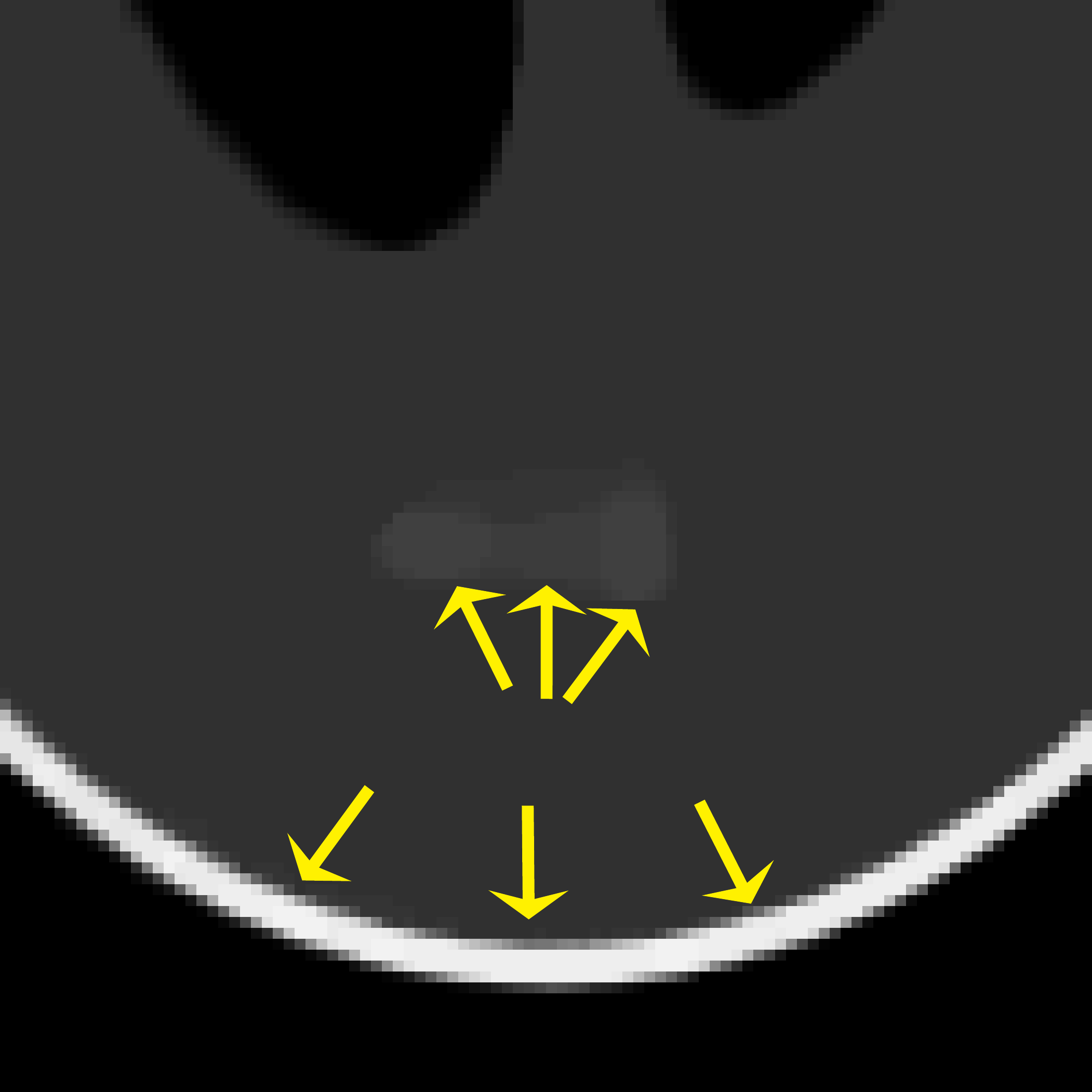}}\hspace{0.005cm}
\subfloat[Haar \cref{HaarModel}]{\label{PhantomHaarZoom}\includegraphics[width=3.00cm]{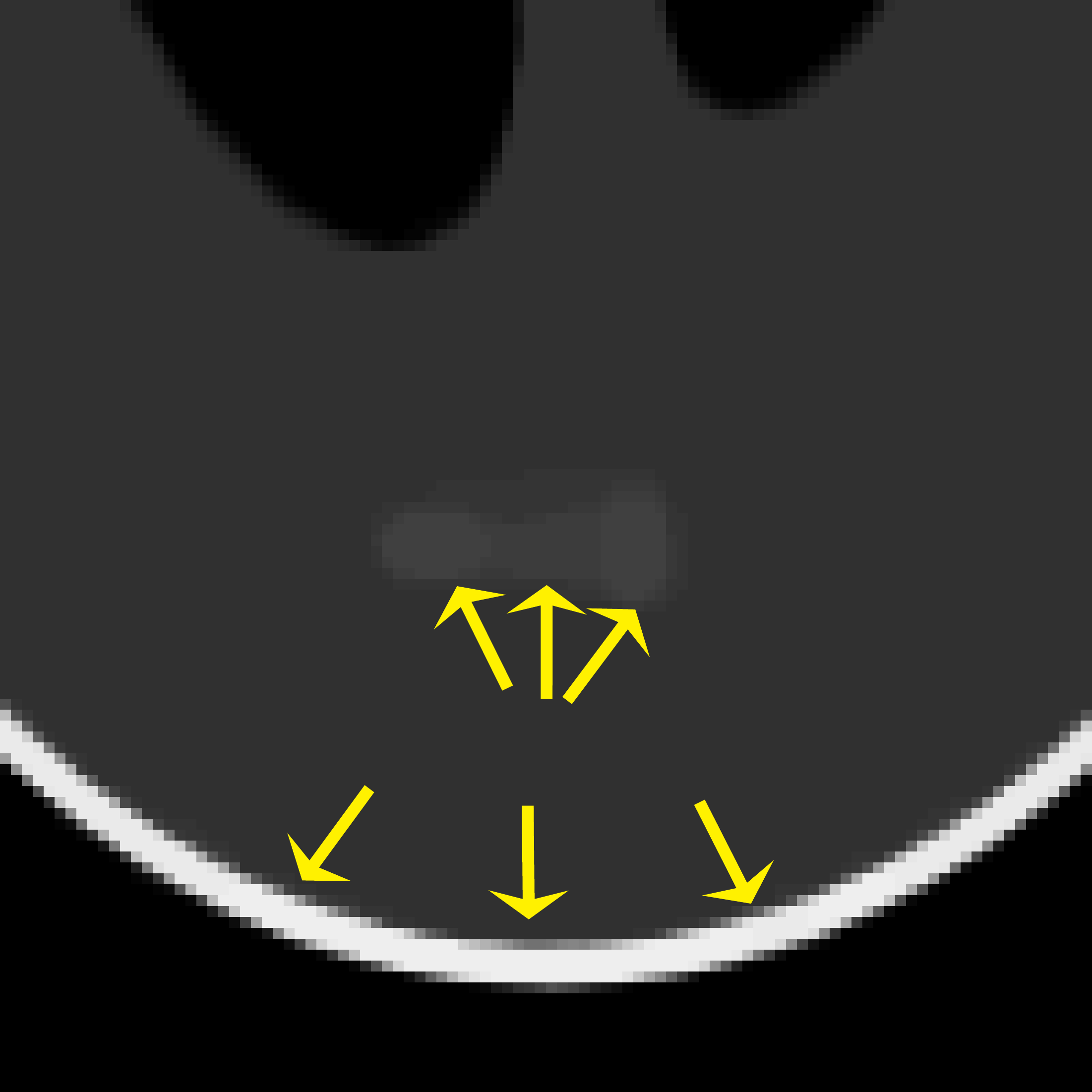}}\hspace{0.005cm}
\subfloat[TLMRI \cref{TLMRI}]{\label{PhantomFRISTZoom}\includegraphics[width=3.00cm]{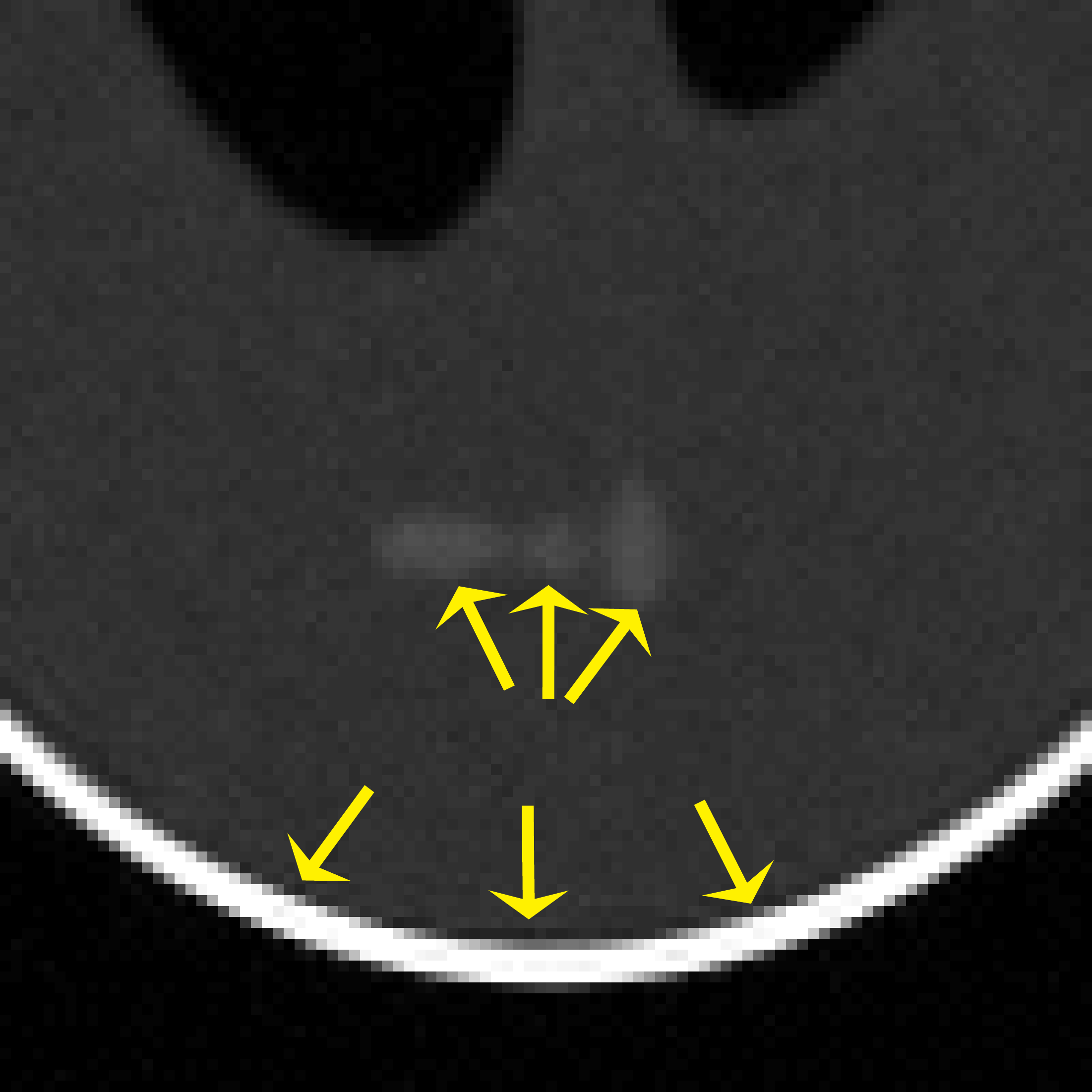}}\vspace{-0.20cm}\\
\subfloat[GIRAF$0$]{\label{PhantomGIRAF0Zoom}\includegraphics[width=3.00cm]{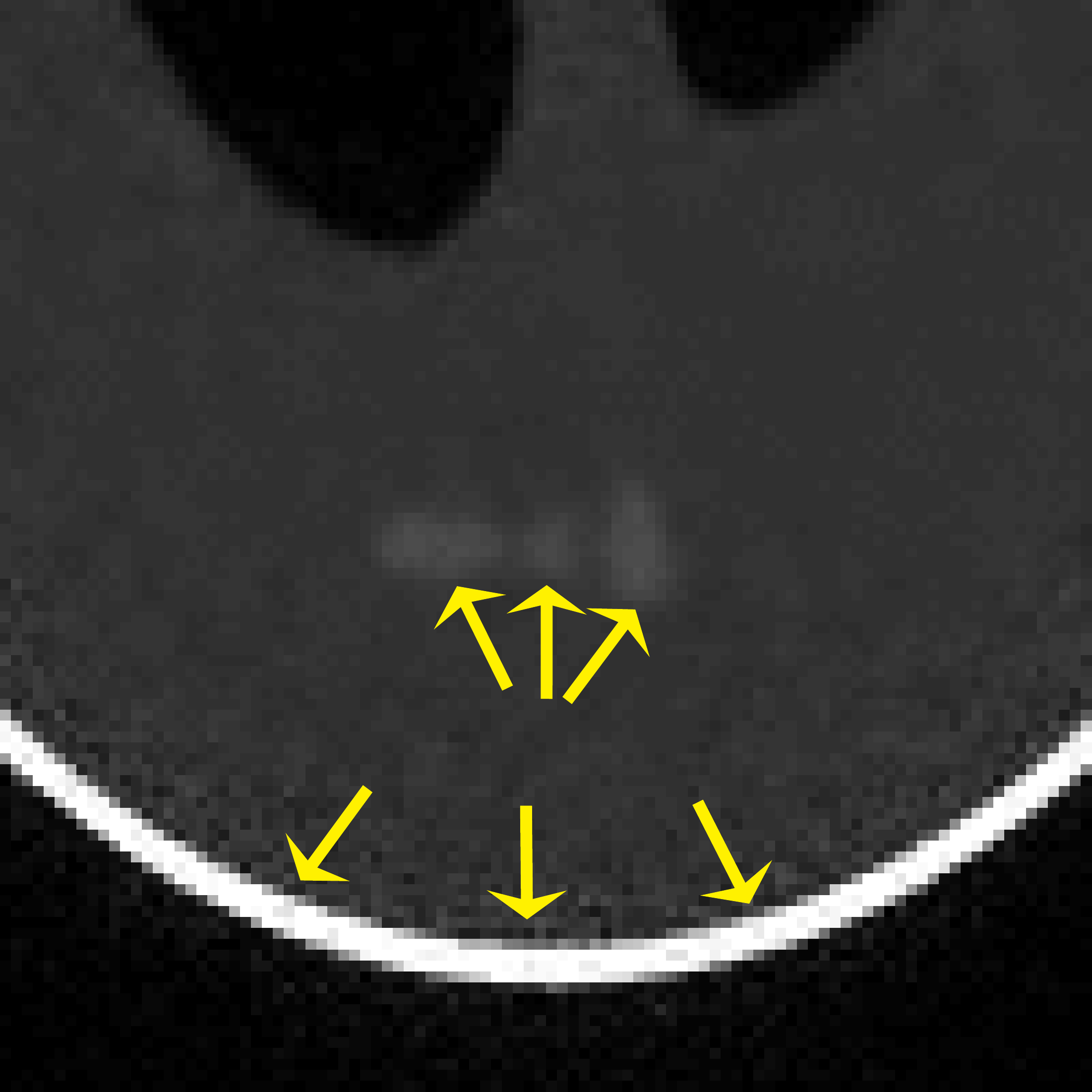}}\hspace{0.005cm}
\subfloat[GIRAF$0.5$]{\label{PhantomGIRAFHalfZoom}\includegraphics[width=3.00cm]{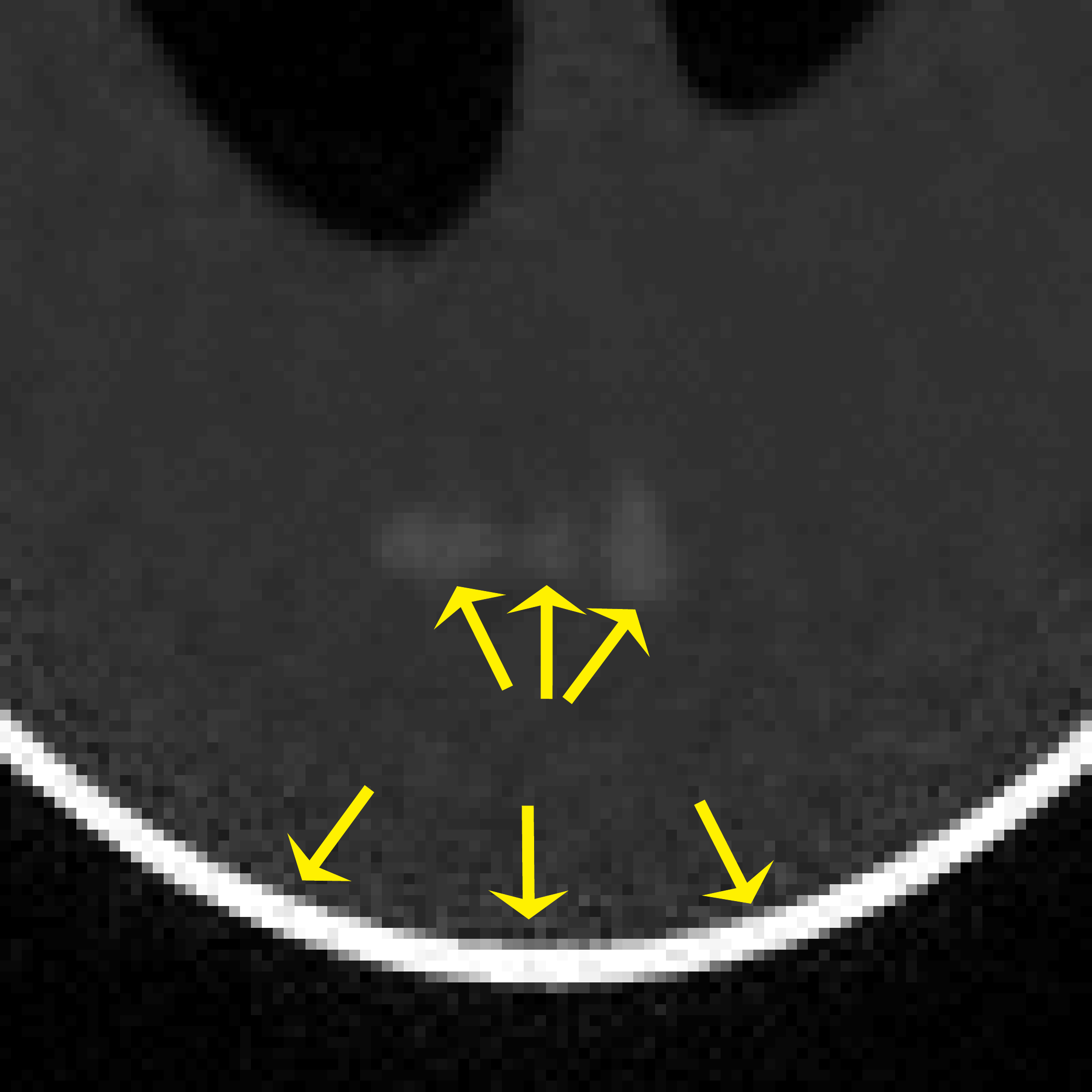}}\hspace{0.005cm}
\subfloat[GIRAF$1$]{\label{PhantomGIRAF1Zoom}\includegraphics[width=3.00cm]{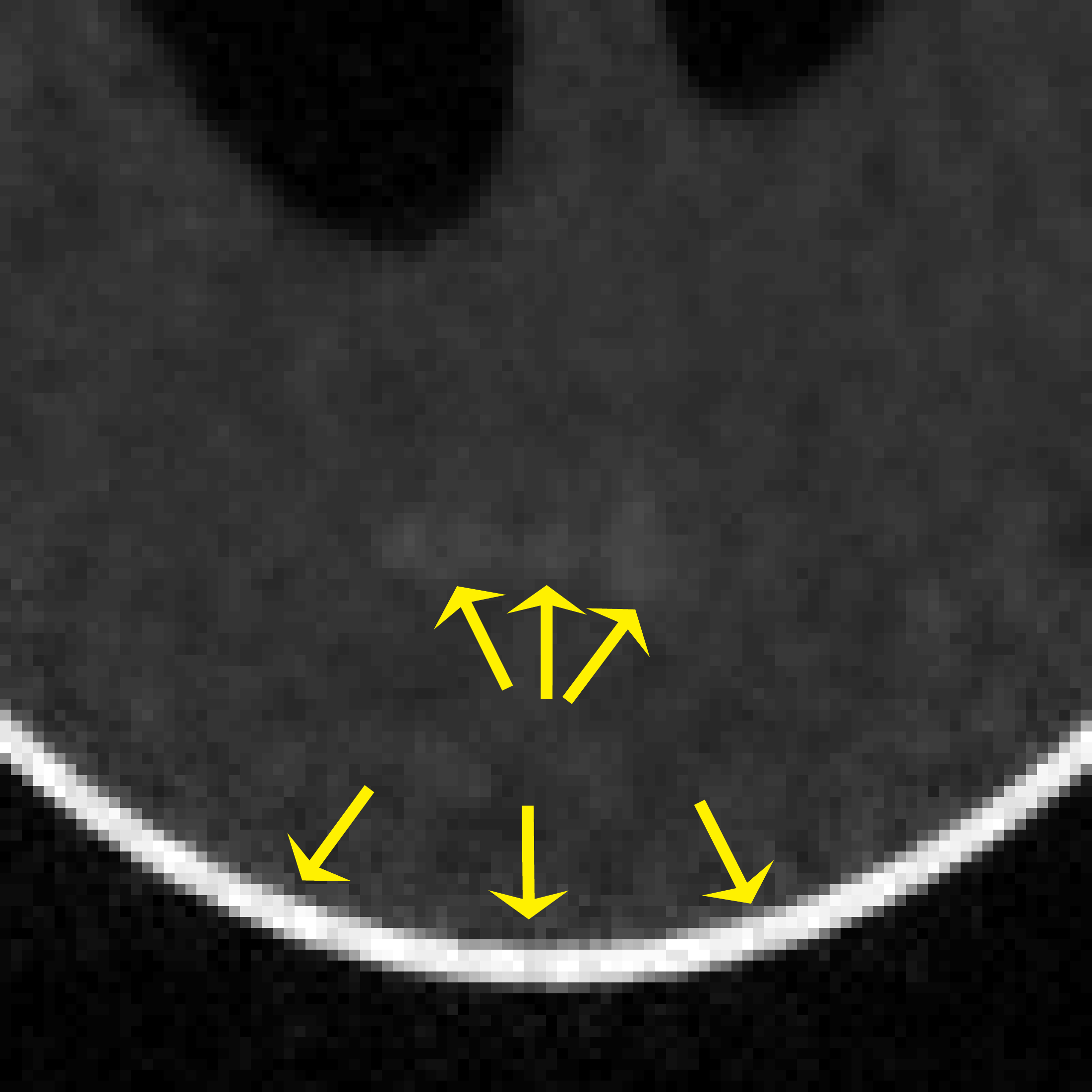}}\hspace{0.005cm}
\subfloat[DDTF \cref{ProposedCSMRIModel}]{\label{PhantomDDTFZoom}\includegraphics[width=3.00cm]{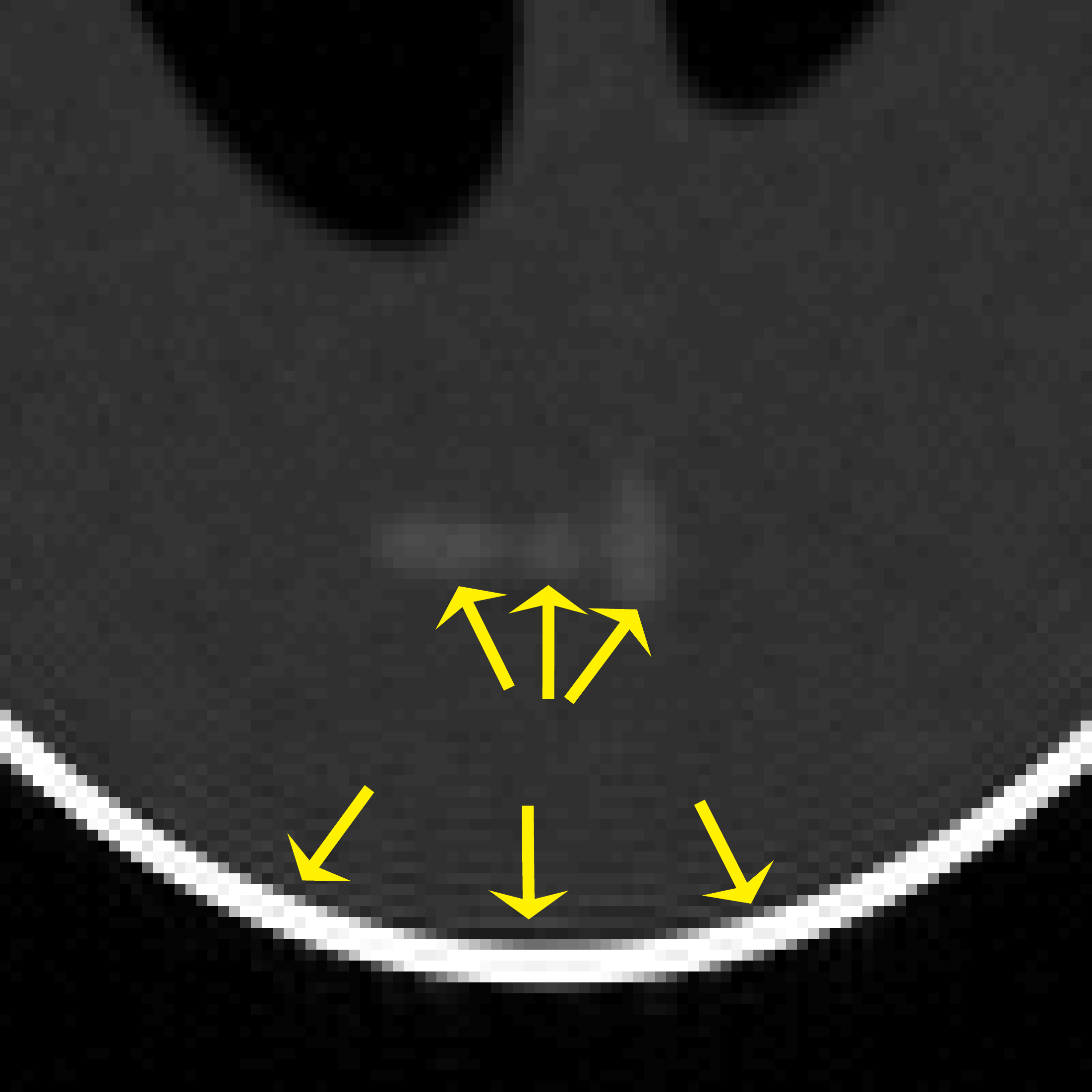}}\vspace{-0.20cm}
\caption{Zoom-in views of \cref{PhantomResults}. The yellow arrows indicate the regions worth noticing.}\label{PhantomResultsZoom}
\end{figure}

\cref{PhantomTable} summarizes the SNR and the HFEN of the aforementioned restoration models, and \cref{PhantomResults} displays the visual comparisons with the zoom-in views in \cref{PhantomResultsZoom} and the error maps in \cref{PhantomResultsErrorMap}, respectively. We can see that the proposed data driven tight frame model \cref{ProposedCSMRIModel} consistently outperforms both the on-the-grid approaches (\cref{TVModel,HaarModel,TLMRI}) and the existing off-the-grid approaches \cref{Schatten} with a smaller error map. Noting that \cref{ProposedCSMRIModel} is an off-the-grid approach, the experimental results also suggest that the off-the-grid approaches have better performance in the CS-MRI due to its ability to reducing the basis mismatch between the true support (or the true singularity) in continuum and the discrete grid. In fact, due to this basis mismatch, we can see from \cref{PhantomTV,PhantomHaar,PhantomTVZoom,PhantomHaarZoom} that the on-the-grid approaches lead to the distortions of three small ellipses, and the errors concentrate on the edges (\cref{PhantomTVError,PhantomHaarError}) compared to the off-the-grid approaches. Even though \cref{TLMRI} is able to restore the three small ellipses, the restored image is still overly smoothed compared to the proposed model with errors concentrated near the edges, as we can see from \cref{PhantomFRISTError,PhantomFRISTZoom}.

It is also worth noting that among the off-the-grid approaches, the proposed DDTF model introduces less artifacts near the edges. In the literature, the noise in the k-space data affects $\rank\left(\bmH\left(\La\bsv\right)\right)$ even in the fully sampled case as the weight matrix $\La$ amplifies the noise in the high frequencies. Hence under such an amplified noise, it is likely that the direct rank minimization leads to the artifacts near the edges corresponding to the high frequencies in the frequency domain, as shown in \cref{PhantomGIRAF0,PhantomGIRAFHalf,PhantomGIRAF1,PhantomGIRAF0Zoom,PhantomGIRAFHalfZoom,PhantomGIRAF1Zoom,PhantomGIRAF0Error,PhantomGIRAFHalfError,PhantomGIRAF1Error}. In contrast, the sparse approximation of $\La\bsv$ can achieve the denoising effect in spite of the amplified noise, leading to the better restoration results with less artifacts near the edges.

\begin{figure}[tp!]
\centering
\hspace{-0.1cm}\subfloat[Fully sampled]{\label{PhantomOriginalk2}\includegraphics[width=3.00cm]{PhantomOriginalk.pdf}}\hspace{0.005cm}
\subfloat[Zero fill]{\label{PhantomUndersamplek}\includegraphics[width=3.00cm]{PhantomUndersamplek.pdf}}\hspace{0.005cm}
\subfloat[TV \cref{TVModel}]{\label{PhantomTVk}\includegraphics[width=3.00cm]{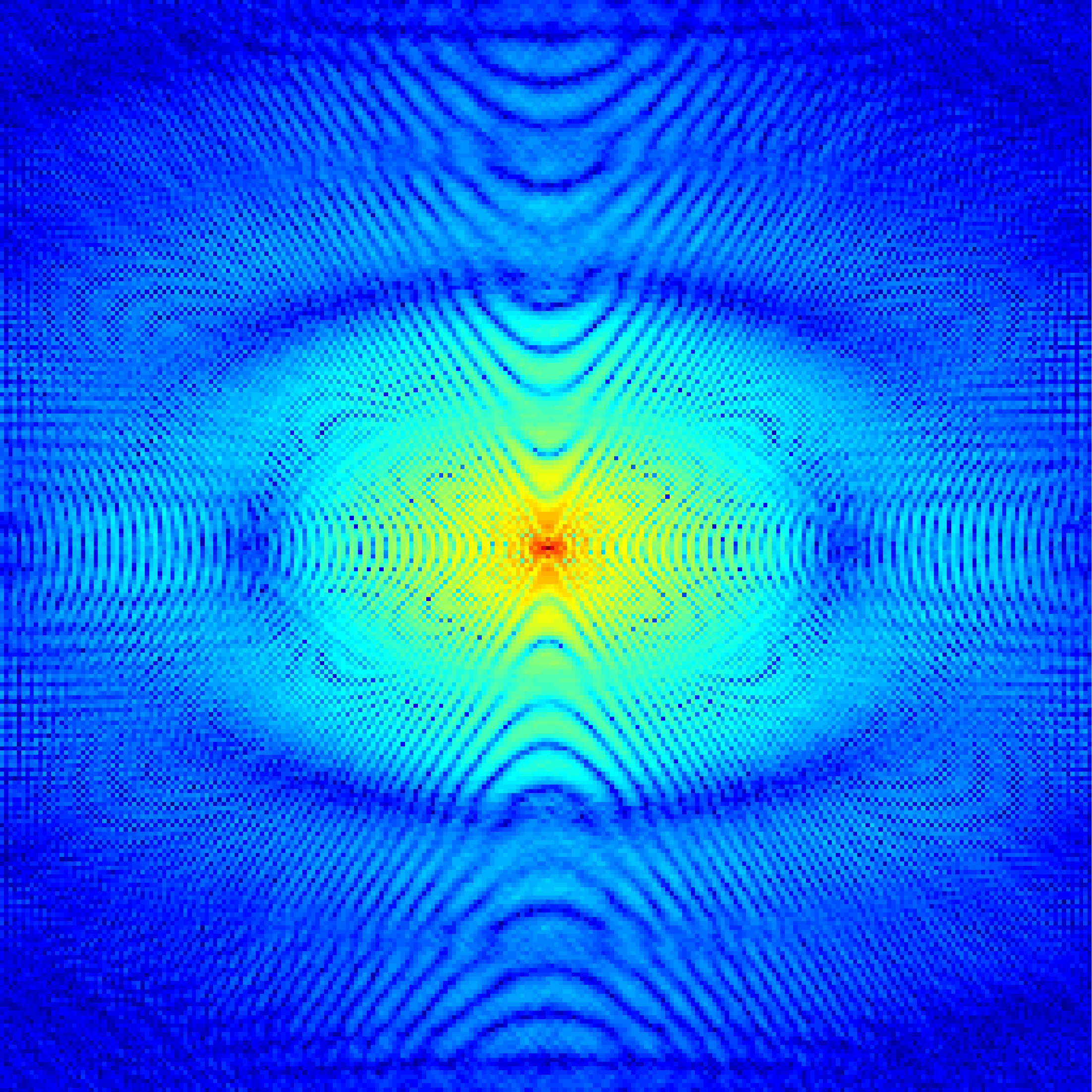}}\hspace{0.005cm}
\subfloat[Haar \cref{HaarModel}]{\label{PhantomHaark}\includegraphics[width=3.00cm]{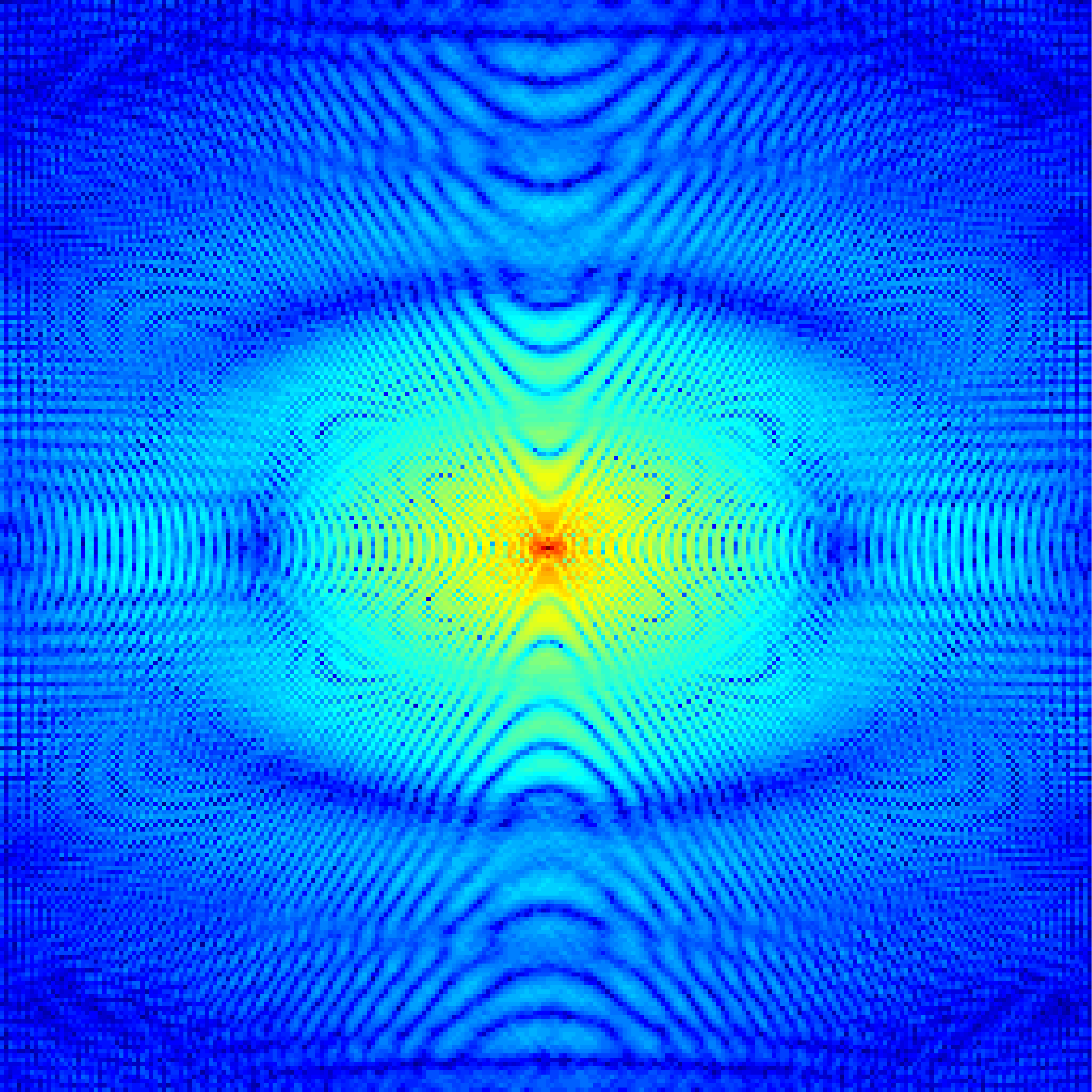}}\hspace{0.005cm}
\subfloat[TLMRI \cref{TLMRI}]{\label{PhantomFRISTk}\includegraphics[width=3.00cm]{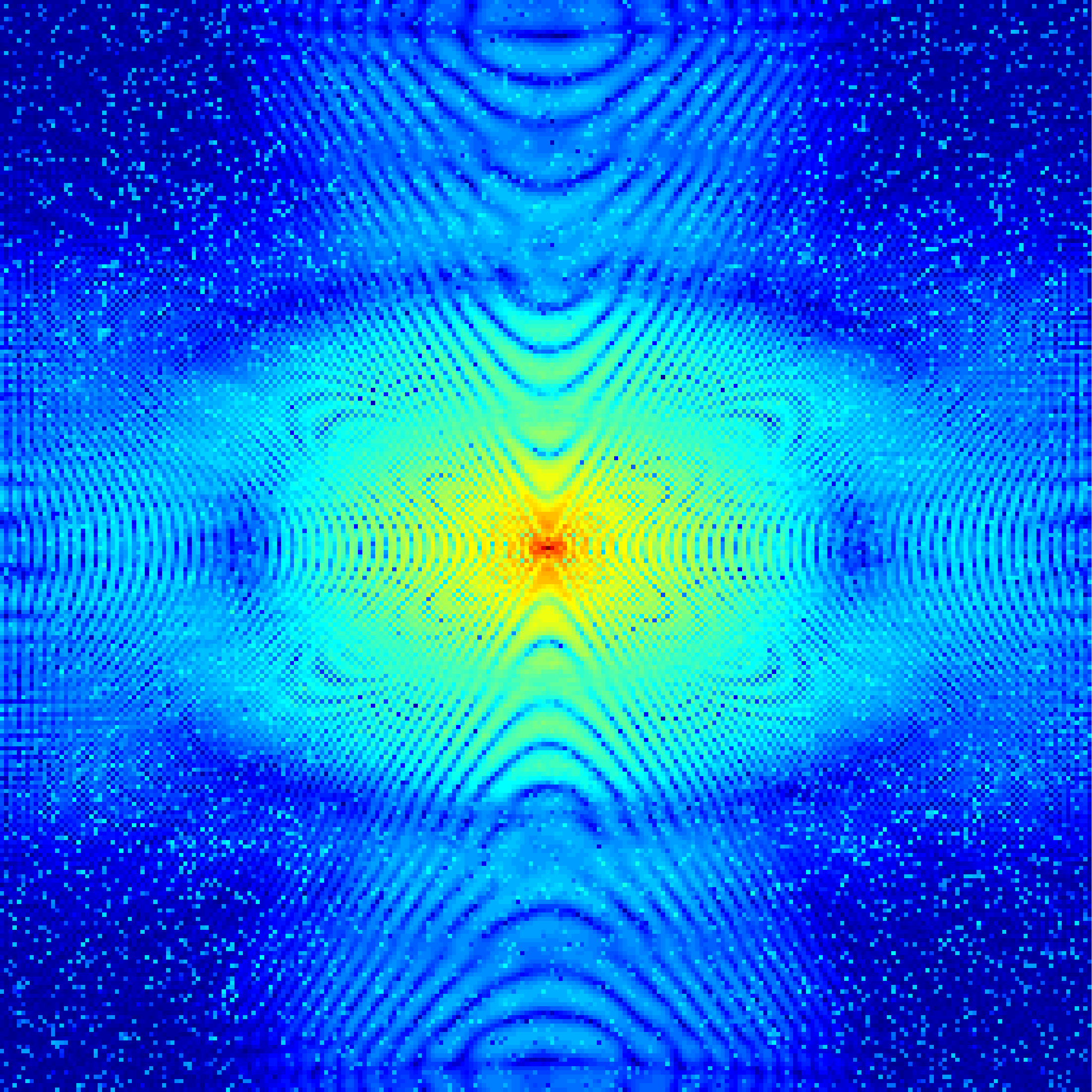}}\vspace{-0.20cm}\\
\subfloat[GIRAF$0$]{\label{PhantomGIRAF0k}\includegraphics[width=3.00cm]{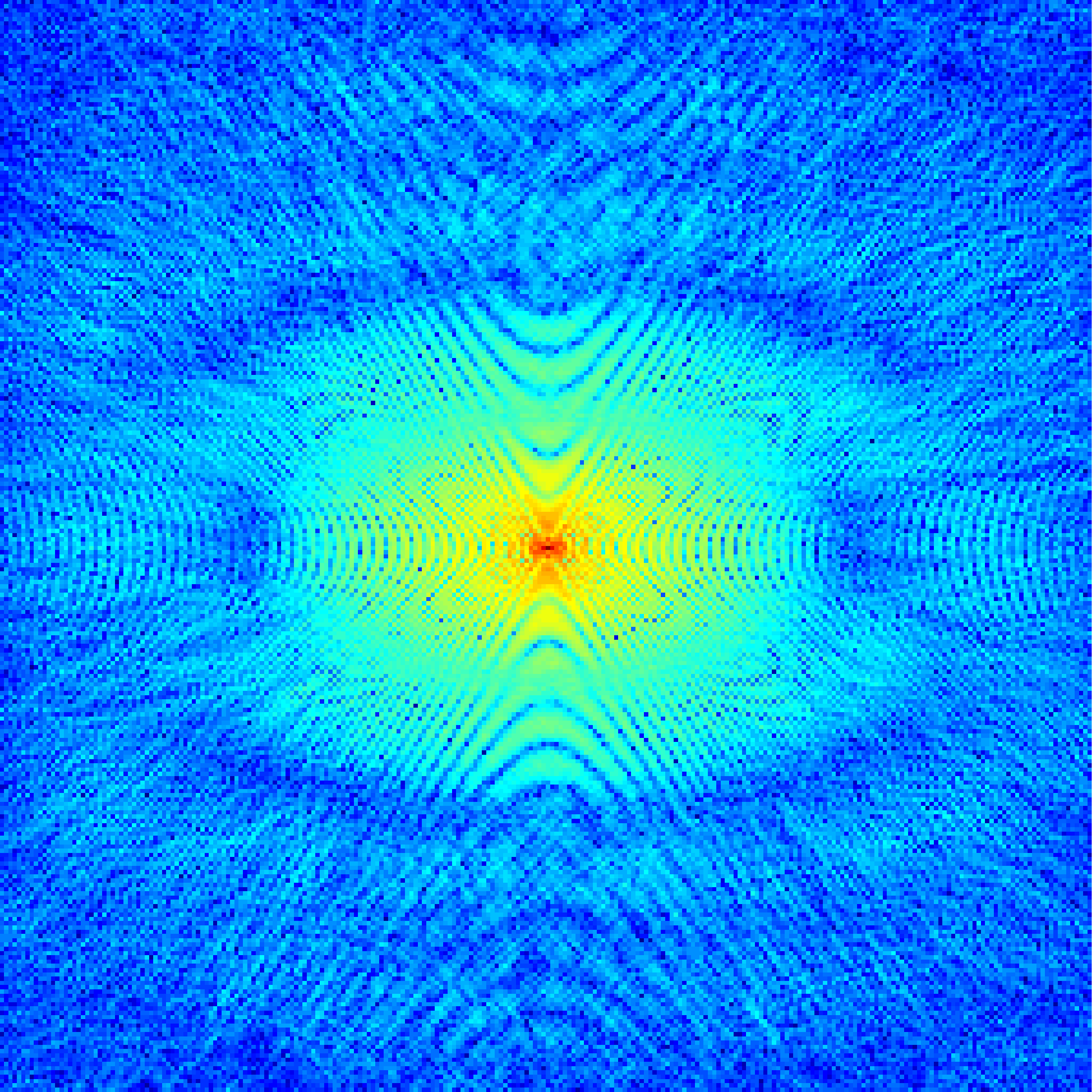}}\hspace{0.005cm}
\subfloat[GIRAF$0.5$]{\label{PhantomGIRAFHalfk}\includegraphics[width=3.00cm]{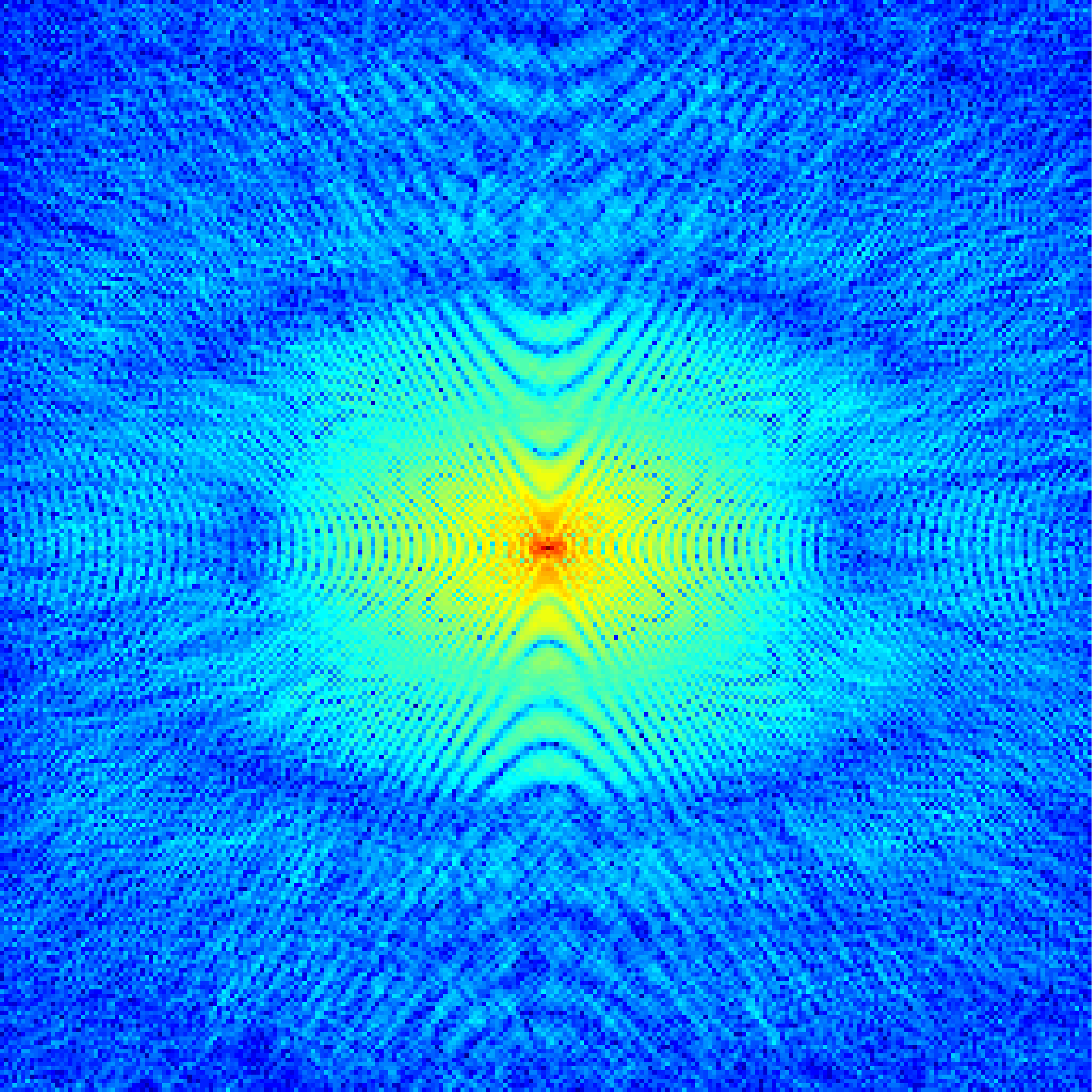}}\hspace{0.005cm}
\subfloat[GIRAF$1$]{\label{PhantomGIRAF1k}\includegraphics[width=3.00cm]{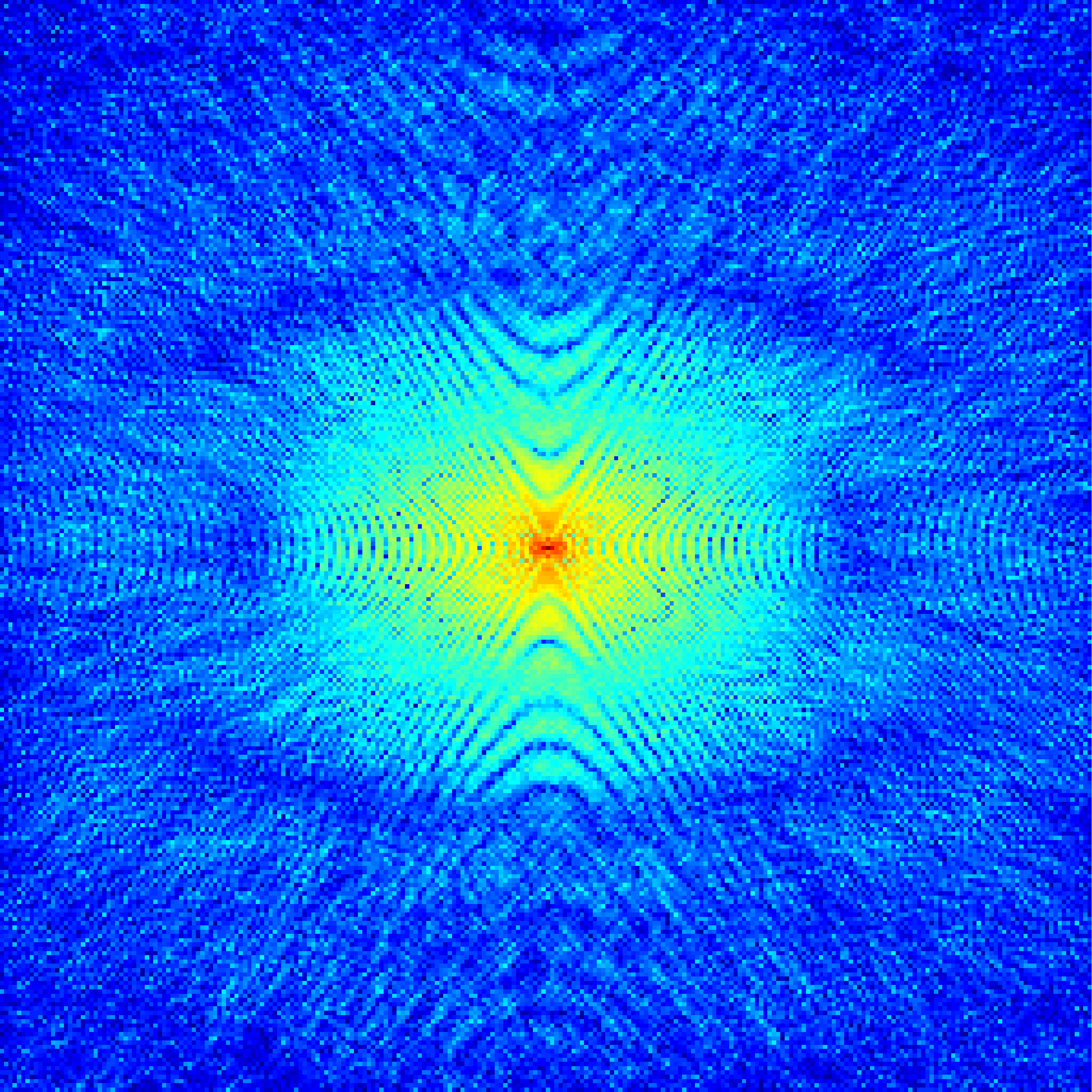}}\hspace{0.005cm}
\subfloat[DDTF \cref{ProposedCSMRIModel}]{\label{PhantomDDTFk}\includegraphics[width=3.00cm]{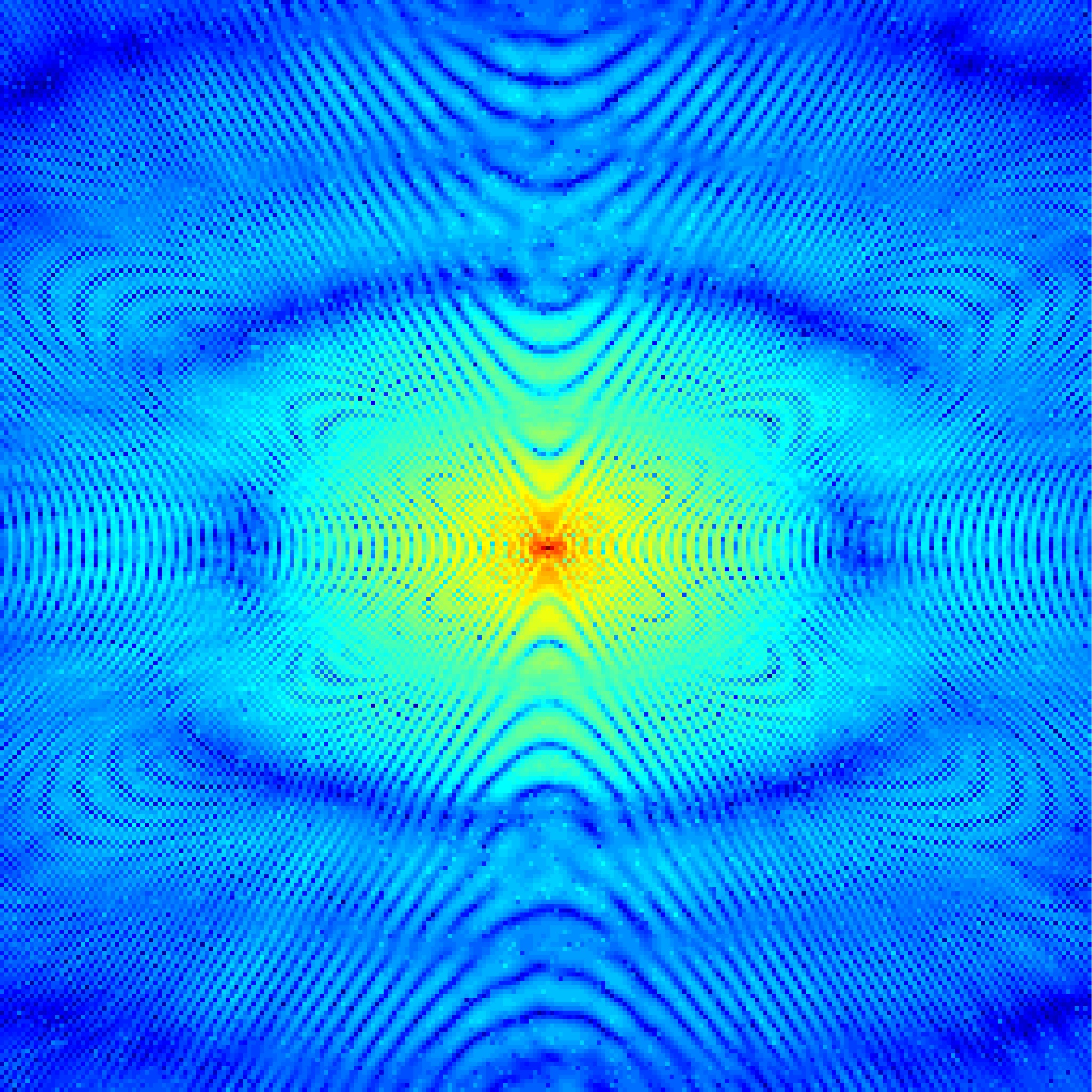}}\vspace{-0.20cm}
\caption{Comparisons of k-space data for the phantom experiments in the log scale. All restored k-space data are displayed in the window level $[0,9]$ for fair comparisons.}\label{PhantomResultsk}
\end{figure}

For further comparisons, we also present the restored k-space data (in the log scale) in \cref{PhantomResultsk}. Note that since the sampling is dense in the low frequencies while the high frequencies are loosely sampled, the restoration qualities depend heavily on the restoration accuracy of high frequency k-space data. Indeed, we can see from \cref{PhantomTVk,PhantomHaark,PhantomFRISTk} that the restored k-space data by \cref{TVModel,HaarModel,TLMRI} decays faster than the original one, which also leads to the inferior restoration results. Even though the GIRAF models are in general able to restore the high frequency part better than the on-the-grid approaches, they still fail to restore the dominant structures on the high frequencies, as shown in \cref{PhantomGIRAF0k,PhantomGIRAFHalfk,PhantomGIRAF1k}. In contrast, the proposed model \cref{ProposedCSMRIModel} is able to restore the high frequency k-space data in spite of the loose sampling, which also results in the improvements over the existing approaches. In summary, our proposed DDTF CS-MRI model shows the overall better restoration quality in both the indices (SNR and HFEN) and the visual quality.

Finally, we note that even though \cref{ProposedCSMRIModel} requires more iterations than \cref{Schatten} for convergence, the number of iterations is less than \cref{TVModel,HaarModel,TLMRI}. We also note that the BCD algorithm of \cref{TLMRI} fails to meet the stopping criterion. It is likely that due to the subsequence convergence, the BCD algorithm for \cref{TLMRI} behaves inferior to other models, whose algorithms guarantee the global convergence. Moreover, it takes \cref{TLMRI} approximately $9\mathrm{hrs}$ to perform $600$ iterations. It is most likely that the generation of patch and the projection of a coefficient matrix onto the $\ell_0$ ball are most time consuming, affecting the CPU time. Meanwhile, despite the superior restoration performance, it takes \cref{ProposedCSMRIModel} approximately $24\mathrm{mins}$ to meet the stopping criterion due to a number of discrete convolutions. Even though this is shorter than \cref{TLMRI}, the proposed approach is still time-consuming compared to \cref{TVModel,HaarModel,Schatten}. Nevertheless, our approach is able to gain on the quality of the restored image.

\subsection{Real MR image experiments}\label{RealMRExperiments}

\begin{figure}[tp!]
\centering
\hspace{-0.1cm}\subfloat[Fully sampled]{\label{RealOriginalk}\includegraphics[width=3.0cm]{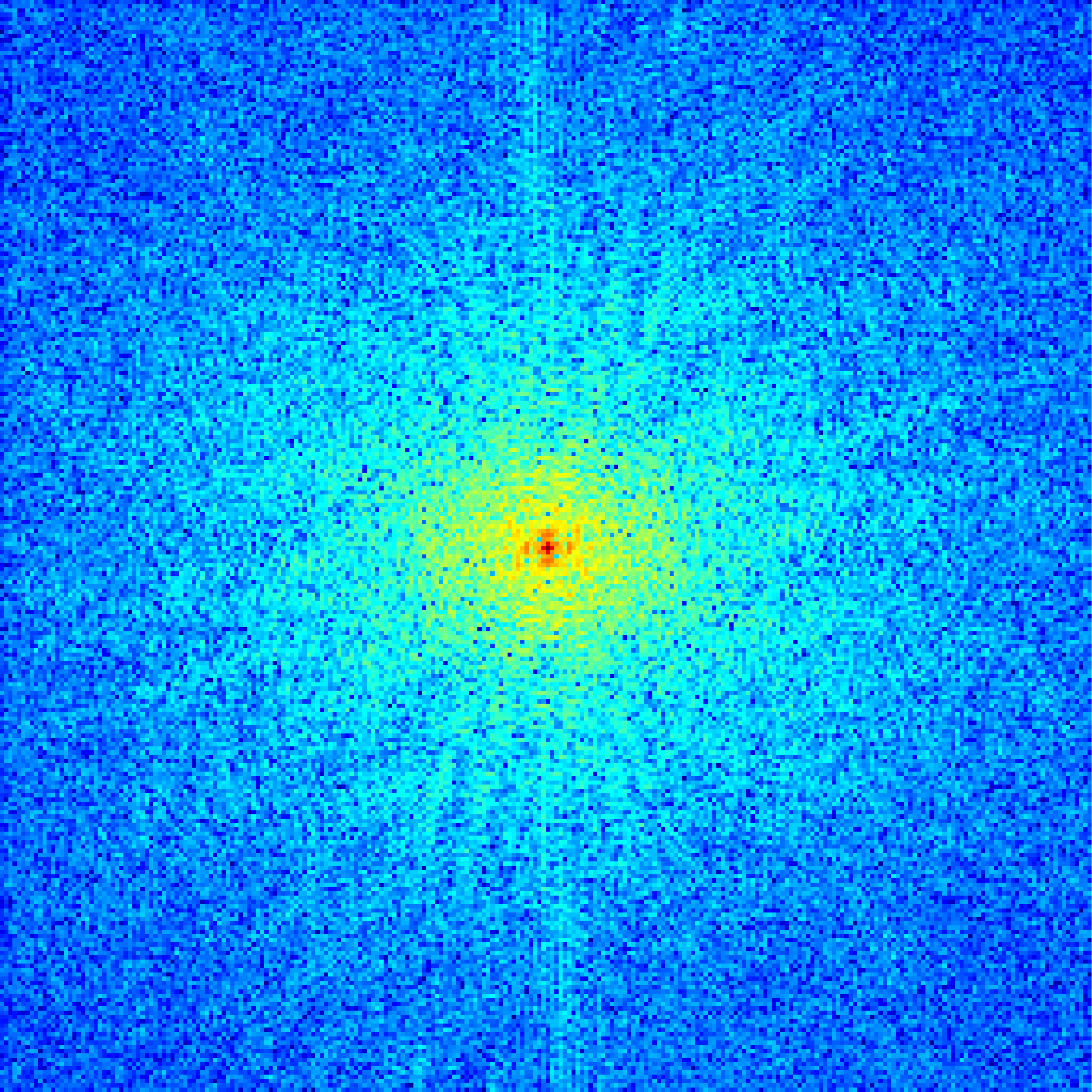}}\hspace{0.005cm}
\subfloat[Ground truth]{\label{RealOriginal}\includegraphics[width=3.0cm]{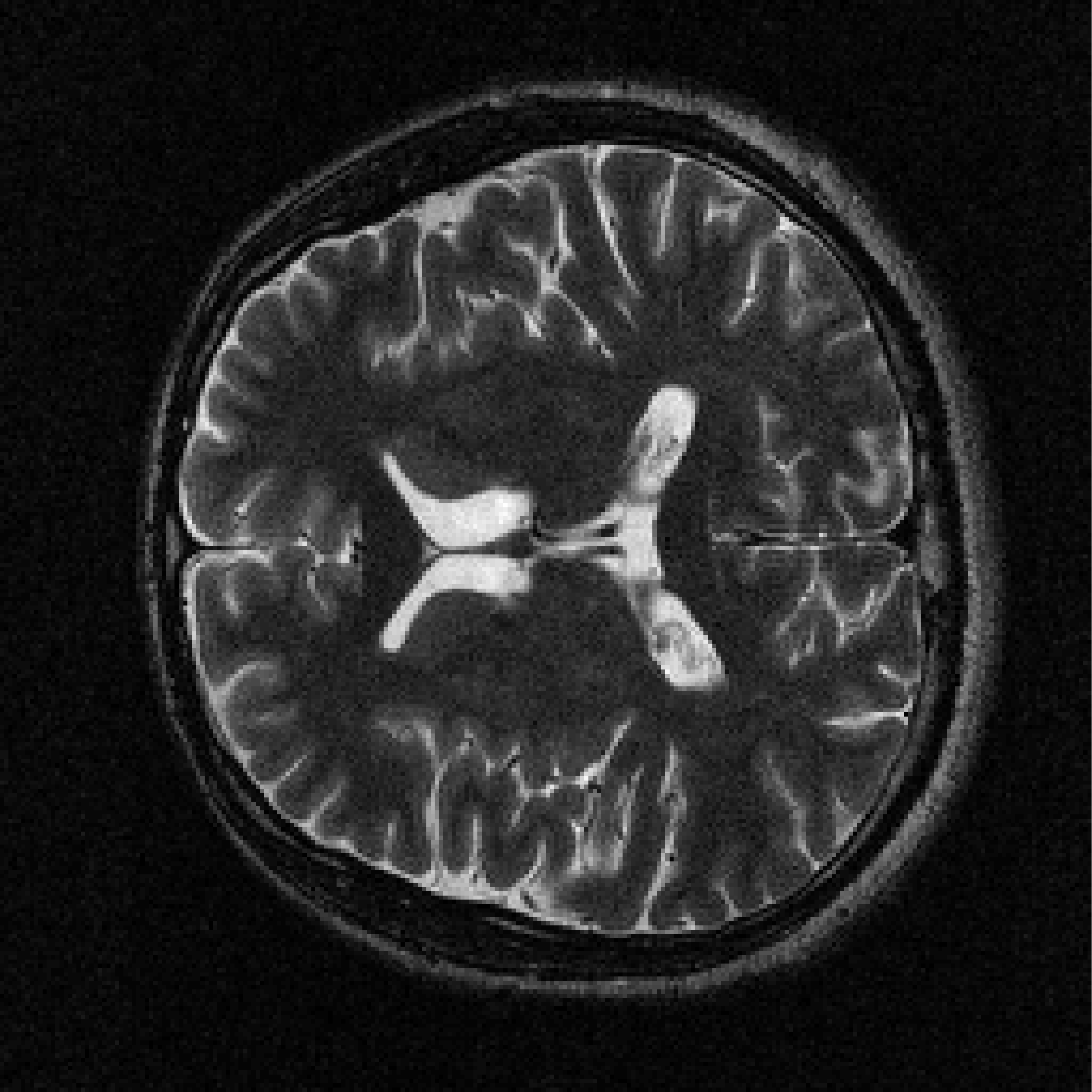}}\hspace{0.005cm}
\subfloat[Sample mask]{\label{RealMask}\includegraphics[width=3.0cm]{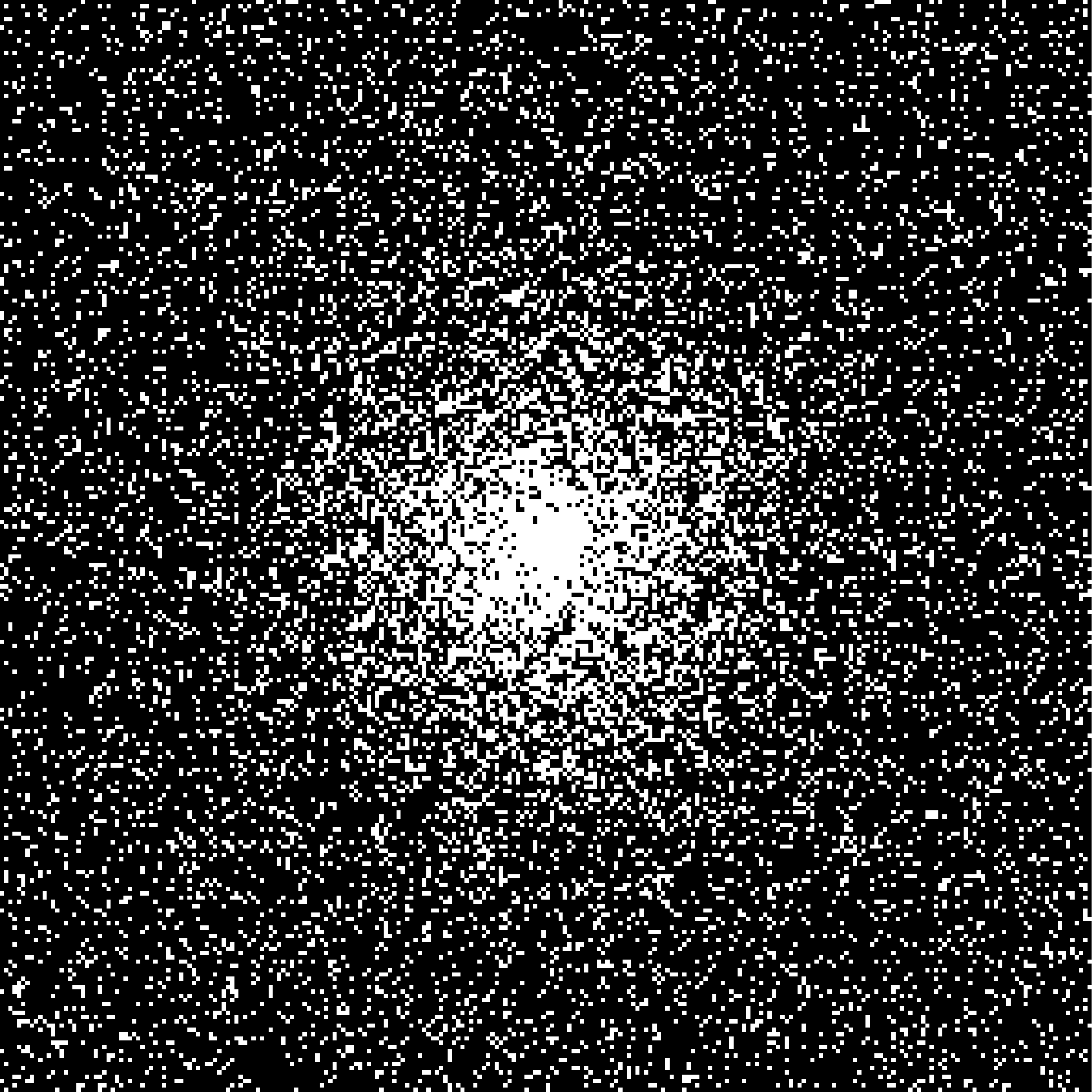}}\hspace{0.005cm}
\subfloat[Undersampled]{\label{RealUndersamplek2}\includegraphics[width=3.0cm]{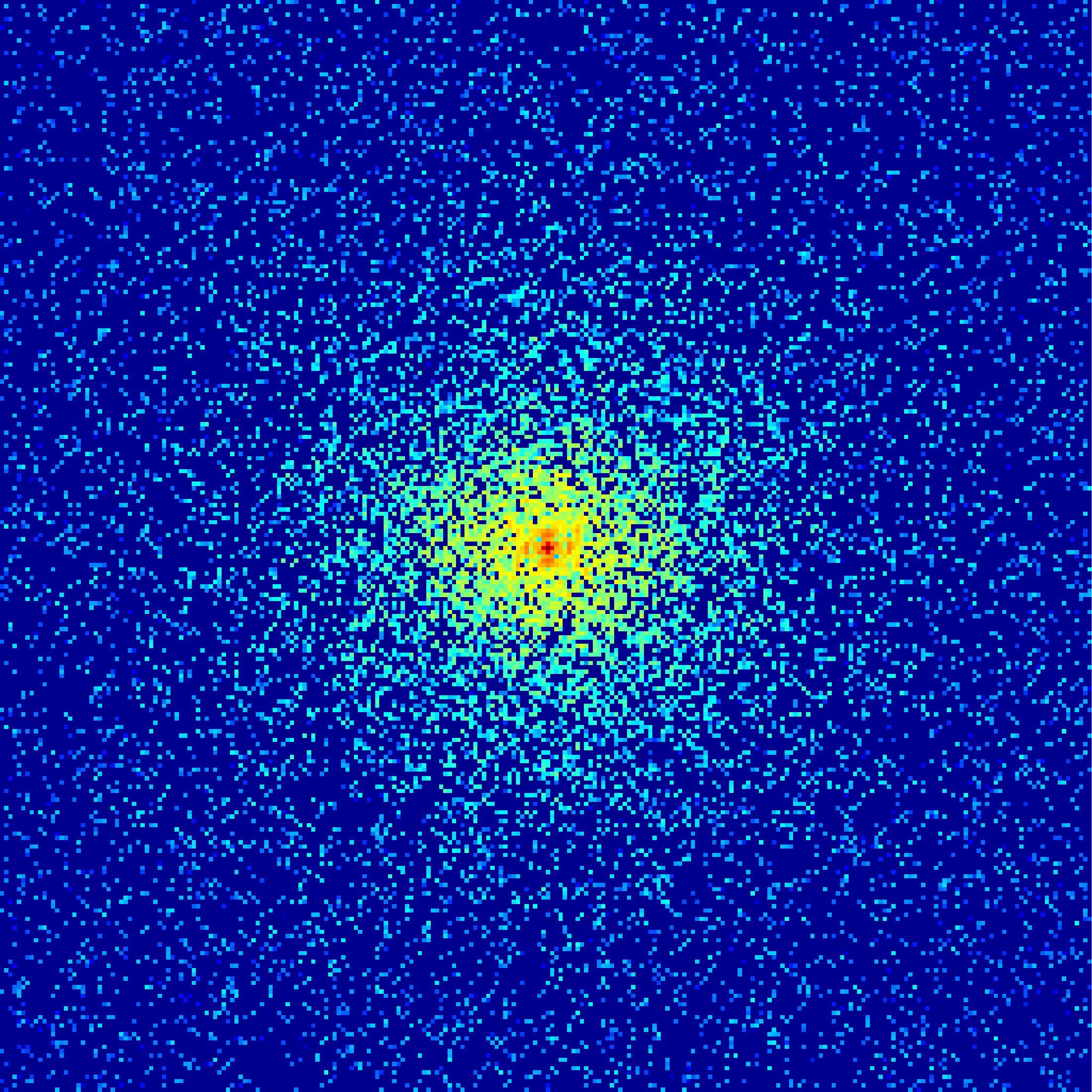}}\vspace{-0.20cm}
\caption{Dataset for the real MR experiments. Fully sampled k-space data, its inverse DFT as a ground truth, the undersampling mask, and the undersampled k-space data.}\label{RealDataSet}
\end{figure}

The real MR image experiments use the k-space data which is obtained from a fully sampled $4$-coil acquisition, and then compressed into a single virtual coil using the SVD technique in \cite{T.Zhang2013}. Since the data from the single virtual coil is complex-valued in the image domain with smoothly varying phase, we further correct the phase using the method described in \cite{G.Ongie2016}. More concretely, we first perform the inverse DFT of the zero padded k-space data, canceling out the phase in the image domain, and passing back to the frequency domain. Then as in the phantom experiments, we generate $20\%$ undersampled k-space data using the variable density sampling, and further add the complex white Gaussian noise so that the resulting SNR of the samples is approximately $25\mathrm{dB}$; see \cref{RealDataSet}.

\begin{table}[tp!]
\centering
\caption{Comparison of signal-to-noise ratio, high frequency error norm, iteration numbers, and CPU time for the real MR experiments. MaxIter denotes that the iteration reached the maximum number without satisfying the stopping criterion.}\label{RealTable}
\vspace{-0.20cm}
\begin{tabular}{|c|c|c|c|c|c|c|c|c|}
\hline
\multirow{2}{*}{Indices}&\multirow{2}{*}{Zero fill}&\multirow{2}{*}{TV \cref{TVModel}}&\multirow{2}{*}{Haar \cref{HaarModel}}&\multirow{2}{*}{TLMRI \cref{TLMRI}}&\multicolumn{3}{|c|}{GIRAF \cref{Schatten}}&\multirow{2}{*}{DDTF \cref{ProposedCSMRIModel}}\\ \cline{6-8}
&&&&&$p=0$&$p=0.5$&$p=1$&\\ \hline\hline
SNR&$11.09$&$14.27$&$14.31$&$\textbf{16.36}$&$15.61$&$15.50$&$13.55$&$16.22$\\ \hline
HFEN&$0.6050$&$0.3477$&$0.3393$&$\textbf{0.2141}$&$0.2788$&$0.2874$&$0.4200$&$0.2456$\\ \hline
NIters&$\cdot$&$277$&$364$&MaxIter&$16$&$18$&$18$&$382$\\ \hline
Time($\mathrm{s}$)&$\cdot$&$3.3$&$22.7$&$50097.5$&$75.3$&$104.3$&$98.3$&$16657.4$\\ \hline
\end{tabular}
\end{table}

\begin{figure}[tp!]
\centering
\hspace{-0.1cm}\subfloat[Fully sampled]{\label{RealOriginal2}\includegraphics[width=3.00cm]{RealOriginal.pdf}}\hspace{0.005cm}
\subfloat[Zero fill]{\label{RealZeroPad}\includegraphics[width=3.00cm]{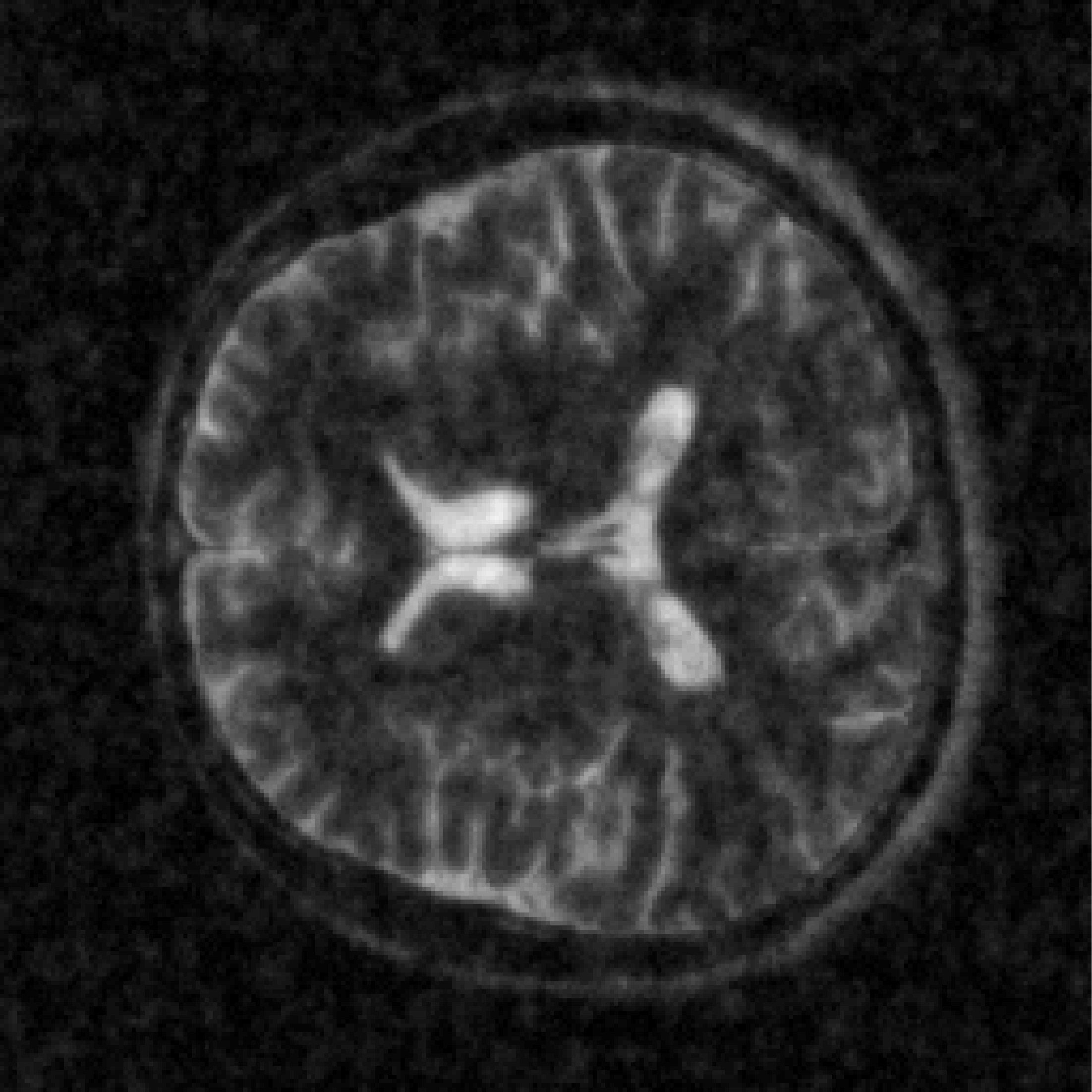}}\hspace{0.005cm}
\subfloat[TV \cref{TVModel}]{\label{RealTV}\includegraphics[width=3.00cm]{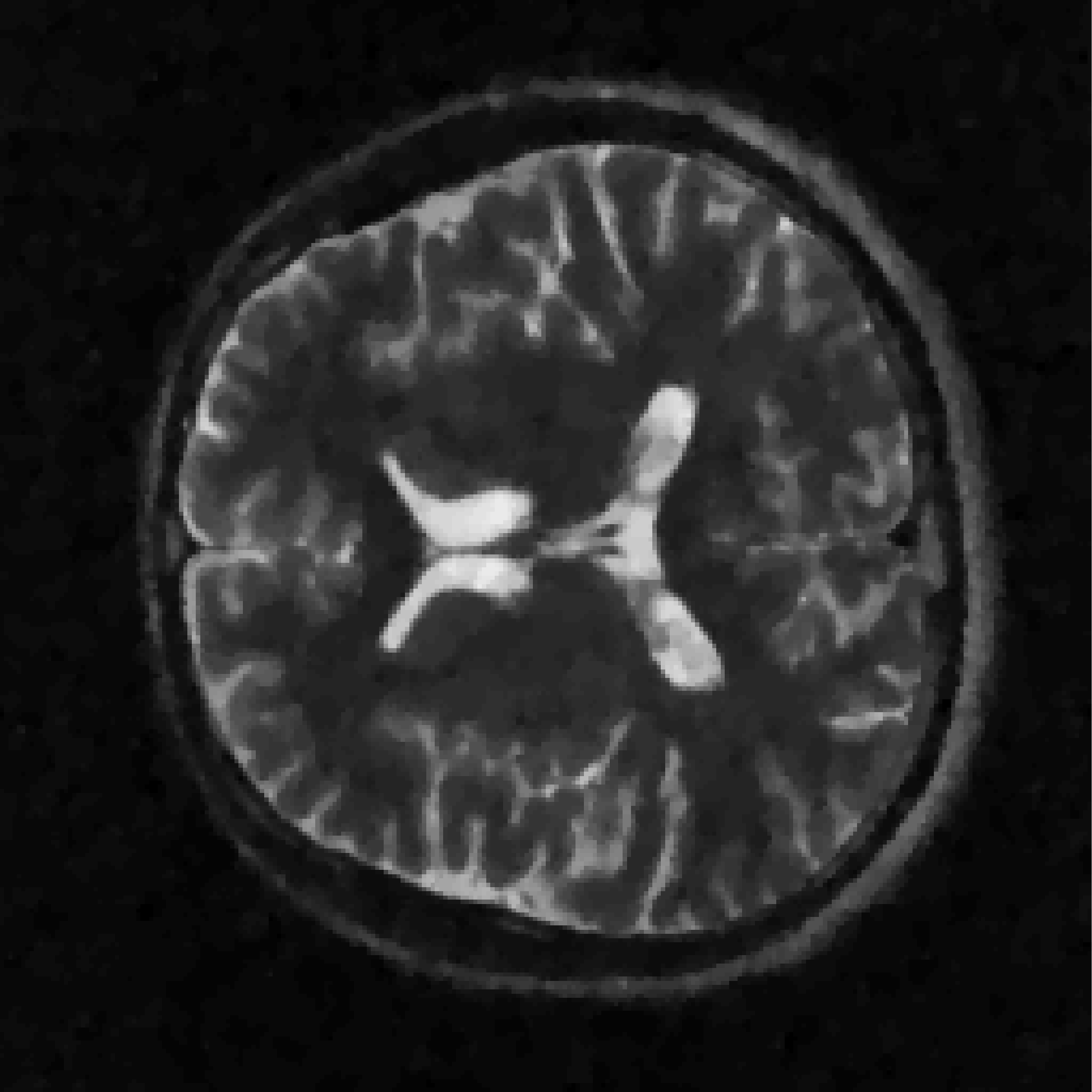}}\hspace{0.005cm}
\subfloat[Haar \cref{HaarModel}]{\label{RealHaar}\includegraphics[width=3.00cm]{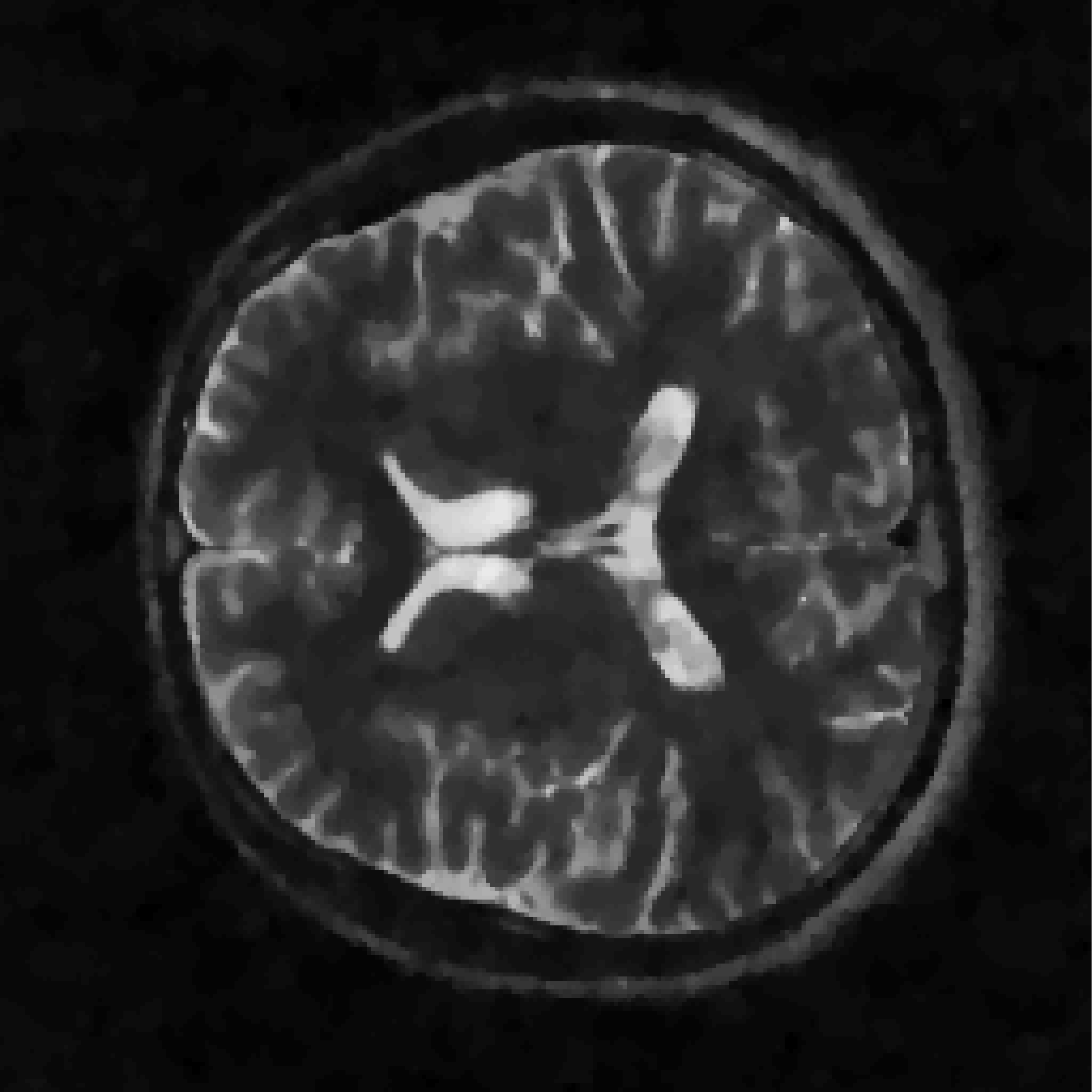}}\hspace{0.005cm}
\subfloat[TLMRI \cref{TLMRI}]{\label{RealTLMRI}\includegraphics[width=3.00cm]{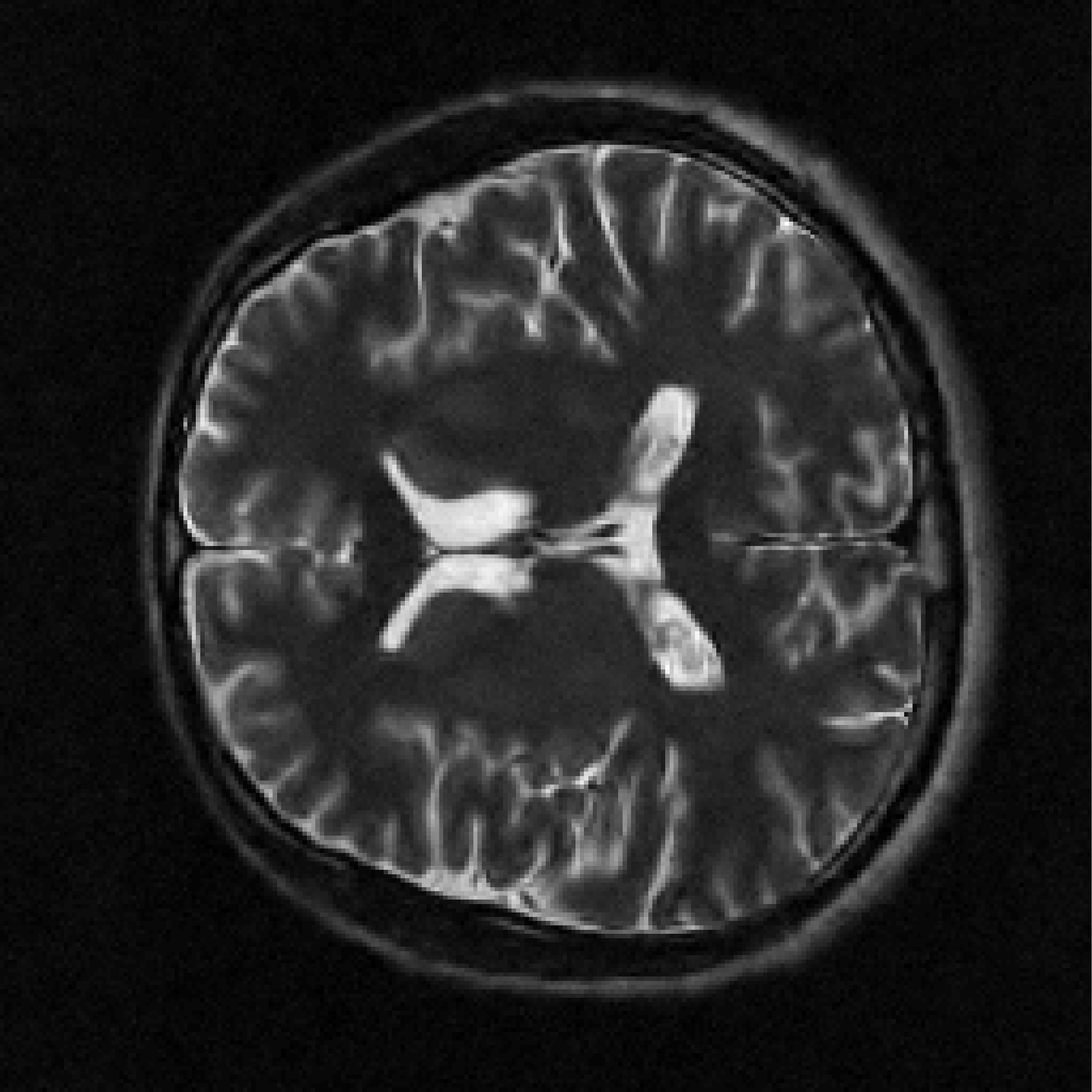}}\vspace{-0.20cm}\\
\subfloat[GIRAF$0$]{\label{RealGIRAF0}\includegraphics[width=3.00cm]{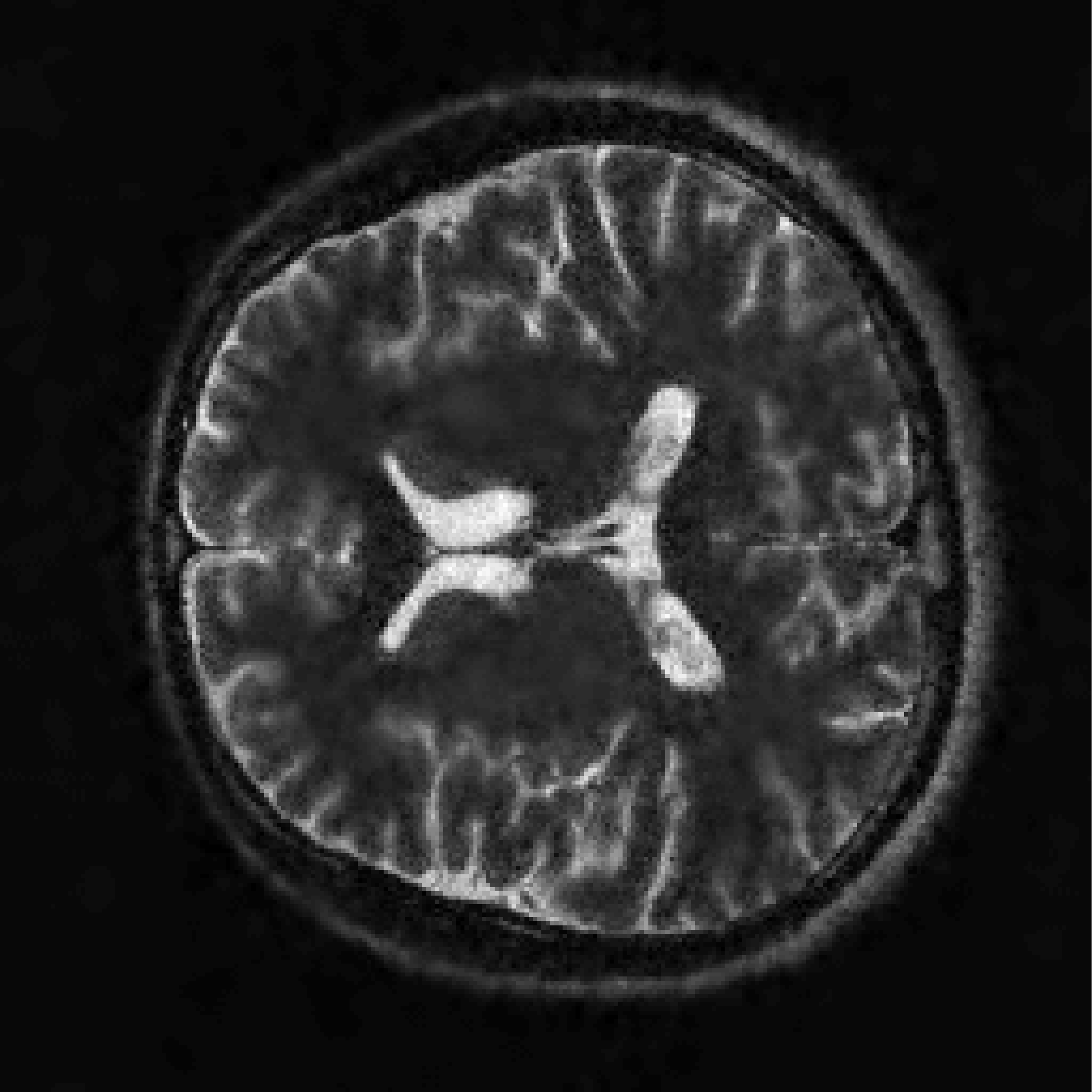}}\hspace{0.005cm}
\subfloat[GIRAF$0.5$]{\label{RealGIRAFHalf}\includegraphics[width=3.00cm]{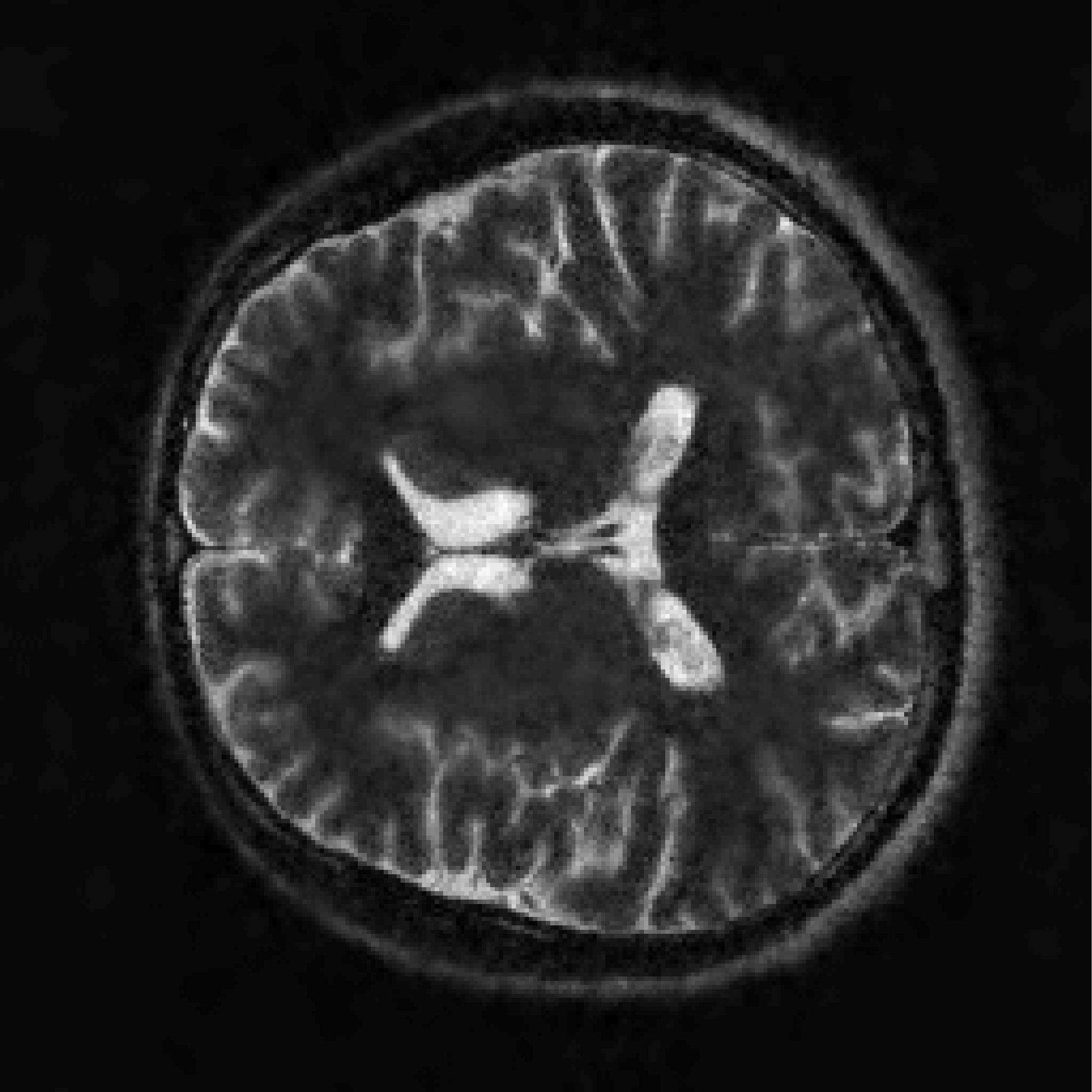}}\hspace{0.005cm}
\subfloat[GIRAF$1$]{\label{RealGIRAF1}\includegraphics[width=3.00cm]{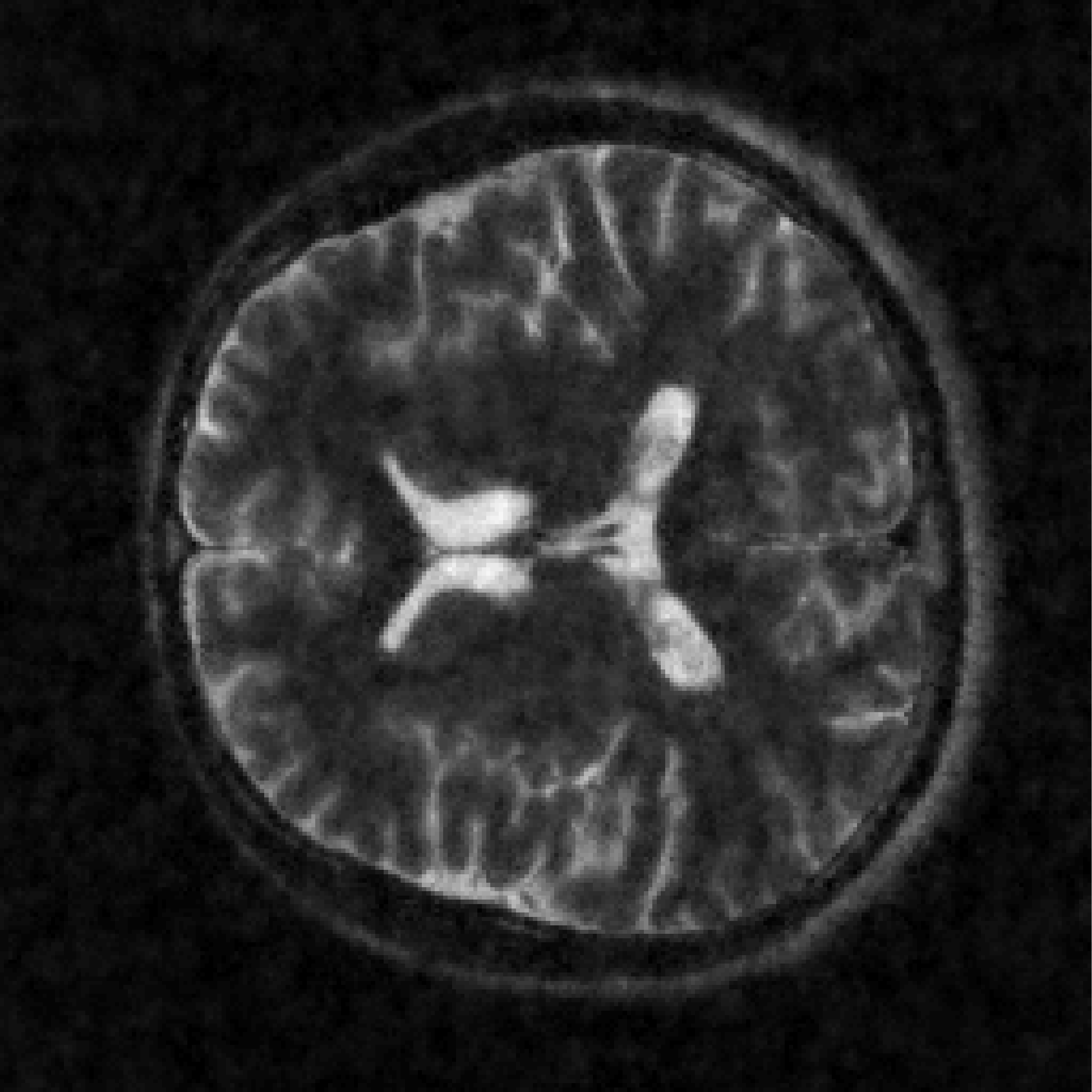}}\hspace{0.005cm}
\subfloat[DDTF \cref{ProposedCSMRIModel}]{\label{RealDDTF}\includegraphics[width=3.00cm]{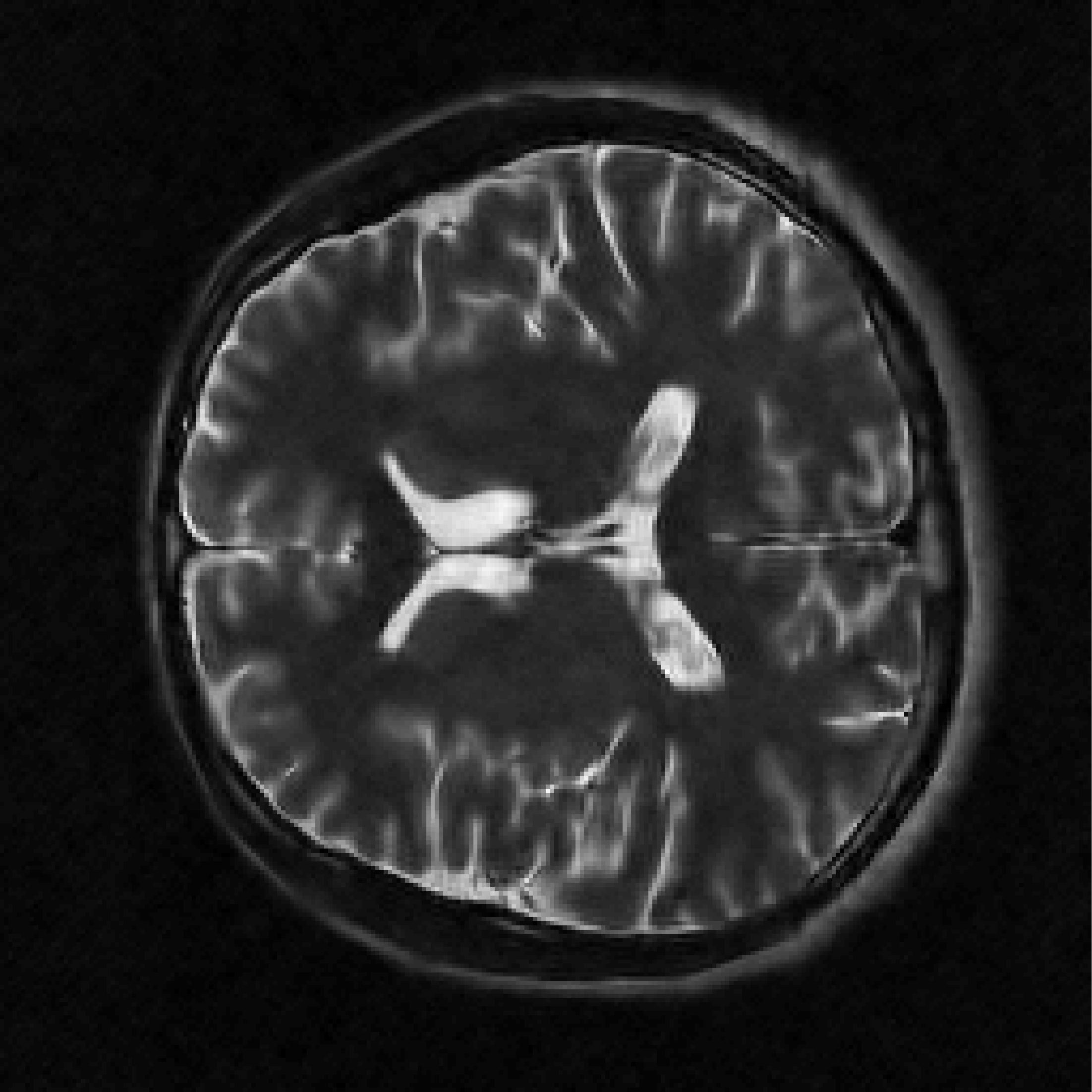}}\vspace{-0.20cm}
\caption{Visual comparisons of each restoration models for the real MR experiments.}\label{RealResults}
\end{figure}

\begin{figure}[tp!]
\centering
\hspace{-0.1cm}\subfloat[Sample mask]{\label{RealMask2}\includegraphics[width=3.00cm]{RealMask.pdf}}\hspace{0.005cm}
\subfloat[Zero fill]{\label{RealZeroPadError}\includegraphics[width=3.00cm]{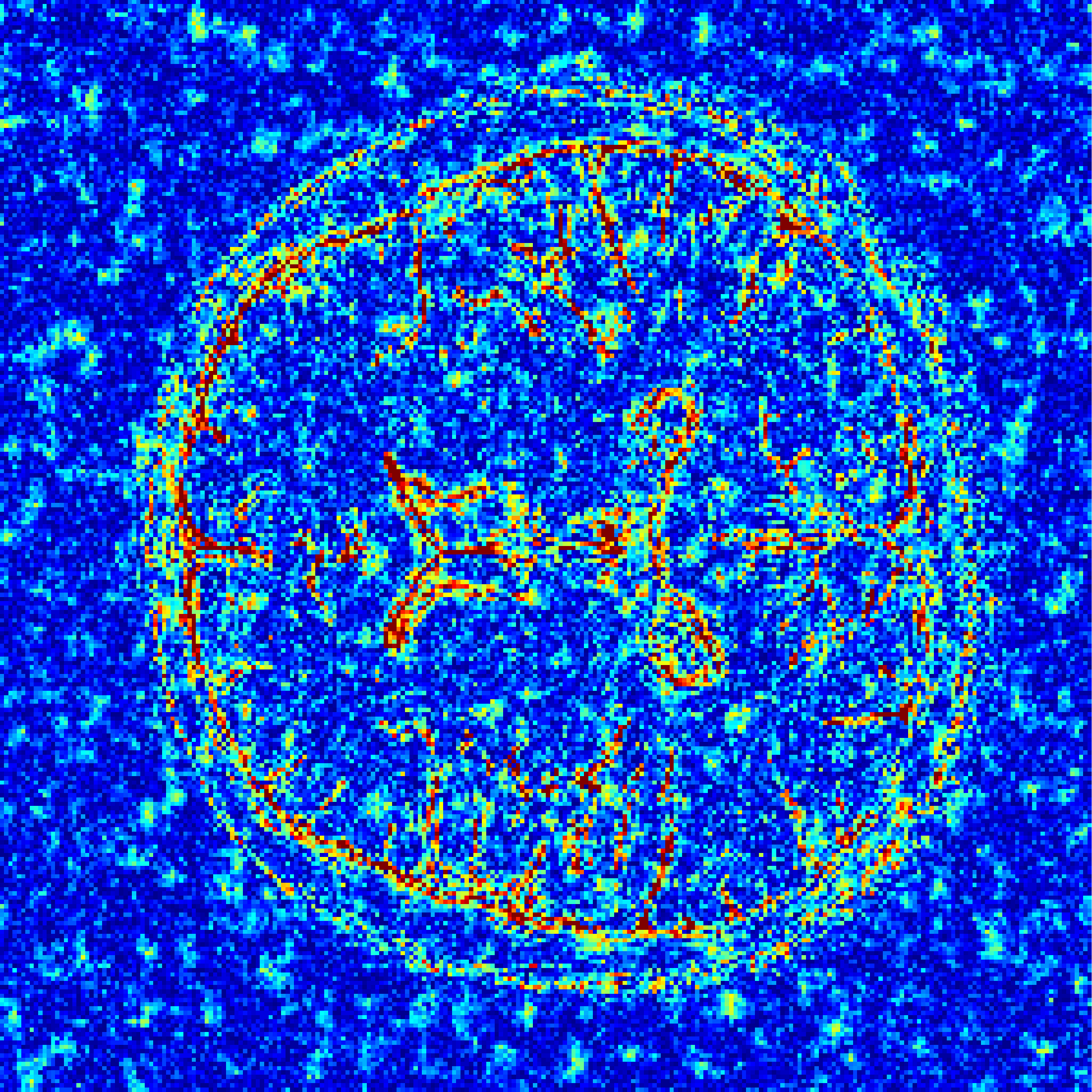}}\hspace{0.005cm}
\subfloat[TV \cref{TVModel}]{\label{RealTVError}\includegraphics[width=3.00cm]{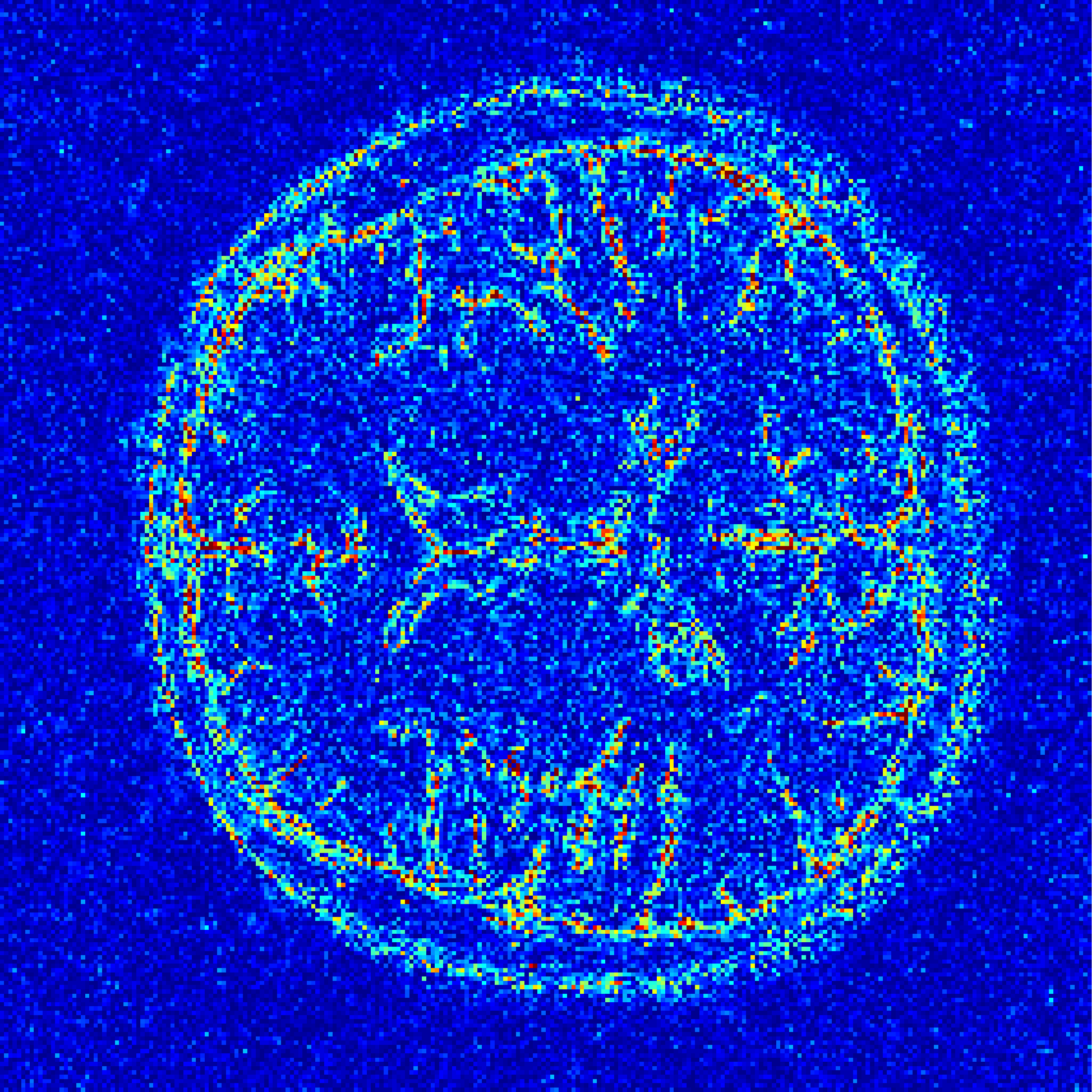}}\hspace{0.005cm}
\subfloat[Haar \cref{HaarModel}]{\label{RealHaarError}\includegraphics[width=3.00cm]{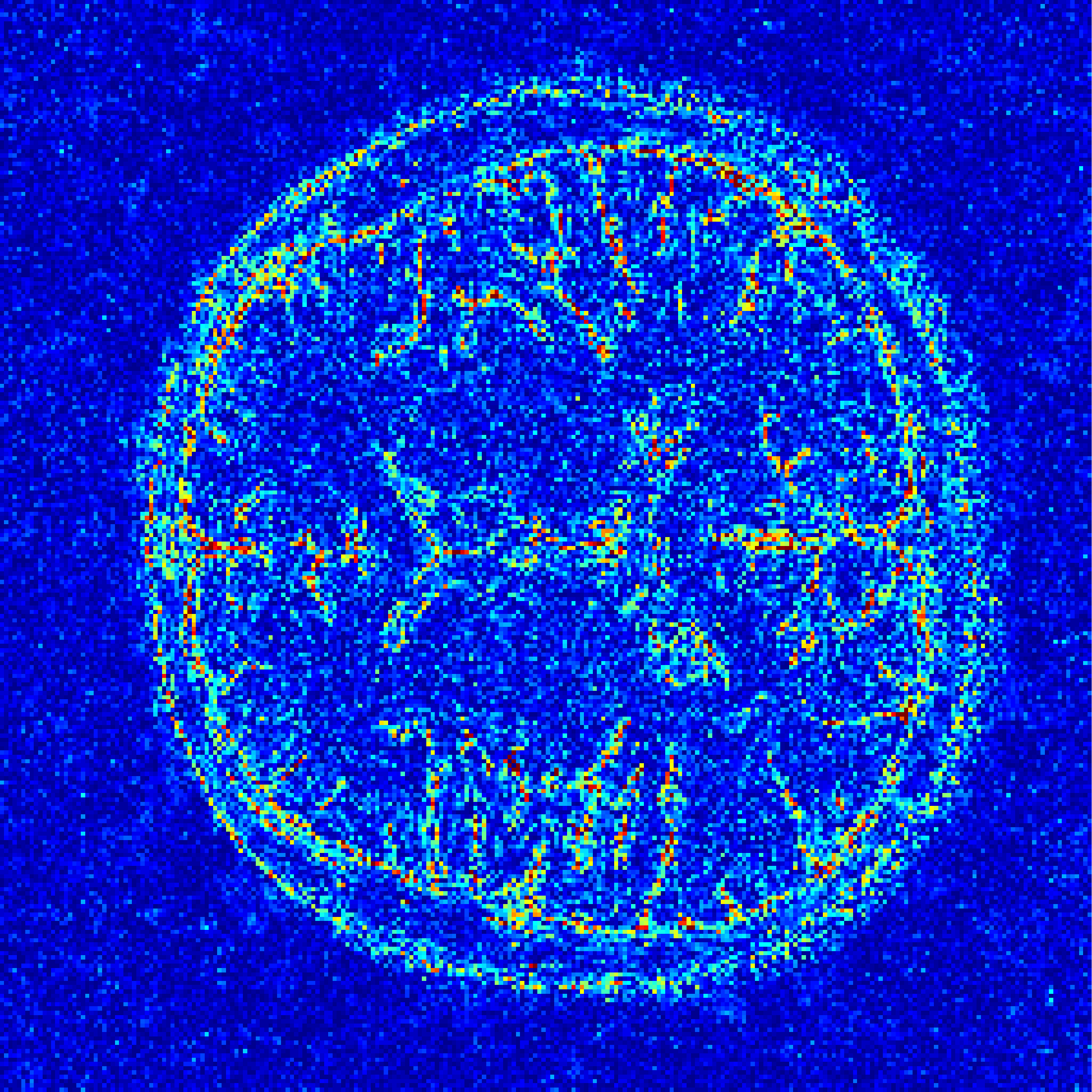}}\hspace{0.005cm}
\subfloat[TLMRI \cref{TLMRI}]{\label{RealTLMRIError}\includegraphics[width=3.00cm]{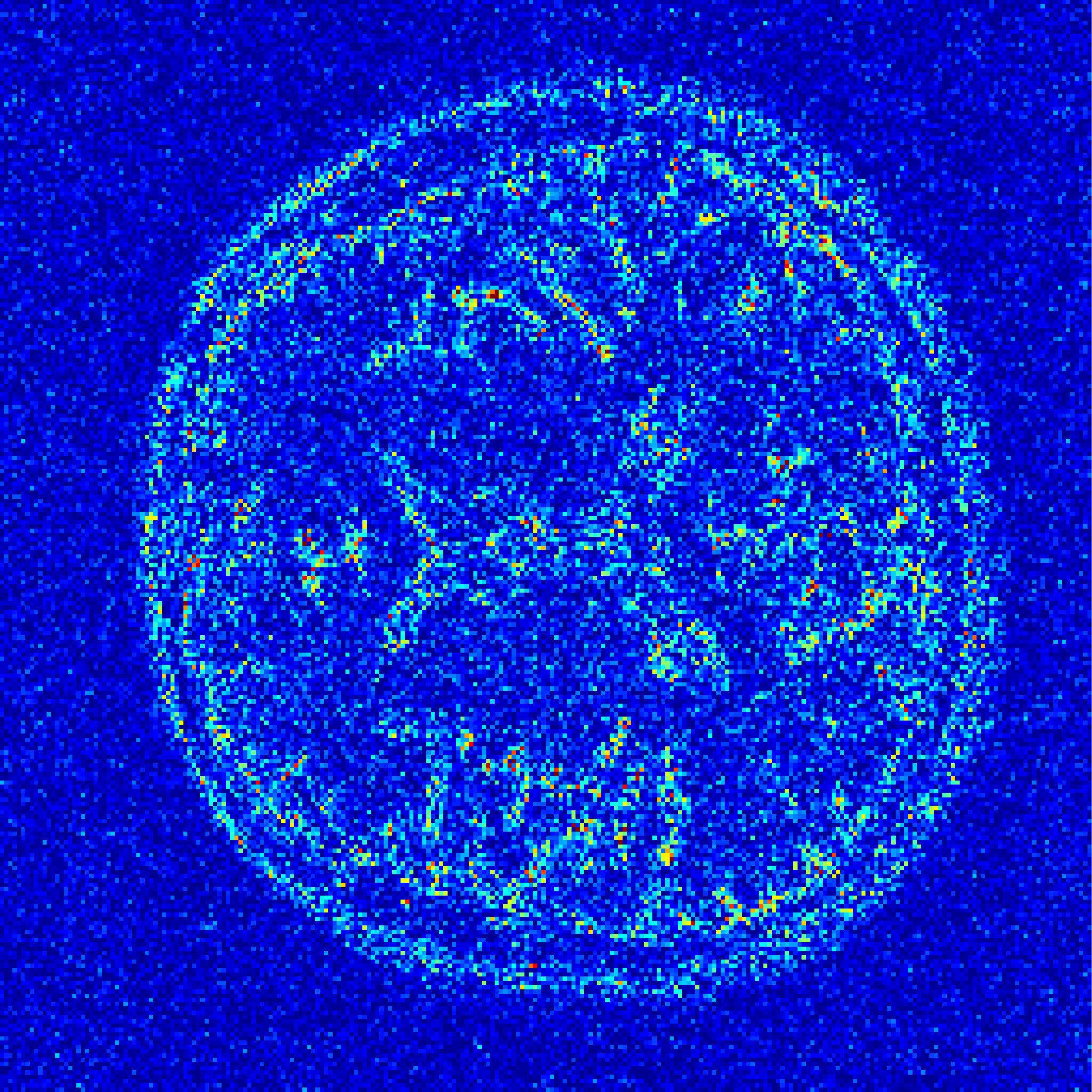}}\vspace{-0.20cm}\\
\subfloat[GIRAF$0$]{\label{RealGIRAF0Error}\includegraphics[width=3.00cm]{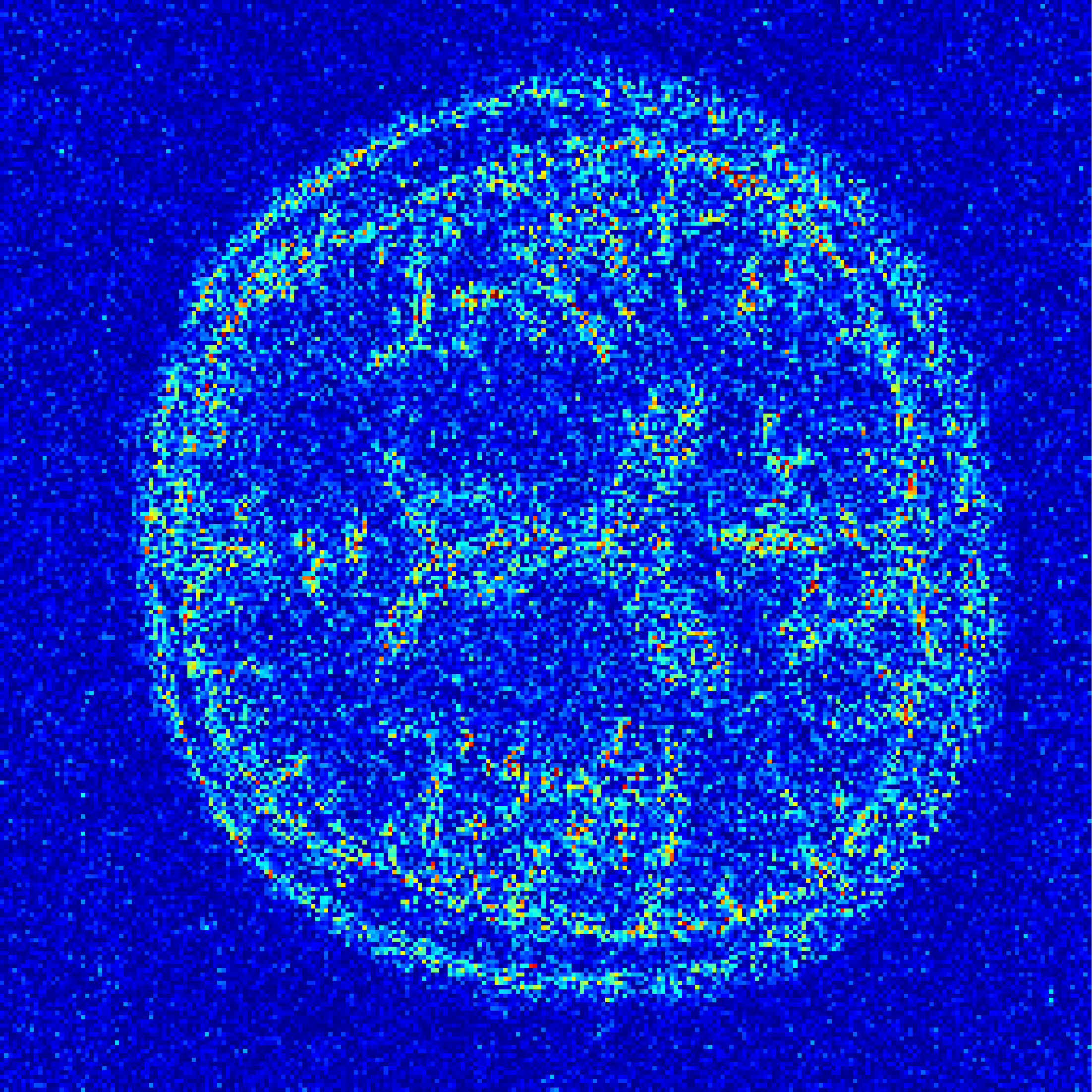}}\hspace{0.005cm}
\subfloat[GIRAF$0.5$]{\label{RealGIRAFHalfError}\includegraphics[width=3.00cm]{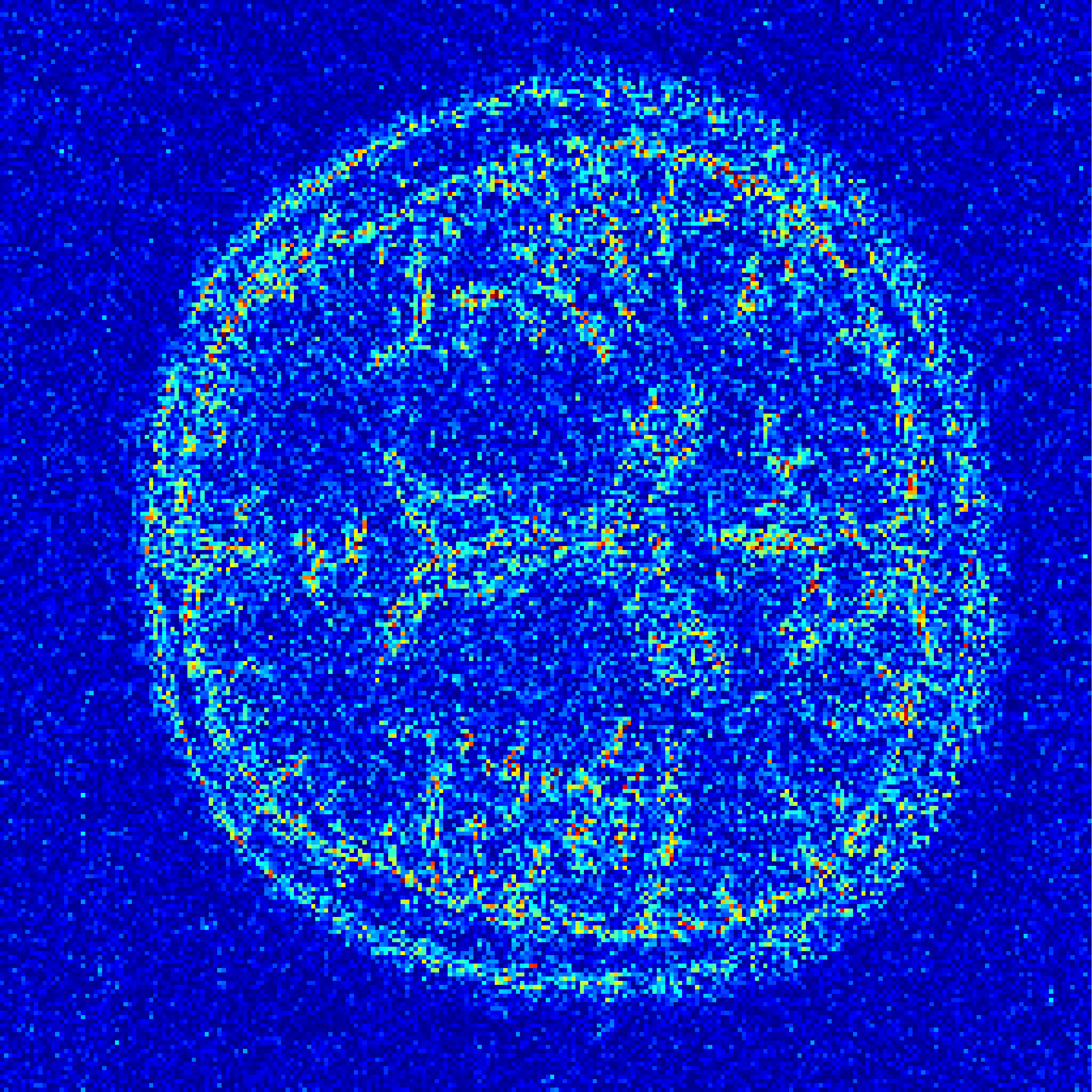}}\hspace{0.005cm}
\subfloat[GIRAF$1$]{\label{RealGIRAF1Error}\includegraphics[width=3.00cm]{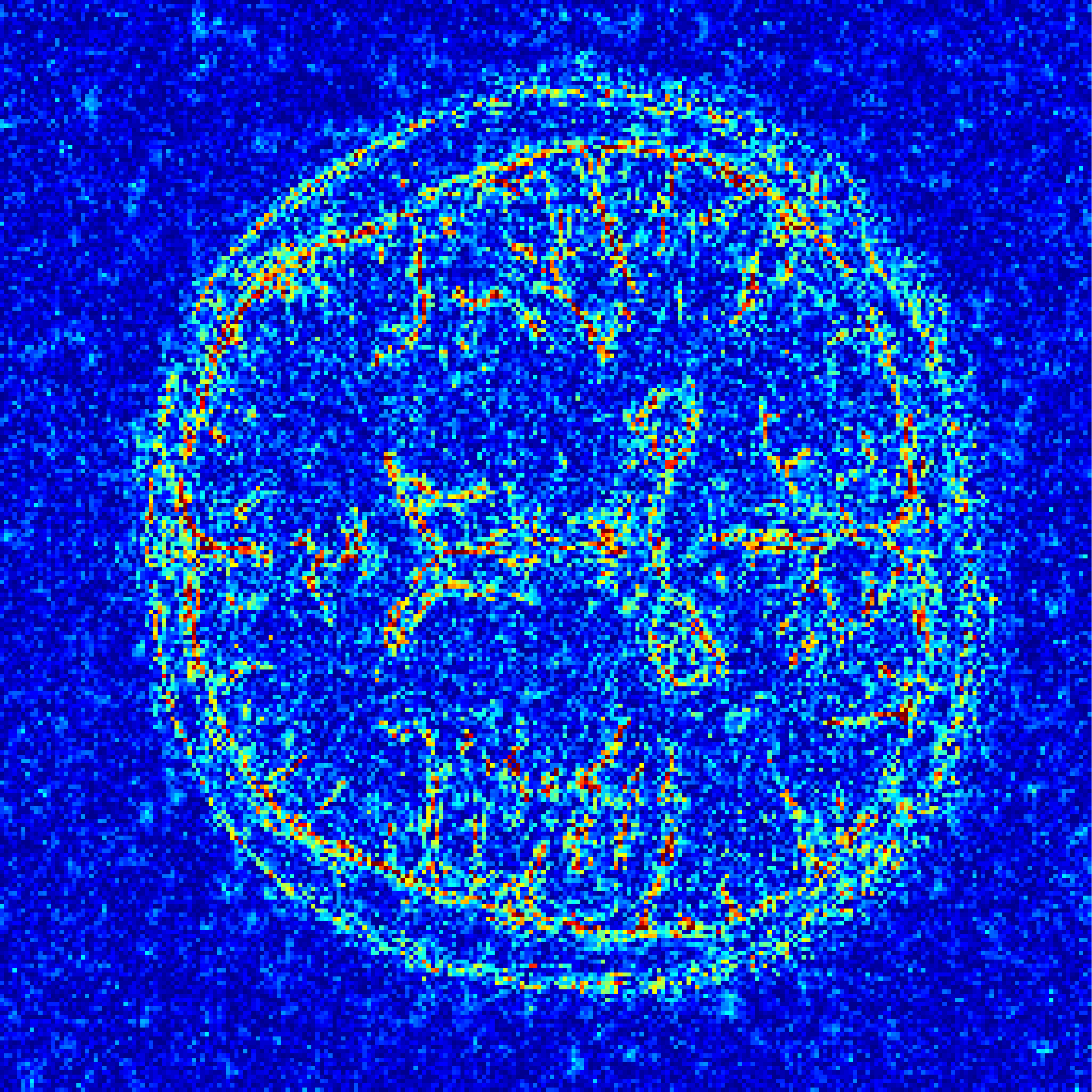}}\hspace{0.005cm}
\subfloat[DDTF \cref{ProposedCSMRIModel}]{\label{RealDDTFError}\includegraphics[width=3.00cm]{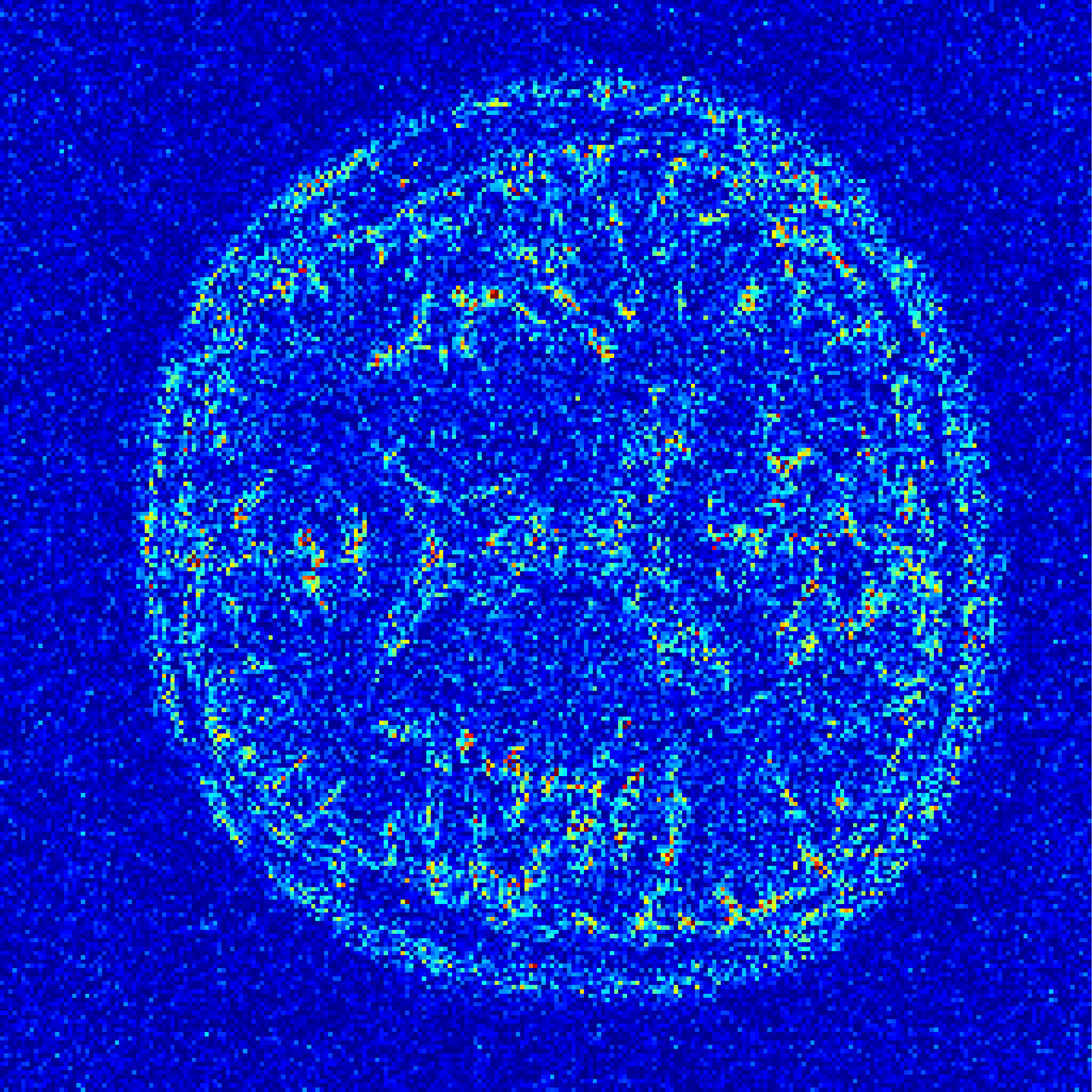}}\vspace{-0.20cm}
\caption{Comparisons of error maps for the phantom experiments.}\label{RealResultsErrorMap}
\end{figure}

\begin{figure}[tp!]
\centering
\hspace{-0.1cm}\subfloat[Fully sampled]{\label{RealOriginalk2}\includegraphics[width=3.00cm]{RealOriginalk.pdf}}\hspace{0.005cm}
\subfloat[Zero fill]{\label{RealUndersamplek}\includegraphics[width=3.00cm]{RealUndersamplek.pdf}}\hspace{0.005cm}
\subfloat[TV \cref{TVModel}]{\label{RealTVk}\includegraphics[width=3.00cm]{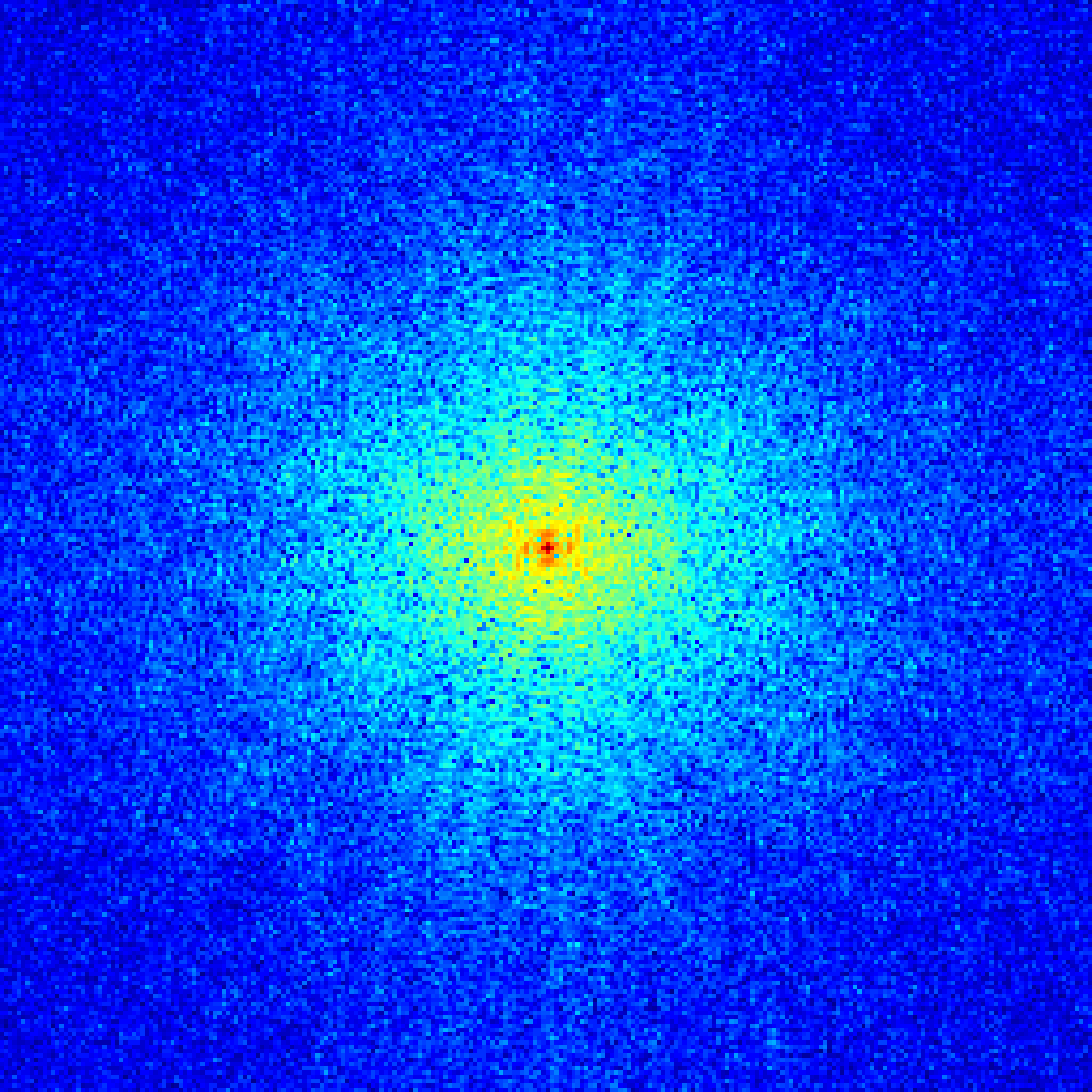}}\hspace{0.005cm}
\subfloat[Haar \cref{HaarModel}]{\label{RealHaark}\includegraphics[width=3.00cm]{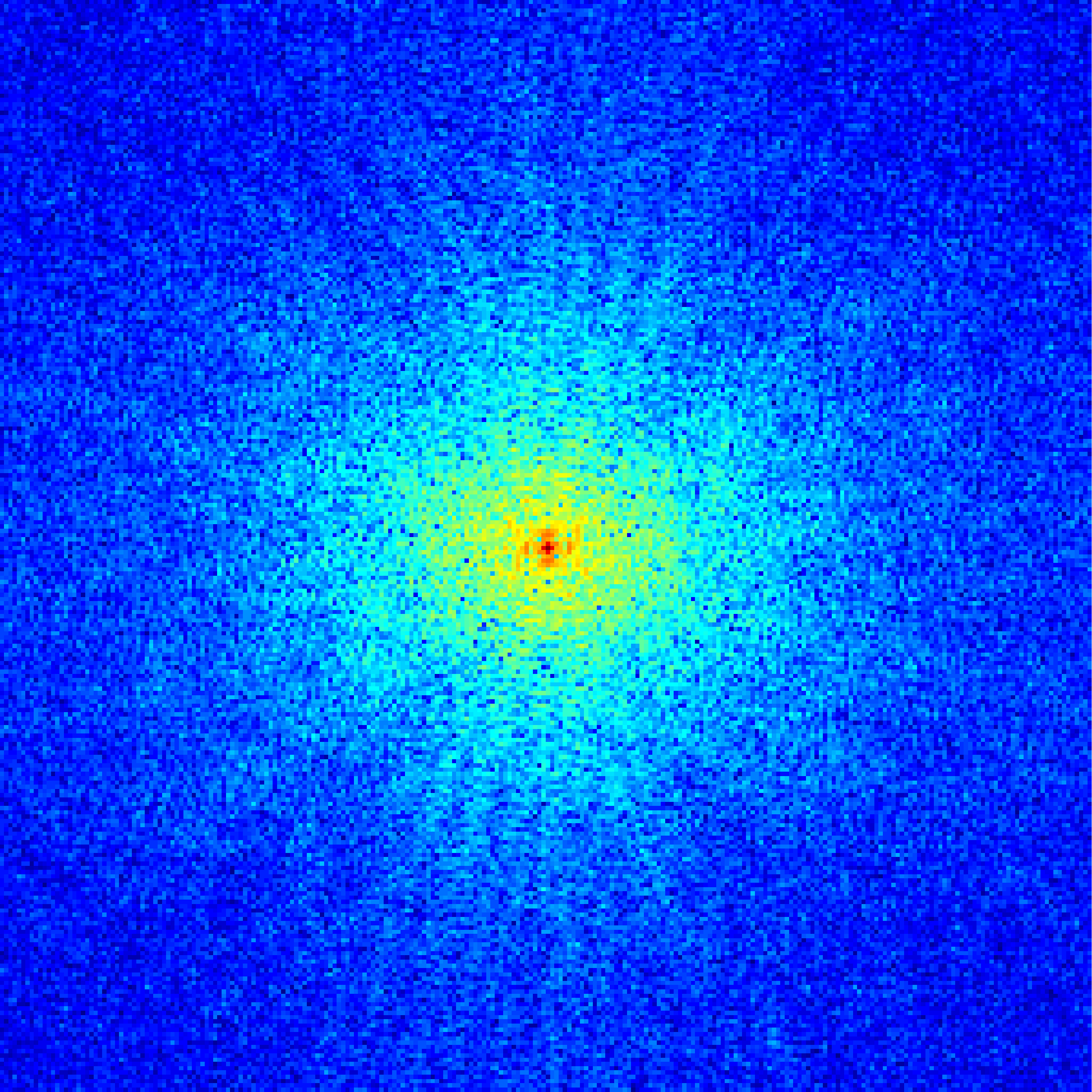}}\hspace{0.005cm}
\subfloat[TLMRI \cref{TLMRI}]{\label{RealTLMRIk}\includegraphics[width=3.00cm]{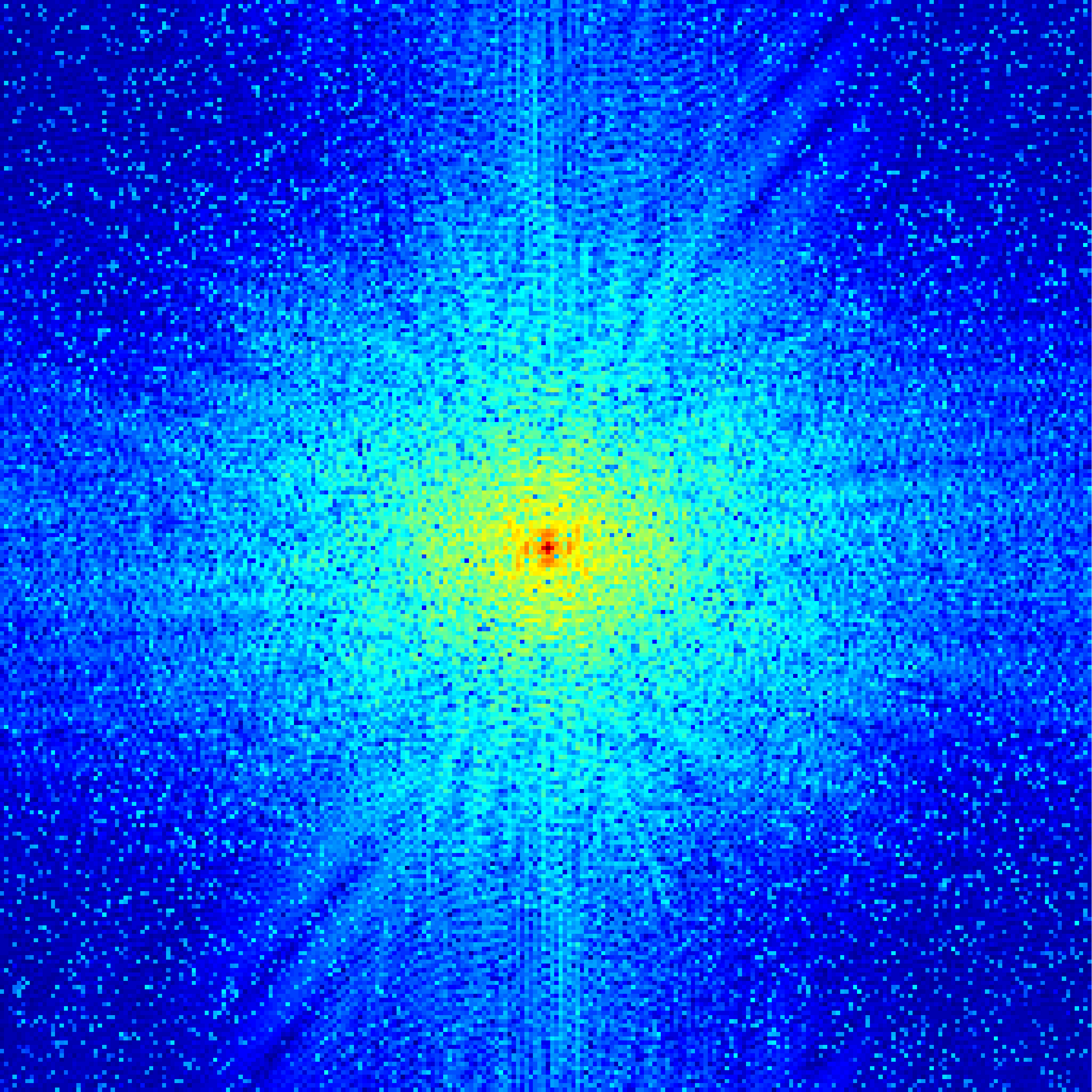}}\vspace{-0.20cm}\\
\subfloat[GIRAF$0$]{\label{RealGIRAF0k}\includegraphics[width=3.00cm]{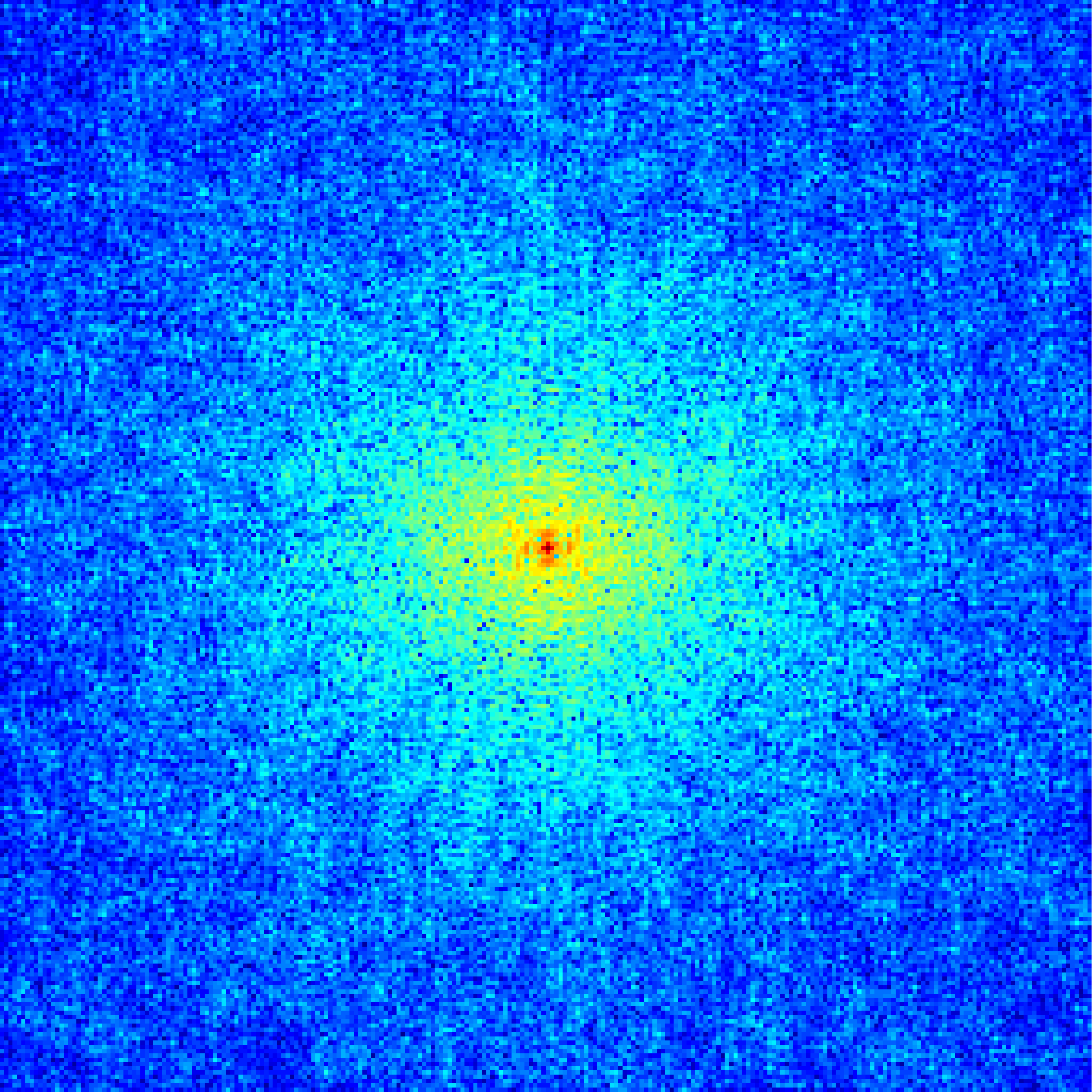}}\hspace{0.005cm}
\subfloat[GIRAF$0.5$]{\label{RealGIRAFHalfk}\includegraphics[width=3.00cm]{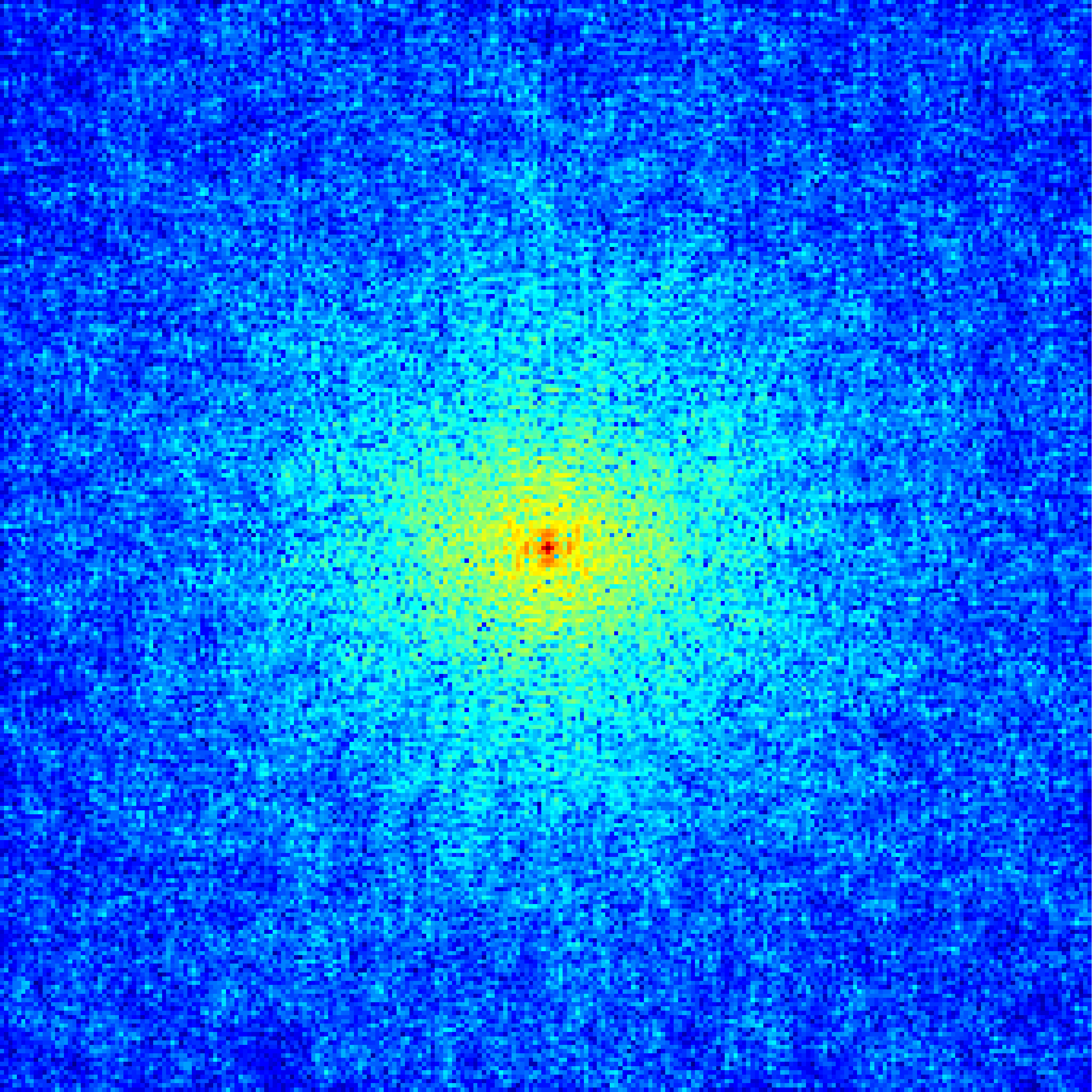}}\hspace{0.005cm}
\subfloat[GIRAF$1$]{\label{RealGIRAF1k}\includegraphics[width=3.00cm]{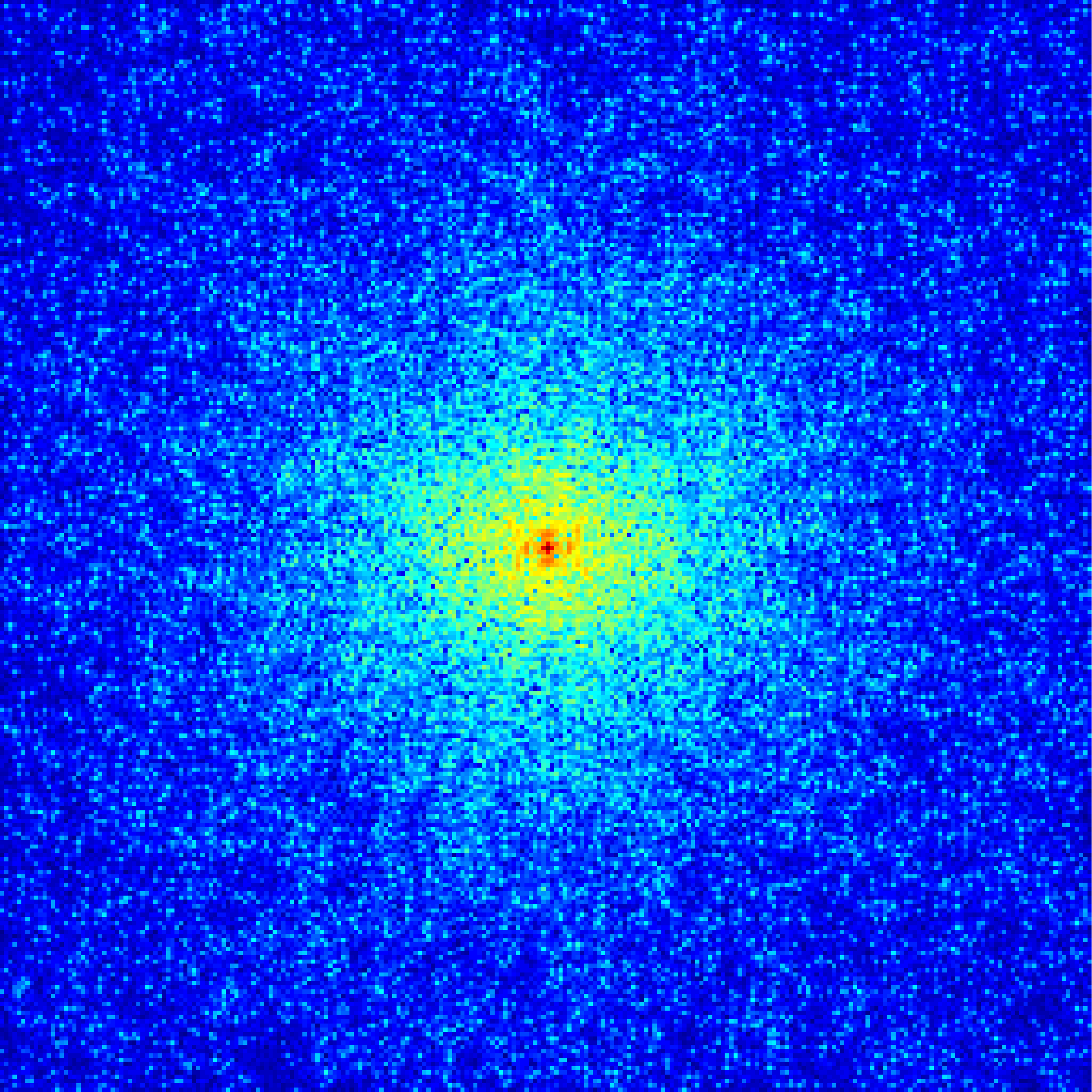}}\hspace{0.005cm}
\subfloat[DDTF \cref{ProposedCSMRIModel}]{\label{RealDDTFk}\includegraphics[width=3.00cm]{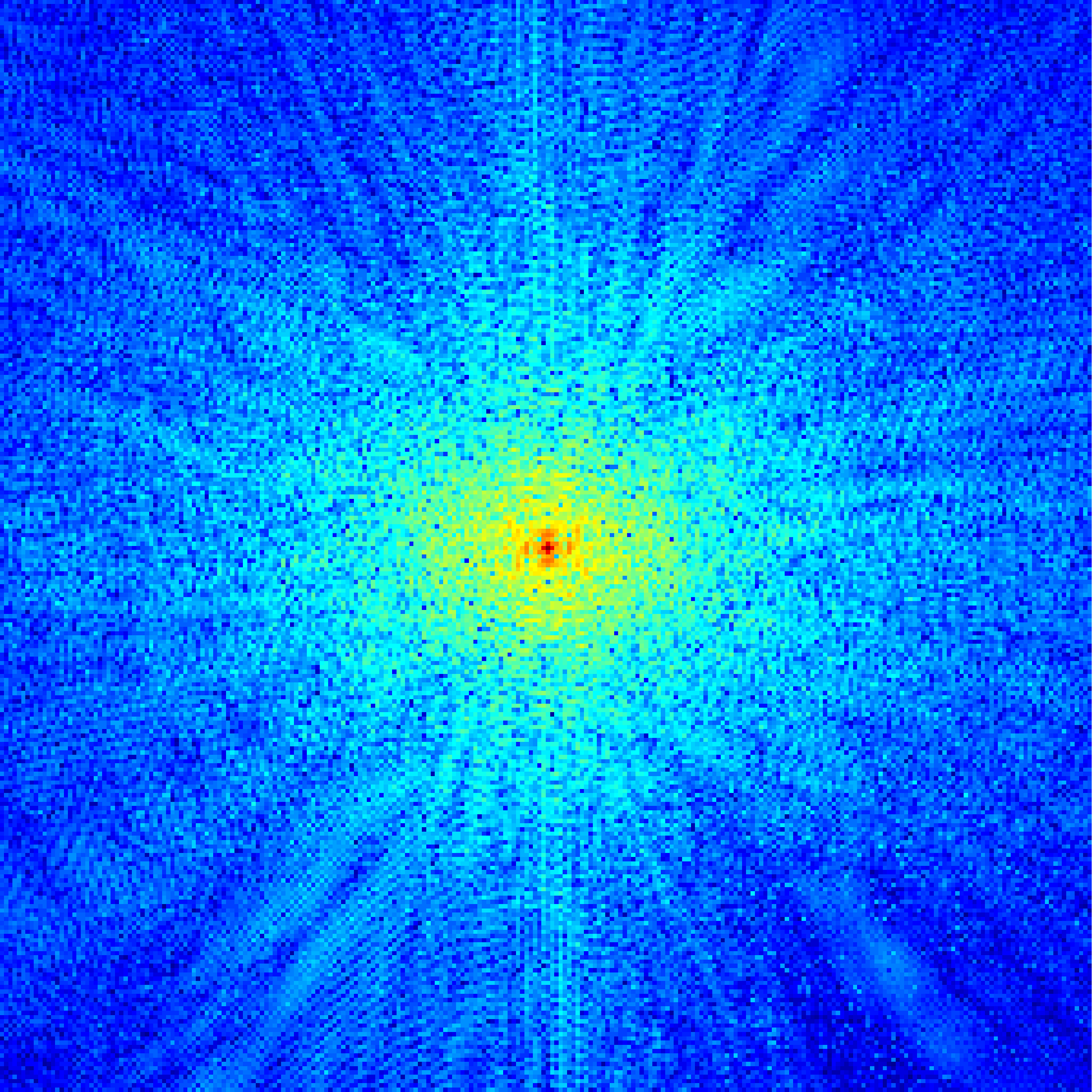}}\vspace{-0.20cm}
\caption{Comparisons of k-space data for the phantom experiments in the log scale. All restored k-space data are displayed in the window level $[0,9]$ for fair comparisons.}\label{RealResultsk}
\end{figure}

\begin{figure}[tp!]
\centering
\hspace{-0.1cm}\subfloat[Fully sampled]{\label{RealOriginalZoom}\includegraphics[width=3.00cm]{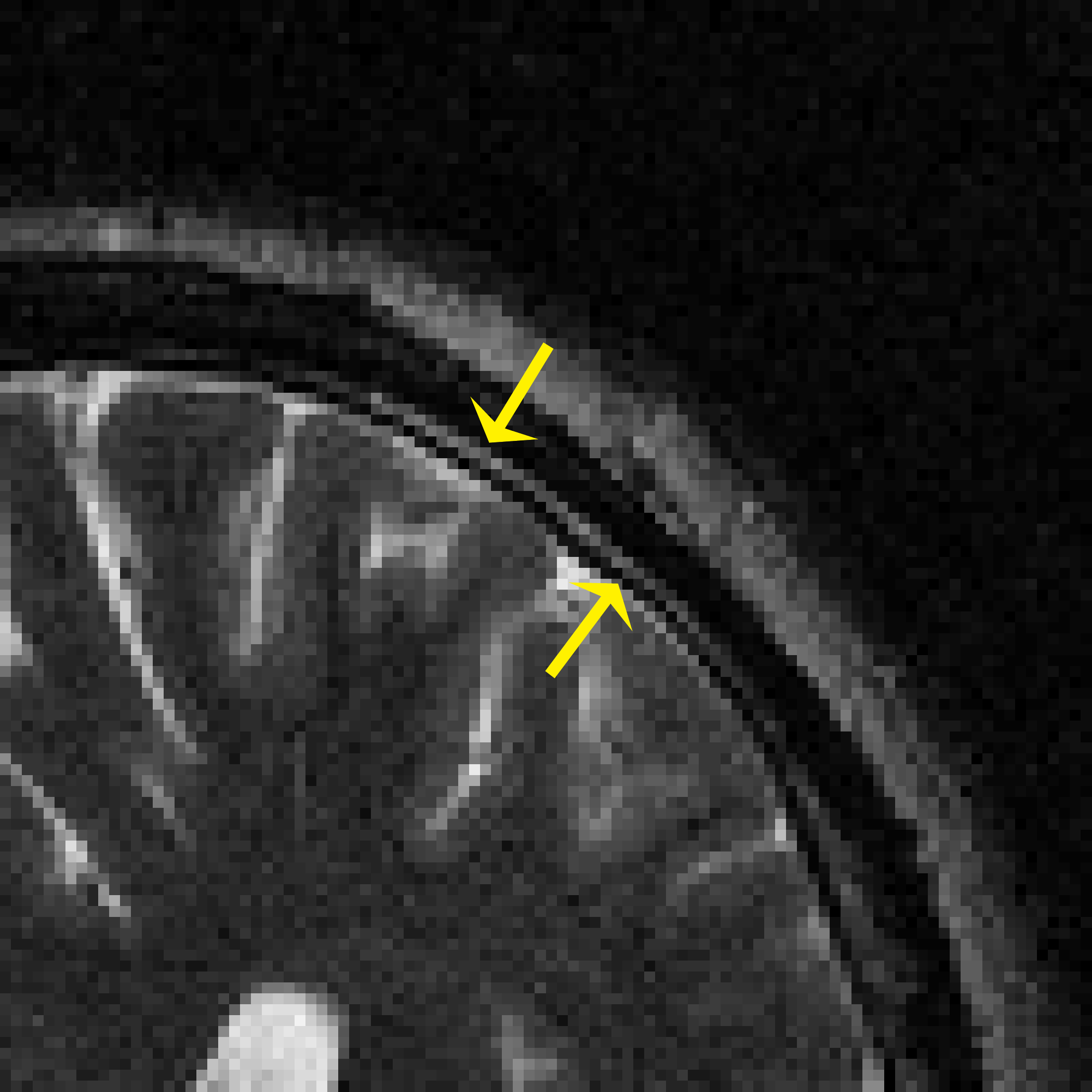}}\hspace{0.005cm}
\subfloat[Zero fill]{\label{RealZeroPadZoom}\includegraphics[width=3.00cm]{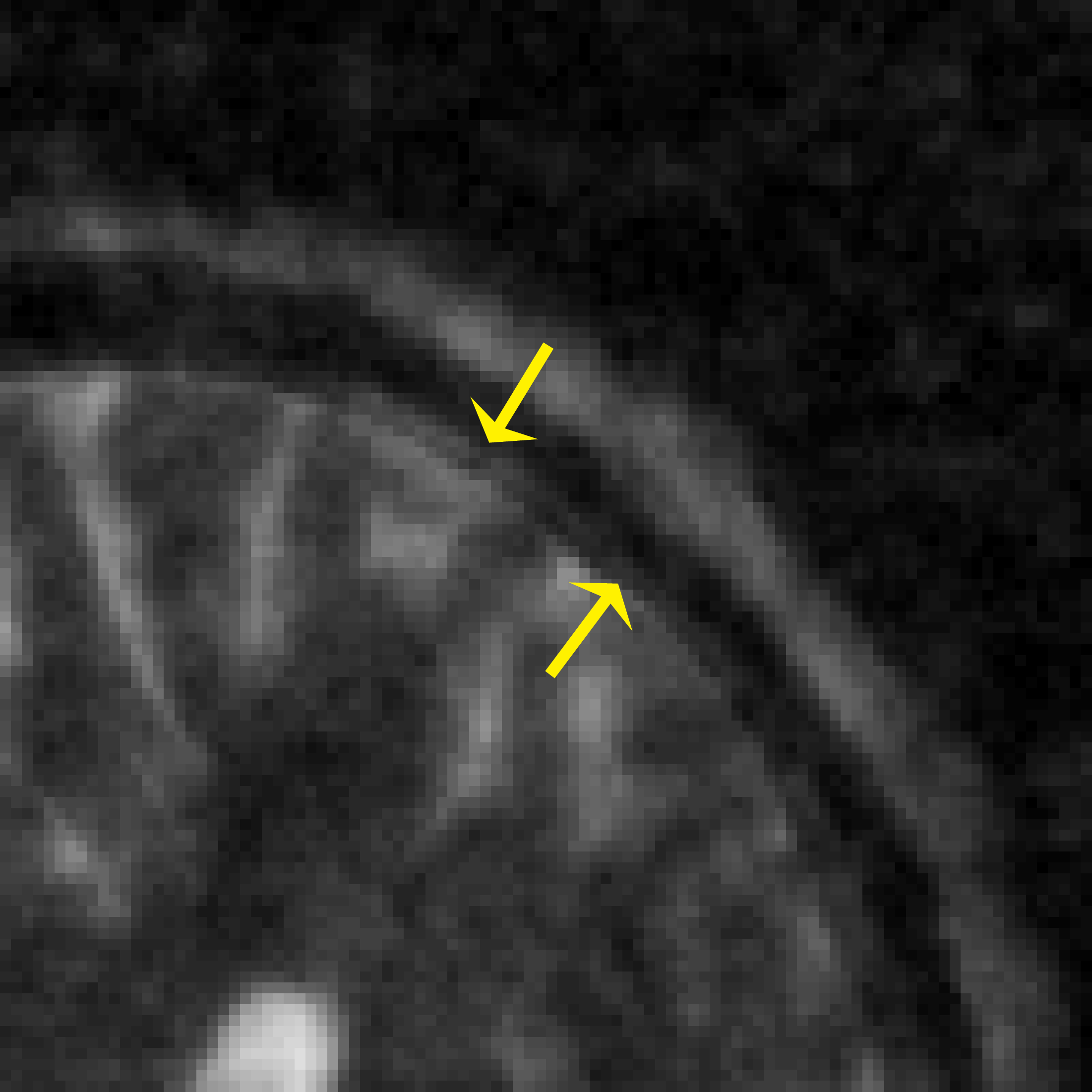}}\hspace{0.005cm}
\subfloat[TV \cref{TVModel}]{\label{RealTVZoom}\includegraphics[width=3.00cm]{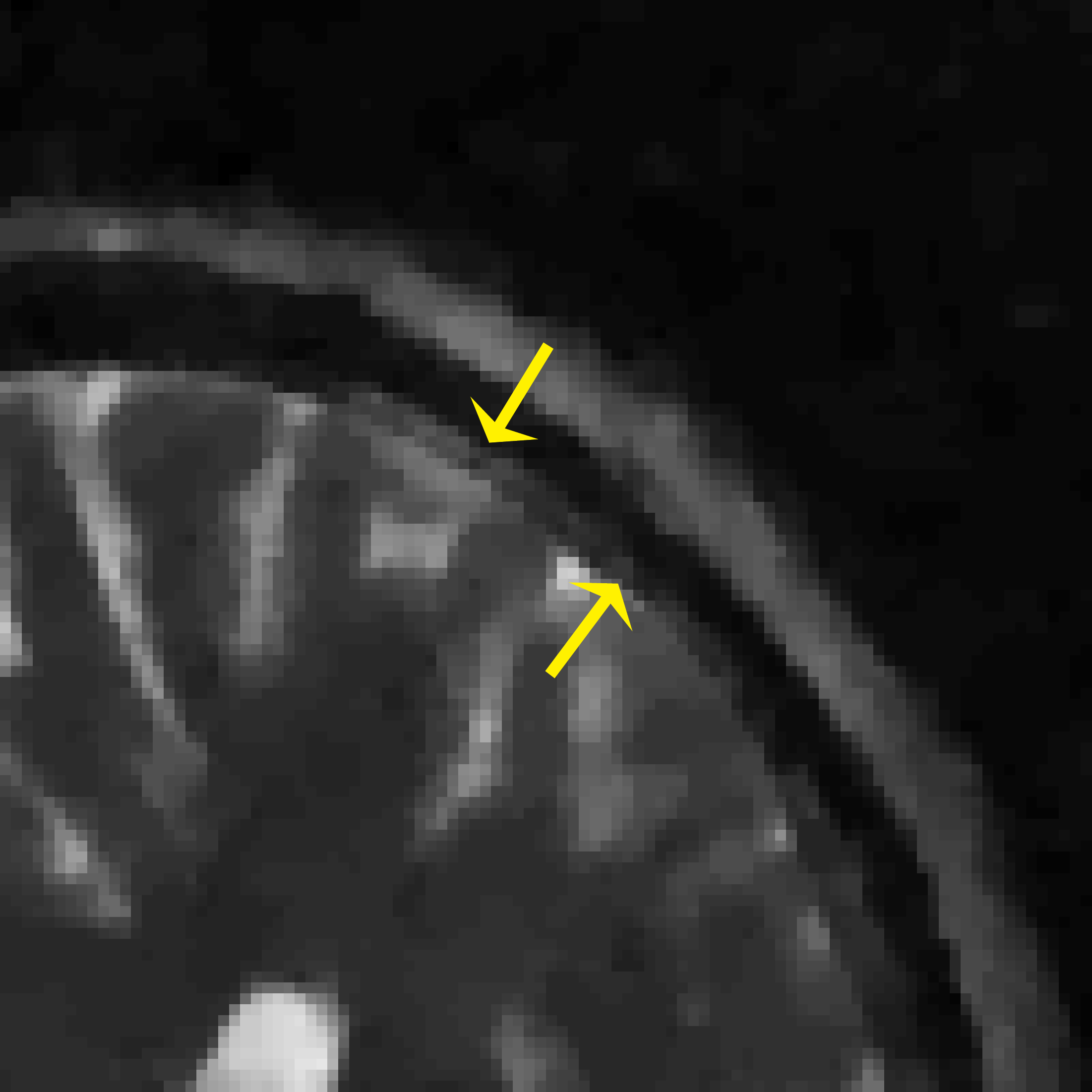}}\hspace{0.005cm}
\subfloat[Haar \cref{HaarModel}]{\label{RealHaarZoom}\includegraphics[width=3.00cm]{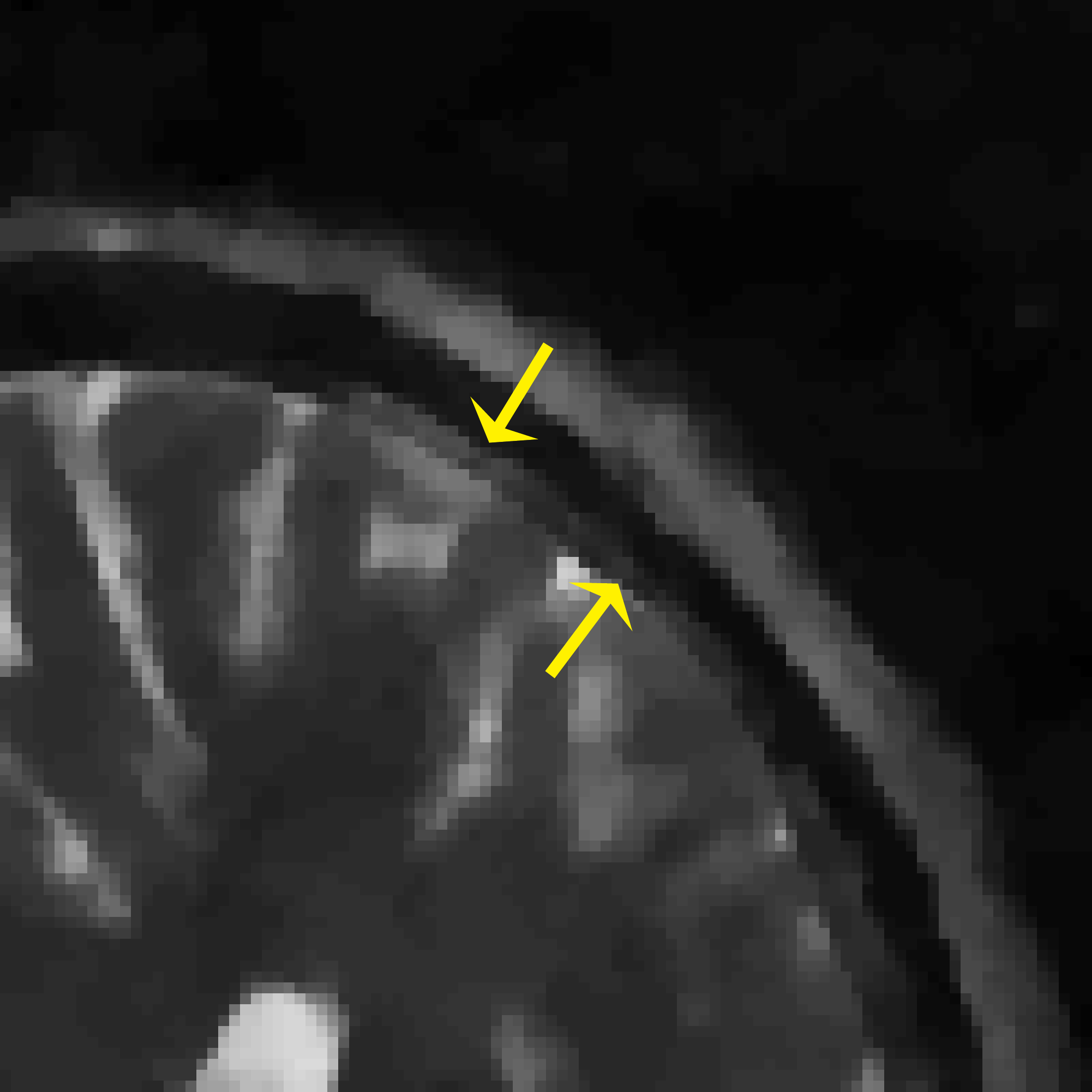}}\hspace{0.005cm}
\subfloat[TLMRI \cref{TLMRI}]{\label{RealTLMRIZoom}\includegraphics[width=3.00cm]{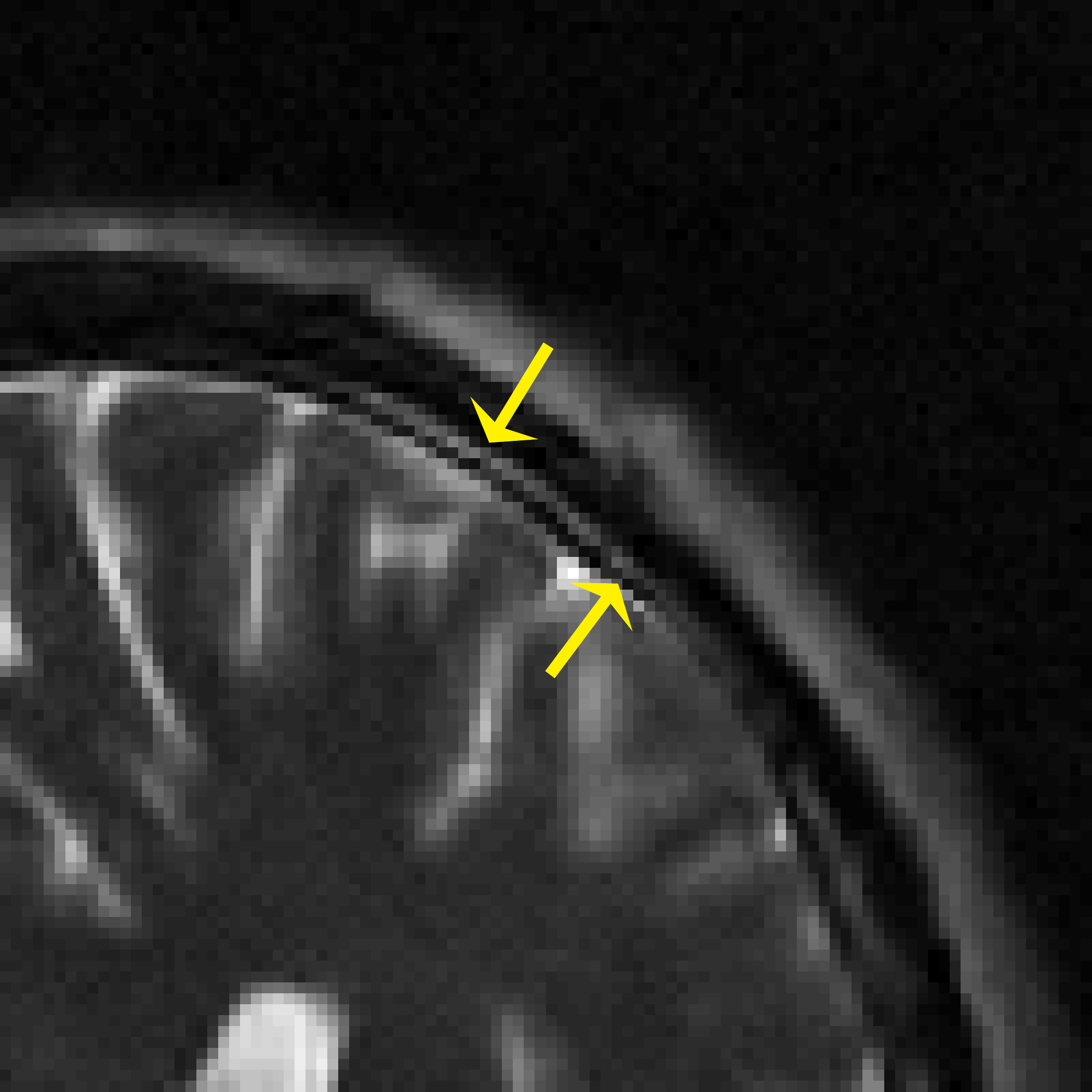}}\vspace{-0.20cm}\\
\subfloat[GIRAF$0$]{\label{RealGIRAF0Zoom}\includegraphics[width=3.00cm]{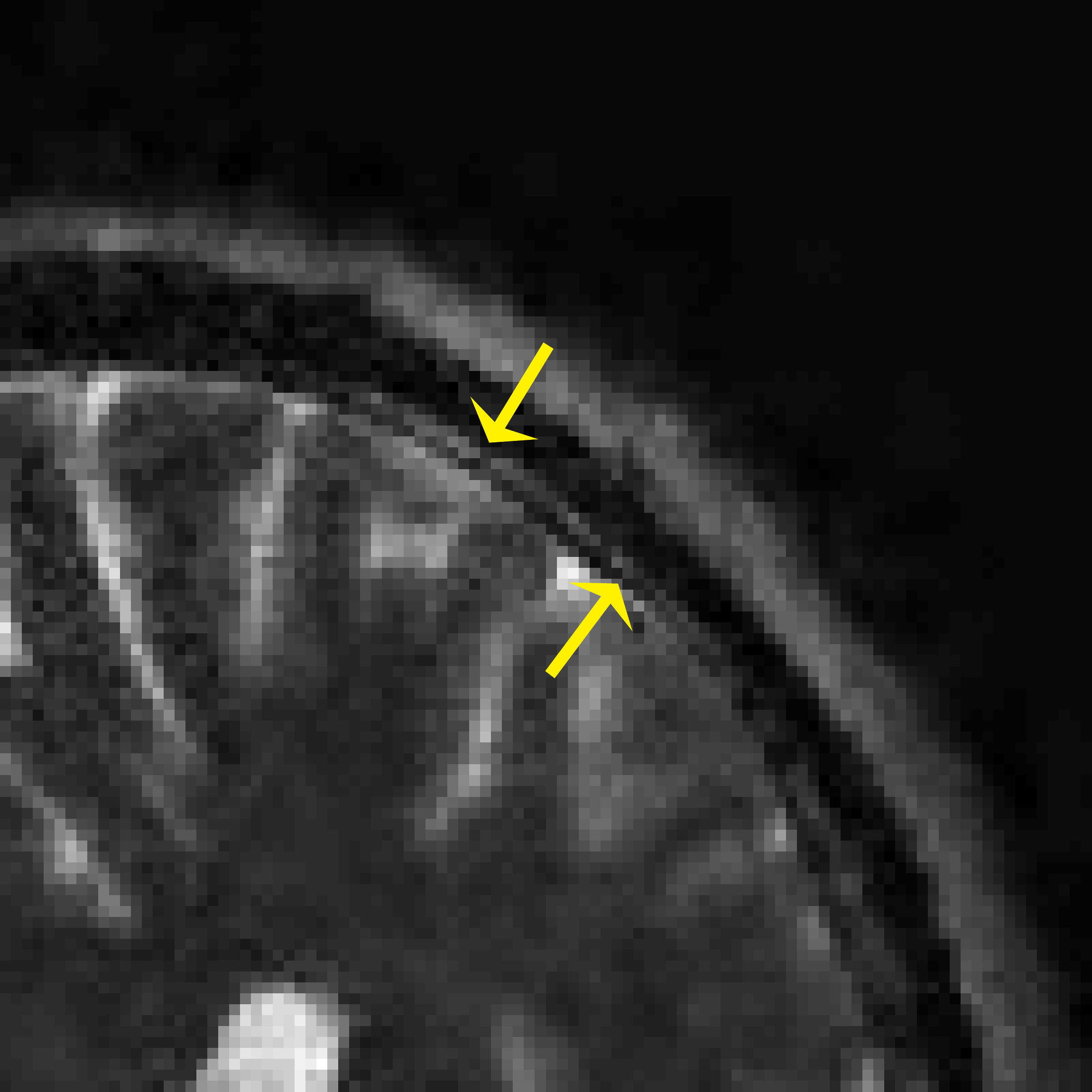}}\hspace{0.005cm}
\subfloat[GIRAF$0.5$]{\label{RealGIRAFHalfZoom}\includegraphics[width=3.00cm]{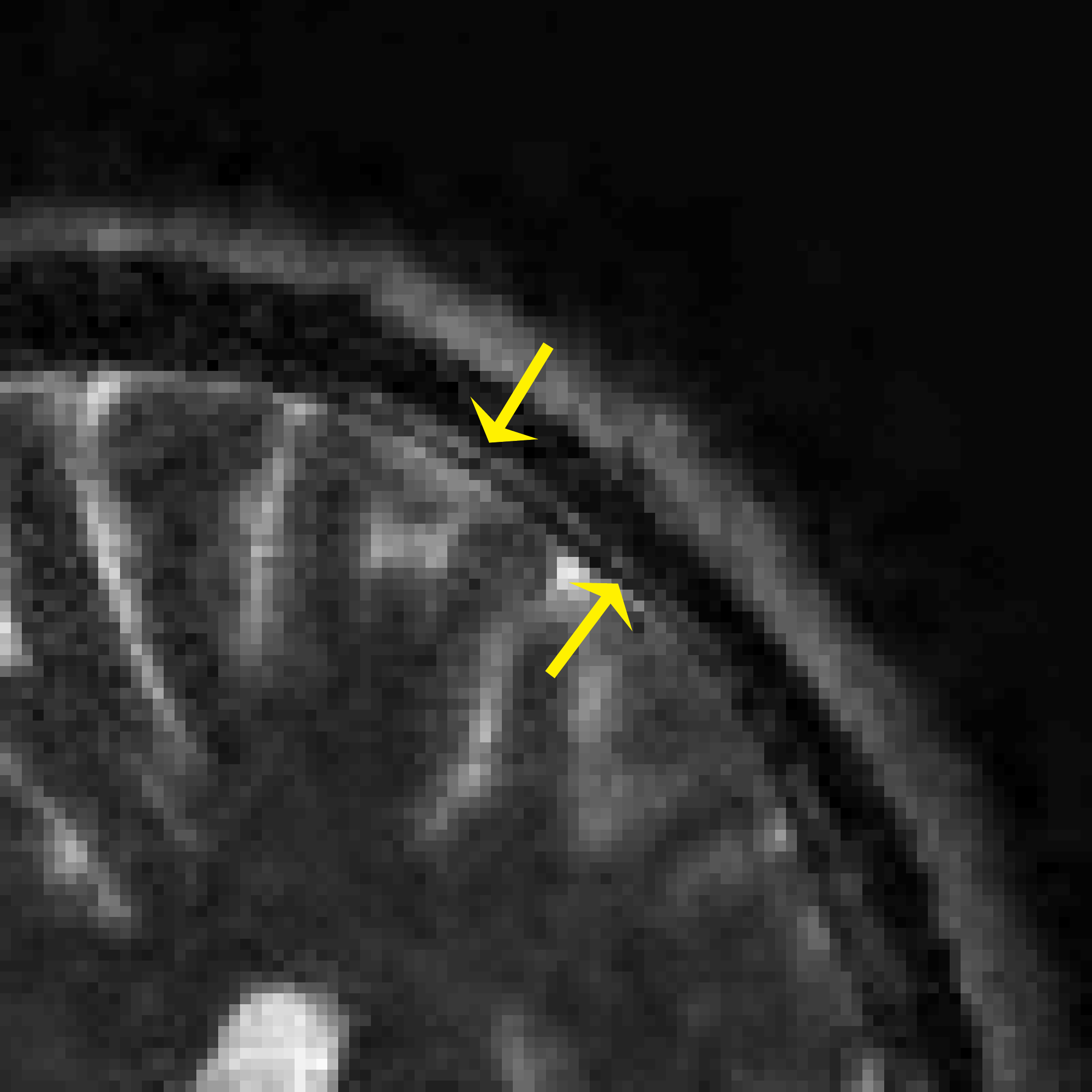}}\hspace{0.005cm}
\subfloat[GIRAF$1$]{\label{RealGIRAF1Zoom}\includegraphics[width=3.00cm]{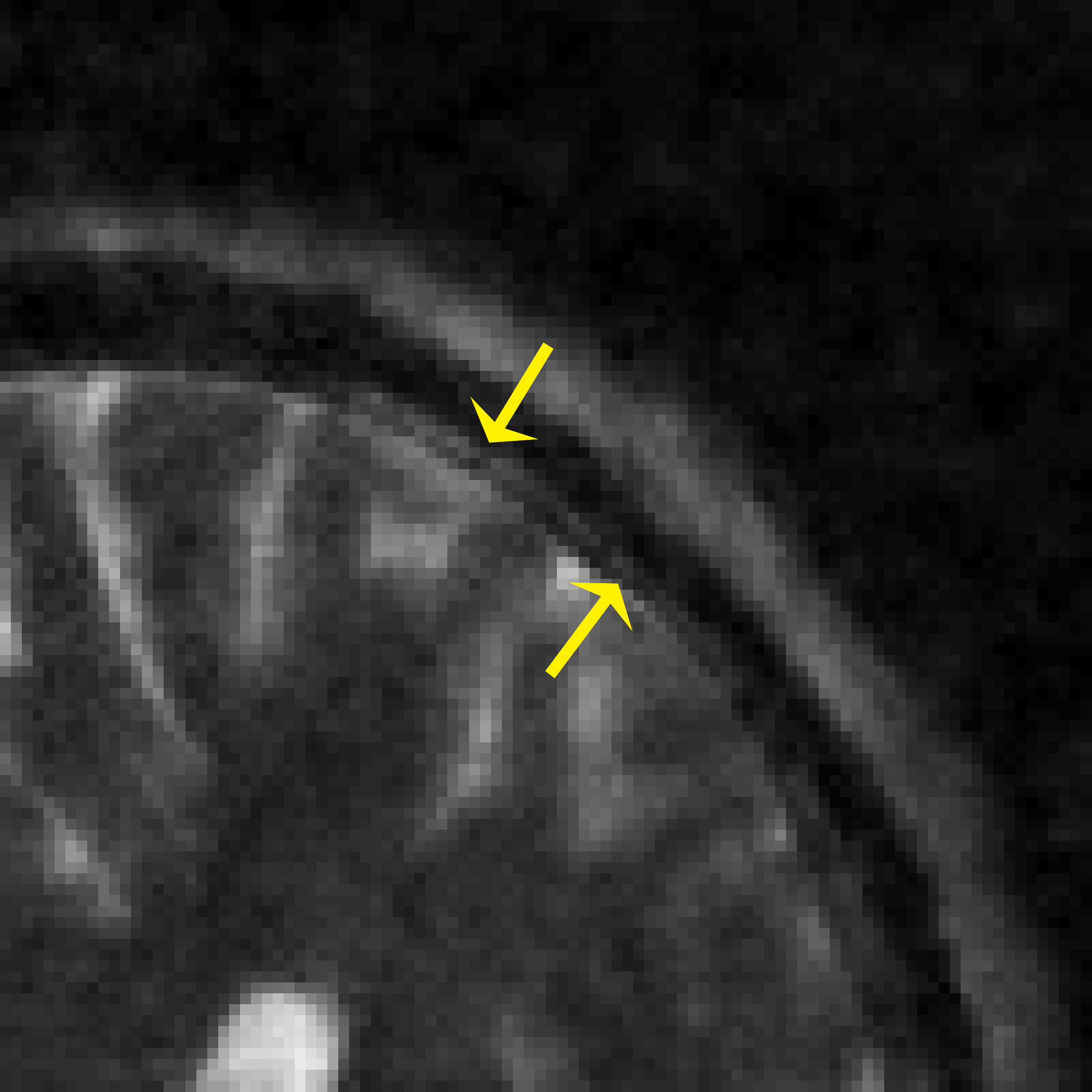}}\hspace{0.005cm}
\subfloat[DDTF \cref{ProposedCSMRIModel}]{\label{RealDDTFZoom}\includegraphics[width=3.00cm]{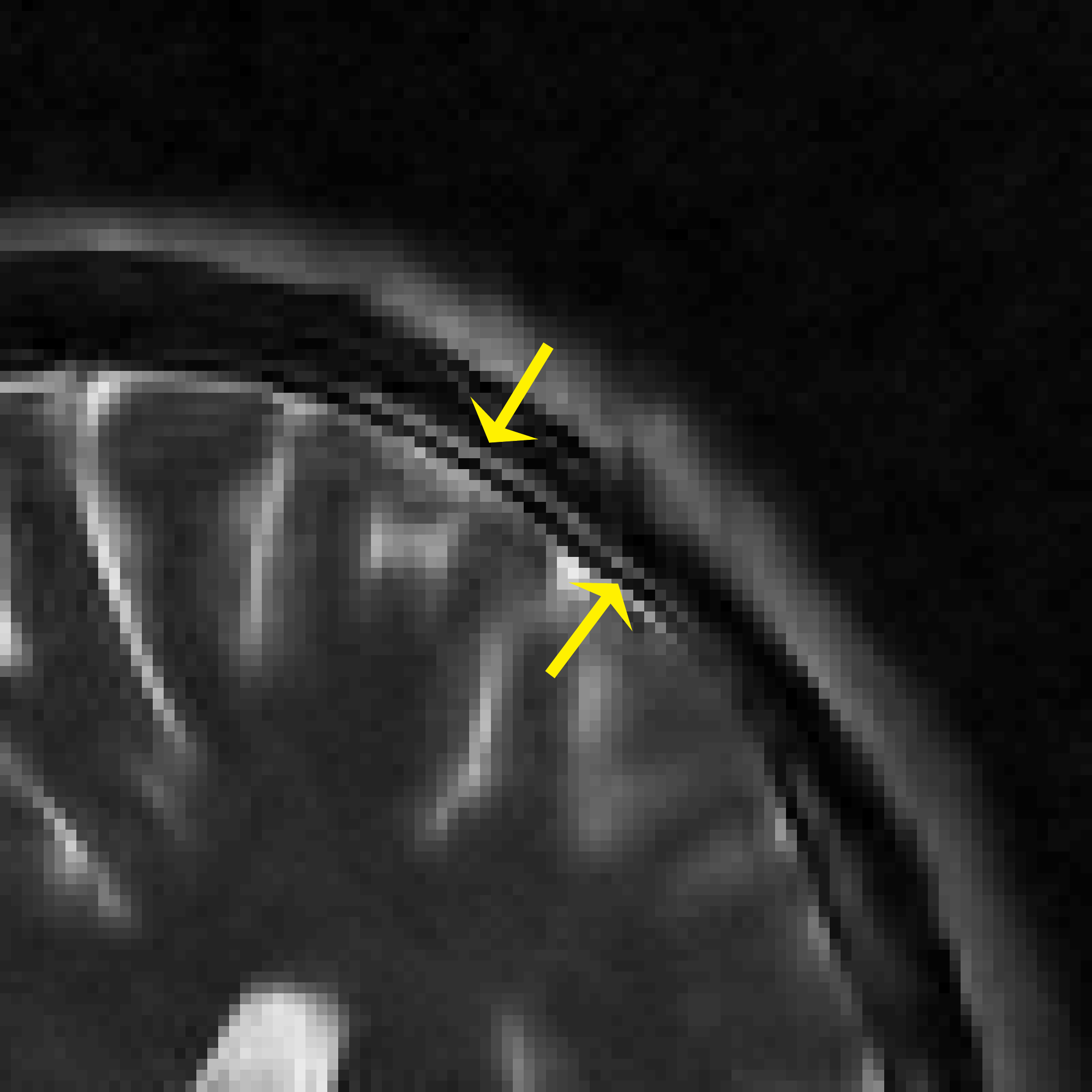}}\vspace{-0.20cm}
\caption{Zoom-in views of \cref{RealResults}. The yellow arrows indicate the regions worth noticing.}\label{RealResultsZoom}
\end{figure}

\cref{RealTable} summarizes the SNR, the HFEN, the number of iterations and the CPU time of the restoration results. For visual comparisons, the restored images, the error maps, the restored k-space data, and the zoom-in views are presented in \cref{RealResults,RealResultsk,RealResultsErrorMap,RealResultsZoom}, respectively. We can observe that our proposed model \cref{ProposedCSMRIModel} outperforms \cref{TVModel,HaarModel,Schatten} as in the phantom experiment, but \cref{TLMRI} performs slightly better than \cref{ProposedCSMRIModel} in the real MR experiment. Nonetheless, the difference of SNR between \cref{ProposedCSMRIModel,TLMRI} is $0.14\mathrm{dB}$, and visually they appear nearly indistinguishable, as shown in \cref{RealTLMRI,RealDDTF,RealTLMRIZoom,RealDDTFZoom}. In addition, even though our approach is relatively time-consuming (approximately $4\mathrm{hrs}$ and $30\mathrm{mins}$) compared to \cref{TVModel,HaarModel,Schatten}, its CPU time is far shorter than \cref{TLMRI} which spends $14\mathrm{hrs}$ but fails to meet stopping criterion. Hence, we can conclude that our proposed approach \cref{ProposedCSMRIModel} performs at least comparably to the TLMRI model \cref{TLMRI} with shorter computational times than \cref{TLMRI}.

Finally, we observe that \cref{ProposedCSMRIModel} for the real MR experiment requires $382$ iterations to meet the stopping criterion, which is more than the piecewise constant phantom image, as shown in \cref{PhantomTable,RealTable}. It should be noted that the original k-space data originally contains the thermal noise. It is also possible that compressing a multi coil k-space data into a virtual single coil data generates an error in the k-space data. In any case, it might fail to meet the underlying assumption of the proposed model; the fully sampled k-space data to be restored is modeled by \cref{MRForward} with the piecewise constant proton density $u$ in \cref{uModel}. In fact, since the proposed data driven tight frame based approach is a relaxation of a low rank (two-fold) Hankel matrix completion, the future work will require the analysis on the convergence rate of \cref{Alg1}, which is likely to depend on the incoherence bound (e.g. \cite{J.F.Cai2019}) of the two-fold Hankel matrix constructed by the original k-space data.

\section{Conclusion and future directions}\label{Conclusion}

In this paper, we propose a new off-the-grid regularization based CS-MRI reconstruction model for the piecewise constant image restoration in the two dimensional FRI framework \cite{G.Ongie2016}\cite{G.Ongie2016}. Our proposed model is inspired by the observation that the SVD of a Hankel matrix to some extent corresponds to an adaptive tight frame system which can sparsely represent a given image. This motivates us to adopt the sparse approximation by the data driven tight frames as an alternative to the structured low rank matrix completion. Finally, the numerical experiments show that our approach outperforms both the conventional on-the-grid approaches and the existing the low rank Hankel matrix models and their relaxations.

For the future work, we plan to provide a rigorous theoretical framework on our observations. Specifically, we need to rigorously analyze whether the data driven tight frame can indeed reflect the true frequency information of spectrally sparse signals (Such signals correspond to low rank Hankel matrices). In addition to the analysis on the data driven tight frame, we also need to analyze the relation between the convergence rate of \cref{Alg1} and the incoherence bound of the underlying original k-space data. It would be also interesting to find other pseudodifferential operators/Fourier integral operators which can provide more insightful information on the rank of the structured matrix constructed from the Fourier transform of a piecewise constant function. For example, we can attempt to find transformations under which the rank of the structured matrix is related to $J$ in \cref{uModel}.

Finally, to broaden the scope of applications, it is also likely to extend the idea in this paper to 1) the piecewise smooth image restoration framework, such as the total generalized variation model \cite{K.Bredies2010} and the combined first and second order TV model \cite{M.Bergounioux2010,K.Papafitsoros2014}, by considering the higher order derivatives; 2) the ``blind'' multi-coil CS-MRI reconstruction task, whose inverse problem corresponds to a multi-channel blind deconvolution problem with the unknown coil sensitivity map. For the blind multi-coil CS-MRI reconstruction, even though it is not clear for us at this moment how to extend our approach, this is nevertheless definitely a future direction we would like to work on. Motivated by the recent works on the deep learning techniques combined with the structured low rank matrix approaches \cite{M.Jacob2019,T.H.Kim2019,A.Pramanik2019a,A.Pramanik2019}, we are also interested in developing a deep learning framework based on our approach.

\appendix
\section{Proof of \cref{Th1}}\label{ProofTh1} Let $u(\x)$ be a piecewise constant function defined as \cref{uModel} and $\bsv=\msF(u)\big|_{L^{-1}\OO}$. Assume that $\rank\left(\bmH\left(\La\bsv\right)\right)=r$, and consider its full SVD
\begin{align}\label{HankelSVD}
\bmH\left(\La\bsv\right)=\left[\begin{array}{c}
\bmH\left(\La_1\bsv\right)\\
\bmH\left(\La_2\bsv\right)
\end{array}\right]=\bX\Sig\bY^*=\sum_{j=1}^{M_2}\Sig^{(j,j)}\bX^{(:,j)}\left(\bY^{(:,j)}\right)^*,
\end{align}
with $\Sig^{(1,1)}\geq\Sig^{(2,2)}\geq\cdots\geq\Sig^{(r,r)}>0$ and $\Sig^{(j,j)}=0$ for $j>r$. Consider $\a_j=M_2^{-1/2}\bY^{(:,j)}$ by reformulating $\bY^{(:,j)}\in\C^{M_2}$ in \cref{HankelSVD} into a $K_1\times K_2$ filter supported on $\KK$. We note that
\begin{align}\label{Stiefel}
\sum_{j=1}^{M_2}\left(M_2^{-1/2}\bY^{(k,j)}\right)\left(M_2^{-1/2}\bY^{(l,j)}\right)^*=M_2^{-1}\delta[k-l],~~~~~k,l=1,\cdots,M_2.
\end{align}
Taking summation along the diagonal and rearranging into the two dimensional multi-indices, we can write \cref{Stiefel} as
\begin{align*}
\sum_{j=1}^{M_2}\sum_{\bsl\in\KK}\a_j[\bk+\bsl]\overline{\a}_j[\bsl]=\dde[\bk].
\end{align*}
Hence, $\bmW$ and $\bmW^*$ in \cref{OurAnalysis,OurSynthesis} defined by using these $\a_1,\cdots,\a_{M_2}$ satisfies
\begin{align*}
\bmW^*\bmW\left(\La\bsv\right)=\sum_{j=1}^{M_2}\bmS_{\overline{\a}_j}\left(\bmS_{\a_j[-\cdot]}\left(\La\bsv\right)\right)=\La\bsv,
\end{align*}
which shows that the filters $\a_1,\cdots,\a_{M_2}$ form a tight frame system.

To complete the proof, we note that
\begin{align*}
\Sig^{(j,j)}\bX^{(:,j)}=\bmH\left(\La\bsv\right)\bY^{(:,j)}=\left(\bY^{(:,j)}[-\cdot]\left(\La\bsv\right)\right)\big|_{\OO:\KK},
\end{align*}
where $\ast$ is performed by reformulating $\bY^{(:,j)}\in\C^{M_2}$ into a $K_1\times K_2$ filter. Then \cref{HankelSVD} implies that for $j=r+1,\cdots,M_2,$ we have
\begin{align*}
\left(\bmS_{\a_j[-\cdot]}\left(\La\bsv\right)\right)[\bk]=\0,~~~~~~~\bk\in\OO:\KK.
\end{align*}
This completes the proof.

\section*{Acknowledgments} The authors would like to thank Dr. Greg Ongie in the Department of Statistics at the University of Chicago, an author of \cite{G.Ongie2018,G.Ongie2015,G.Ongie2016,G.Ongie2017}, for making the data sets as well as the MATLAB toolbox available so that the experiments can be implemented. The authors also thank the anonymous reviewers for their constructive suggestions and comments that helped tremendously with improving the presentation of this paper.

\bibliographystyle{siamplain}

\end{document}